\documentclass[10pt,letterpaper,reqno]{amsart}
\usepackage{color}
\usepackage[latin1]{inputenc}
\usepackage{amsmath}
\usepackage{amsfonts}
\usepackage{amsthm}
\usepackage{epsfig}
\usepackage{graphicx}
\usepackage{amssymb}
\usepackage{bm}
\usepackage{esint}
\usepackage{caption}
\usepackage{enumerate}
\usepackage{enumitem}
\usepackage{esvect}
\usepackage{nccmath}
\usepackage{mathrsfs}
\usepackage{verbatim}
\usepackage{url}

\newcommand{\ssubset}{\subset\joinrel\subset}
\usepackage{scalerel,stackengine}
\newcommand\pig[1]{\scalerel*[5pt]{\big#1}{%
		\ensurestackMath{\addstackgap[1.5pt]{\big#1}}}}
\newcommand\pigl[1]{\mathopen{\pig{#1}}}

\newcommand{\bigiotab}{\makebox{\large\ensuremath{\iota}}}

\usepackage[a4paper,bindingoffset=0.2in,%
left=1.2in,right=1.2in,top=1.2in,bottom=1.2in,%
footskip=.35in]{geometry}

\newcommand{\spt}{{\rm spt}}

\DeclareMathOperator{\restrict}{\llcorner}

\newcommand{\Der}{D}

\newcommand{\der}{\mathbf{D}}

\theoremstyle{plain}

\newtheorem{theorem}{Theorem}[section]
\newtheorem{lemma}[theorem]{Lemma}
\newtheorem{corollary}[theorem]{Corollary}

\newtheorem*{theorem*}{Theorem}
\newtheorem*{corollary*}{Corollary}

\theoremstyle{definition}
\newtheorem{definition}[theorem]{Definition}

\newtheorem{remark}[theorem]{Remark}

\newtheorem*{notation*}{Notation}

\numberwithin{equation}{section}
\numberwithin{figure}{section}

\DeclareRobustCommand{\rchi}{{\mathpalette\irchi\relax}}
\newcommand{\irchi}[2]{\raisebox{\depth}{$#1\chi$}}

\usepackage{xcolor}

\title{Legendrian cycles and Reilly-type variational formulae 
	for $F_nW^{2,n}$-sets}

\author{Paolo Valentini}
\address{Dipartimento di Ingegneria e Scienze dell'Informazione e Matematica, Universit\'a degli Studi dell'Aquila, 67100 L'Aquila, Italy}
\email{paolo.valentini@graduate.univaq.it}

\begin{document}

	\begin{abstract} We construct a natural Legendrian cycle $\mathcal{N}_\mathcal{S}$ associated with any $F_nW^{2,n}$-set $\mathcal{S}$, that is, a closed set locally described as a finite union of graphs of $(C^0\cap W^{2,n})$-regular functions with integer multiplicity. The construction relies on the fact that $\mathcal{S}$ is countably $\mathcal{H}^n$-rectifiable of class $C^2$ and, at $\mathcal{H}^n$-almost every point $p\in\mathcal{S}$, the proximal unit normal bundle at $p$, denoted by $\operatorname{nor}(\mathcal{S},p)$, consists of exactly two antipodal vectors $\{u,-u\}$, even in the presence of overlapping $W^{2,n}$-graphs. As a consequence, we prove Reilly-type variational formulae for the higher-order mean curvature integrals of $\mathcal{S}$, extending the classical results of Reilly to this non-smooth setting.
	\end{abstract}

	\maketitle
	\tableofcontents
	\paragraph*{\small MSC-classes 2020:}{\small 53C65, 49Q15, 49Q20.}
	
	\section{Introduction}
	
	\subsection{Background and motivation}
	
	One of the most fertile lines of inquiry in modern geometric analysis is to 
	understand which classical results of differential geometry survive when the 
	smoothness assumptions on the underlying objects are weakened. This question 
	is not merely of technical interest: it probes the very geometric content of 
	the results themselves, separating those that are intrinsically metric or 
	measure-theoretic in nature from those that depend on the availability of 
	smooth machinery. At the same time, it opens the door to genuinely new phenomena 
	that have no smooth counterpart.
	
	The present paper pursues this philosophy in the context of 
	\emph{higher-order mean curvature integrals} and their \emph{variational 
		properties}. Our objects of study are the so-called \emph{$F_n W^{2,n}$-sets}, 
	a class of geometric objects introduced by Ambrosio, Gobbino and Pallara \cite{AmGobPal} -- following an original idea of De Giorgi -- as local sums of 
	characteristic functions of graphs of $W^{2,p}$-regular functions of $n$ variables with $p>n$, 
	with multiplicities. In the present paper we work in the critical exponent $p = n$, 
	which is the lowest integrability threshold at which the curvature theory of 
	the individual graphs is well-defined. Precisely, a closed set $\mathcal{S} \subset \mathbf{R}^{n+1}$ 
	is a $F_nW^{2,n}$-set if for every point $z \in \mathcal{S}$ there exist a positive integer 
	$q(z)$, an open neighbourhood $U$ of $z$, and a family of hypersurfaces 
	$\Gamma_1, \ldots, \Gamma_{q(z)}$ such that
	\[
	\mathcal{S} \cap U = \bigcup_{i=1}^{q(z)} (\Gamma_i \cap U)\,,
	\]
	where each $\Gamma_i \cap U$ is (up to a $C^2$-diffeomorphism $F$ of $\mathbf{R}^{n+1}$) the graph of 
	a $(C^0 \cap W^{2,n})$-regular function. The integer-valued function 
	$\pmb{\iota} = \sum_{i=1}^{q(z)} \mathbf{1}_{\Gamma_i}$ recording the local 
	multiplicity of $\mathcal{S}$ plays a central role throughout.
	
	The choice of the Sobolev exponent $p = n$ is critical in a precise sense. On the one hand, for every $f \in W^{2,n}_{\textrm{loc}}(\mathbf{R}^n)$ the set $\mathcal{S}(f)$ of points 
	at which $f$ is pointwise twice differentiable has full $\mathcal{L}^n$-measure 
	(more generally, Calder\'on and Zygmund proved that functions in 
	$W^{2,p}_{\textrm{loc}}(\mathbf{R}^n)$ with $p > \frac{n}{2}$ admit a 
	second-order Taylor expansion $\mathcal{L}^n$-a.e.; 
	cf.~\cite{CalderonZygmund} or \cite[Proposition~2.2]{Caffarelli96}). 
	Moreover, $\nabla f$ satisfies the Lusin 
	$(N)$-condition on $\mathcal{S}(f)$ (this follows from Lemma~\ref{W2n functions Lusin}). On the other hand, $W^{2,n}$-functions 
	can exhibit genuinely singular behaviour that is invisible to the classical smooth setting. The Lusin $(N)$-property does not generally hold for the gradient of $W^{2,n}$-functions: Tom\'as Roskovec 
	(cf.\,\cite{Roskovec}), using a Cesari-type construction, provides an 
	example of a function $f \in C^1([-1,1]^n)$ such that $\nabla f$ is a 
	$(C^0 \cap W^{1,n})$-regular vector field that does not satisfy the Lusin 
	$(N)$-property. Furthermore, Tatiana Toro \cite{Toro} constructs a 
	$W^{2,n}$-function with a countable dense set of singular points, and Joseph Fu 
	\cite[p.\,2260]{FuAlexandrov} observes that the gradient of a 
	$W^{2,n}$-function may have a graph dense in $\mathbf{R}^n \times 
	\mathbf{R}^n$. Thus, $F_nW^{2,n}$-sets sit precisely at the boundary 
	between classical regularity and genuinely non-smooth behaviour, making 
	them both challenging and geometrically rich.
	
	\subsection{The Legendrian cycle associated with a $F_nW^{2,n}$-set}
	
	The key to unlocking the variational geometry of $F_nW^{2,n}$-sets lies in 
	the theory of \emph{Legendrian cycles}, a class of integer-multiplicity locally 
	rectifiable currents developed by Fu \cite{Fu98} as the natural framework for 
	Lipschitz--Killing curvature integrals and their variational properties. 
	Recall that an integer-multiplicity locally rectifiable $n$-current $T$ 
	of $\mathbf{R}^{n+1} \times \mathbf{R}^{n+1}$ with support in 
	$\mathbf{R}^{n+1} \times \mathbf{S}^n$ is called a \emph{Legendrian cycle} if 
	$\partial T = 0$ and $T \llcorner \alpha = 0$, where $\alpha$ is the 
	contact $1$-form on $\mathbf{R}^{n+1} \times \mathbf{R}^{n+1}$.
	
	In our preceding work \cite{SantilliValentini}, we proved that the proximal 
	unit normal bundle $\operatorname{nor}(\Gamma_f)$ of 
	$\Gamma_f = \mathrm{graph}(f)$, for an arbitrary $W^{2,n}$-function $f$, 
	carries a canonical structure of \emph{Legendrian cycle} 
	(cf.\ \cite[Theorem~3.9]{SantilliValentini}). This result -- whose proof 
	rests on fine properties of $W^{2,n}$-functions -- was then used to 
	extend Reilly's variational formulae and Alexandrov's sphere theorem to 
	hypersurfaces locally described as $W^{2,n}$-graphs.
	
	The present paper addresses the highly delicate problem of 
	constructing a Legendrian cycle for the entire class of $F_nW^{2,n}$-sets, 
	i.e.\ for objects that are locally a \emph{finite union of overlapping} 
	$W^{2,n}$-regular graphs, each carrying its own multiplicity. Two fundamental properties -- both established in Lemma \ref{Ufset14}, which constitutes the technical core of the paper -- are the following. First,
	\[
	\mathcal{H}^n\bigl(\mathcal{S}\setminus N_2(\mathcal{S})\bigr)=0\,,
	\]
	that is, for $\mathcal{H}^n$-almost every $p\in\mathcal{S}$, the proximal unit normal bundle $\operatorname{nor}(\mathcal{S},p)$ at $p$ consists of exactly two antipodal vectors $\{u,-u\}$, even in the presence of overlapping $W^{2,n}$-graphs at the point $p$. Second, $\mathcal{S}$ is countably $\mathcal{H}^n$-rectifiable of order $2$. As a consequence (cf.\,Lemma \ref{LemmaSelection}), the multivalued map
	$$
	\textup{nor}(\mathcal{S}, \pmb{\cdot}) : \mathcal{S} \to \mathcal{P}(\mathbf{S}^n)
	$$
	admits an $\mathcal{H}^n$-measurable selection 
	$\nu_\mathcal{S} : \mathcal{S} \to \mathbf{S}^n$, such that
	$$
	\nu_\mathcal{S}(p) \in \textup{Nor}^n(\mathcal{H}^n \restrict \mathcal{S}, p) 
	\quad \text{for $\mathcal{H}^n$-a.e.\ $p \in \mathcal{S}$}\,.
	$$
	Furthermore, by virtue of the second-order rectifiability of $\mathcal{S}$, we deduce that $\nu_{\mathcal{S}}$ is $(\mathcal{H}^n \restrict \mathcal{S})$-approximately differentiable $\mathcal{H}^n$-a.e. with symmetric approximate differential (cf. Lemma~\ref{lem approx diff unit normal}). This forms the foundation of the entire curvature theory.
	
	With this structural result in hand, we construct in Theorem \ref{ReillyWset1} the 
	\emph{Legendrian cycle $\mathcal{N}_{\mathcal{S}}$ associated with $\mathcal{S}$} as the integer-multiplicity 
	current
	$$\mathcal{N}_{\mathcal{S}}=\big(\pmb\bigiotab\circ\pi_0\circ(\Psi_F|\textup{nor}(\mathcal{S}'))^{-1}\big)\,\big(\mathcal{H}^n\restrict\textup{nor}(\mathcal{S})\big)\wedge\vv{\xi}_{\mathcal{S}}\,.$$
	The same theorem yields the key representation formula
	{\allowdisplaybreaks\begin{align*}
			&\big(\mathcal{N}_{\mathcal{S}} \restrict \varphi_{n-k}\big)(\phi)&\\ \nonumber
			&\quad\quad\quad={n \choose k}\int_{\mathcal{S}} \big[\phi\big(x, \nu_{\mathcal{S}}(x)\big) +(-1)^k\,\phi\big(x,-\nu_{\mathcal{S}}(x)\big)\big]H_{\mathcal{S}, k}(x)\,\pmb\bigiotab\big(F^{-1}(x)\big)\, d\mathcal{H}^n(x)
	\end{align*}}for every $ \phi \in C^\infty(\mathbf{R}^{n+1} \times \mathbf{R}^{n+1}) $ and $ k\in\{0, \ldots, n \}$, which expresses the action of $\mathcal{N}_{\mathcal{S}}$ on the \emph{Lipschitz-Killing forms 
	$\varphi_{n-k}$} in terms of the $k$-th mean curvature $H_{\mathcal{S},k}$ 
	and the multiplicity $\pmb{\iota}$. This formula is the key to the variational results established 
	in the next section.
	
	\subsection{Reilly-type variational formulae} The construction of $\mathcal{N}_\mathcal{S}$ is not an end in itself: 
	it is the engine that drives the main results of the paper. We prove 
	the following extension of Reilly's variational formulae \cite{Reilly1972} 
	to $F_nW^{2,n}$-sets.
	
	\begin{theorem*}[cf.\,Theorem \ref{ReillyW-sets}]
		Let $\mathcal{S}$ be a compact $F_nW^{2,n}$-set with multiplicity function $\pmb{\iota}$, 
		let $\nu_\mathcal{S}$ be a selected unit-normal vector field on $\mathcal{S}$, and let $ \{F_t\}_{t \in (-\epsilon, \epsilon)} $ be a local variation of $ \mathbf{R}^{n+1} $ with initial velocity vector field $ V $. If $k\in\{1, \ldots , n\}$ is odd, then
		$$\frac{d}{dt}\,\mathcal{A}_{k-1}\big(F_t(\mathcal{S})\big) \Big|_{t =0}=(n-k+1)\int_{\mathcal{S}}  V(x)\bullet\nu_\mathcal{S}(x) \,H_{\mathcal{S}, k}(x)\,\pmb\bigiotab\big(F^{-1}(x)\big)\, d\mathcal{H}^n(x)\,.$$Moreover, if $n$ is even,
		$$\frac{d}{dt}\,\mathcal{A}_n\big(F_t(\mathcal{S})\big) \Big|_{t =0}=0 \, .$$
	\end{theorem*}
	\noindent Here $H_{\mathcal{S},k}$ denotes the $k$-th mean curvature of $\mathcal{S}$ and $\mathcal{A}_{k-1}(\mathcal{S}) = \int_{\mathcal{S}} H_{\mathcal{S},k-1}\, 
	\pmb{\iota}\circ F^{-1}\, d\mathcal{H}^n$
	is the $(k-1)$-th total curvature weighted by multiplicity (cf.\,Definitions \ref{rth mean curvature function3} and \ref{A_k}). 
	
	For $k=1$ the formula reads
	$$\frac{d}{dt}\,\mathcal{A}_{0}\big(F_t(\mathcal{S})\big) \Big|_{t =0}=n\int_{\mathcal{S}}  V(x)\bullet\nu_\mathcal{S}(x) \,H_{\mathcal{S}, 1}(x)\,\pmb\bigiotab\big(F^{-1}(x)\big)\, d\mathcal{H}^n(x)\,,$$
	which is precisely the \emph{first variation of the weighted area 
		functional} $\mathcal{A}_{0}(\mathcal{S}) = \int_{\mathcal{S}} 
	\pmb{\iota}\circ F^{-1}\, d\mathcal{H}^n$. 
	This extends the classical first variation formula for smooth 
	hypersurfaces to $F_nW^{2,n}$-sets with multiplicity, which, 
	as shown by Tatiana Toro \cite{Toro}, may exhibit no classical 
	tangent plane at a countable dense set of points.
	
	In particular, when $\mathcal{S}$ is a single $W^{2,n}$-graph, 
	the weight $\pmb{\iota}\circ F^{-1}$ reduces to $1$ and the 
	formulae of \cite{SantilliValentini} are recovered as a special case.
	
	\subsection{Organisation of the paper}
	
	The paper is organised as follows. In Section~2 we collect the 
	necessary background from geometric measure theory, the theory of 
	Legendrian cycles, the proximal unit normal bundle, and the fine 
	properties of $W^{2,n}$-functions that are used throughout. 
	Section~3 contains the main results: in 3.1 we introduce 
	$F_n W^{2,n}$-sets and establish their fine structural 
	properties (Lemma \ref{Ufset14}); in 3.2 we define the approximate principal 
	curvatures on $F_nW^{2,n}$-sets; in 3.3 we construct the associated 
	Legendrian cycle (Theorems \ref{LegCycSset} and \ref{ReillyWset1}) and prove the Reilly-type 
	variational formulae (Theorem \ref{ReillyW-sets}).
	
	Some of the results of this paper are contained in the author's doctoral thesis \cite{Val25}.

	\section{Notation and background}\label{section: Notation and background}
	Given a set of parameters $\{p_1, \ldots, p_n\}$, we denote by $\textsf{c}(p_1, \ldots, p_n)$ 
	a generic positive constant depending only on $p_1, \ldots, p_n$.
	
	For a function $f : S \rightarrow T$, we define the \emph{graph map} $\overline{f} : S \rightarrow S \times T$ by
	\begin{equation}\label{graph map}
		\overline{f}(x): = \big(x,f(x)\big)\quad\text{for every }x\in S\,.
	\end{equation}
	We also make use of the projection maps
	\begin{equation}\label{projections}
		\pi_0 : \mathbf{R}^{n+1} \times \mathbf{R}^{n+1} \rightarrow \mathbf{R}^{n+1}\quad\text{and}\quad
		\pi_1 : \mathbf{R}^{n+1} \times \mathbf{R}^{n+1} \rightarrow \mathbf{R}^{n+1},
	\end{equation}
	defined by $\pi_0(x,u) := x$ and $\pi_1(x,u) := u$, whenever $x,u\in\mathbf{R}^{n+1}$.
	
	Throughout this paper, the symbol $\bullet$ denotes the scalar product. 
	For a subset $S$ of a Euclidean space, $\overline{S}$ denotes its closure. 
	The symbols $\Der$ and $\nabla$ denote the classical differential and gradient, respectively. 
	We use $\bm{\nabla}f$ and $\der f$ to denote the weak gradient and the approximate 
	differential of a Sobolev function $f$ $($\textup{cf.}\,\cite[Definition 3.70]{AFP00}$)$. 
	The characteristic function of a set $X$ is denoted by $\bm{1}_X$. 
	The Grassmannian of $m$-dimensional subspaces of $\mathbf{R}^k$ is $\mathbf{G}(k,m)$. 
	For a continuous function $f : U \rightarrow \mathbf{R}$ defined on an open set $U$, 
	we denote by $\textsf{Diff}(f)$ the set of points $x \in U$ at which $f$ is pointwise differentiable.
	
	\subsection{Basic notions from geometric measure theory}\addcontentsline{toc}{subsection}{\protect\numberline{}Basic notions from geometric measure theory}
	In this paper we use standard notation from geometric measure theory, 
	for which we refer to~\cite{Fed69}. For the reader's convenience, 
	we recall some basic notions here.
	
	Given $X \subset \mathbf{R}^m$ and $a \in \mathbf{R}^m$, the 
	\emph{tangent cone} of $X$ at $a$, denoted $\textup{Tan}(X,a)$, 
	is the set of all $v \in \mathbf{R}^m$ such that there exists a sequence 
	$\{a_k\}_{k \in \mathbf{N}} \subset X \setminus \{a\}$ with
	$$
	\lim_{k \to \infty} a_k = a 
	\qquad \text{and} \qquad 
	\lim_{k \to \infty} \frac{a_k - a}{|a_k - a|} = \frac{v}{|v|}\,.
	$$
	Equivalently, $\textup{Tan}(X,a)$ is the set of all $v \in \mathbf{R}^m$ 
	such that for every $\epsilon > 0$ there exist $x \in X$ and $r > 0$ 
	with $|x - a| < r$ and $|r(x-a) - v| < \epsilon$. The \emph{normal cone} 
	of $X$ at $a$, with vertex at $0$, is
	$$
	\textup{Nor}(X,a) := \big\{ u \in \mathbf{R}^m : u \bullet v \leq 0 
	\;\text{for all}\; v \in \textup{Tan}(X,a) \big\}\,.
	$$
	
	Given $X \subset \mathbf{R}^m$, $a \in \mathbf{R}^m$, and a positive 
	integer $\mu$, the \emph{$(\mathcal{H}^\mu \restrict X)$-approximate 
		tangent cone} of $X$ at $a$, denoted 
	$\textup{Tan}^\mu(\mathcal{H}^\mu \restrict X, a)$, is the set of all 
	$v \in \mathbf{R}^m$ such that
	\begin{equation}\label{densityTang}
		\Theta^{*\mu}\big(\mathcal{H}^\mu \restrict X \cap 
		\{x : |r(x-a) - v| < \epsilon \;\text{for some}\; r > 0\}, a\big) > 0
	\end{equation}
	for every $\epsilon > 0$. The \emph{$(\mathcal{H}^\mu \restrict X)$-approximate 
		normal cone} of $X$ at $a$, with vertex at $0$, is
	$$
	\textup{Nor}^\mu(\mathcal{H}^\mu \restrict X, a) := \big\{ u \in \mathbf{R}^m : 
	u \bullet v \leq 0 \;\text{for all}\; v \in 
	\textup{Tan}^\mu(\mathcal{H}^\mu \restrict X, a) \big\}\,.
	$$
	
	\begin{remark}\label{propertyTan}
		\textup{%
			\emph{(i)} An equivalent characterization of 
			$\textup{Tan}^\mu(\mathcal{H}^\mu \restrict X, a)$ 
			is given by~\cite[$3.2.16$]{Fed69}
			$$
			\textup{Tan}^\mu(\mathcal{H}^\mu \restrict X, a) = 
			\bigcap\big\{\textup{Tan}(E,a) : E \subseteq X,\, 
			\Theta^\mu(\mathcal{H}^\mu \restrict X \setminus E, a) = 0\big\}\,.
			$$
			In particular, $\textup{Tan}^\mu(\mathcal{H}^\mu \restrict X, a)$ 
			is a subcone of $\textup{Tan}(X,a)$.\\[4pt]
			\emph{(ii)} (\emph{Locality of approximate tangent spaces}).
			Let $X$ and $Y$ be $\mathcal{H}^\mu$-measurable subsets of $\mathbf{R}^m$ 
			with finite $\mathcal{H}^\mu$-measure. Then
			\begin{equation}\label{locality}
				\textup{Tan}^\mu(\mathcal{H}^\mu \restrict X, a) = 
				\textup{Tan}^\mu(\mathcal{H}^\mu \restrict Y, a) 
				\quad \text{for $\mathcal{H}^\mu$-a.e.\ $a \in X \cap Y$}.
		\end{equation}}
	\end{remark}
	
	Let $X \subset \mathbf{R}^m$, let $f$ map a subset of $\mathbf{R}^m$ 
	into $\mathbf{R}^k$, let $\mu$ be a positive integer, and let 
	$a \in \mathbf{R}^m$. We say that $f$ is 
	\emph{$(\mathcal{H}^\mu \restrict X)$-approximately differentiable} at $a$ 
	(cf.~\cite[$3.2.16$]{Fed69}) if there exists a map 
	$g : \mathbf{R}^m \rightarrow \mathbf{R}^k$, pointwise differentiable at $a$, 
	such that $f(a) = g(a)$ whenever $a \in \textup{dmn}(f)$ and
	$$
	\Theta^\mu\big(\mathcal{H}^\mu \restrict X \cap 
	\{b : f(b) \neq g(b)\}, a\big) = 0\,.
	$$
	In this case (see~\cite[$3.2.16$]{Fed69}), $f$ determines the restriction 
	of $Dg(a)$ to the approximate tangent cone 
	$\textup{Tan}^\mu(\mathcal{H}^\mu \restrict X, a)$, and we define
	$$
	\textup{ap}\, Df(a) := Dg(a)\big|\textup{Tan}^\mu(\mathcal{H}^\mu \restrict X, a)\,.
	$$
	
	Let $X \subset \mathbf{R}^m$ and let $\mu$ be a positive integer. 
	We say that $X$ is \emph{countably $\mathcal{H}^\mu$-rectifiable} if there 
	exist countably many $\mu$-dimensional $C^1$-submanifolds $\Sigma_i$ of 
	$\mathbf{R}^m$ such that
	\begin{equation}\label{rectCondition}
		\mathcal{H}^\mu\Big(X \setminus \bigcup_{i=1}^\infty \Sigma_i\Big) = 0\,.
	\end{equation}
	We say that $X$ is \emph{locally $\mathcal{H}^\mu$-rectifiable} if~\eqref{rectCondition} 
	holds and $\mathcal{H}^\mu(X \cap K) < \infty$ for every compact 
	$K \subset \mathbf{R}^m$, and simply \emph{$\mathcal{H}^\mu$-rectifiable} 
	if~\eqref{rectCondition} holds and $\mathcal{H}^\mu(X) < \infty$. 
	Finally, $X$ is \emph{countably $\mathcal{H}^\mu$-rectifiable of class $k$} 
	if~\eqref{rectCondition} holds with the $\Sigma_i$ being 
	$\mu$-dimensional $C^k$-submanifolds.
	
	It is well known that if $X$ is countably $\mathcal{H}^\mu$-rectifiable 
	with $\mathcal{H}^\mu(X) < \infty$, then $\textup{Tan}^\mu(\mathcal{H}^\mu 
	\restrict X, a)$ is a $\mu$-dimensional plane at $\mathcal{H}^\mu$-a.e.\ 
	$a \in X$, and every Lipschitz function $f : X \rightarrow \mathbf{R}^k$ 
	admits an $(\mathcal{H}^\mu \restrict X)$-approximate differential
	$$
	\textup{ap}\, Df(a) : \textup{Tan}^\mu(\mathcal{H}^\mu \restrict X, a) 
	\rightarrow \mathbf{R}^k
	$$
	at $\mathcal{H}^\mu$-a.e.\ $a \in X$. At such points, for each 
	$h \in \{1, \ldots, k\}$, the \emph{$(\mathcal{H}^h \restrict X)$-approximate 
		Jacobian} of $f$ at $a$ is defined as
	$$
	J_h^X f(a) := \sup\big\{ \big|[{\textstyle\bigwedge_h \textup{ap}\,Df(a)}](\xi)\big| 
	: \xi \in {\textstyle\bigwedge_h \textup{Tan}^\mu(\mathcal{H}^\mu \restrict X, a)},\, 
	|\xi| = 1 \big\}
	$$
	(see~\eqref{natural map} for the definition of 
	${\textstyle\bigwedge_h \textup{ap}\,Df(a)}$). The approximate Jacobian 
	appears naturally in the area and coarea formulas for $f$ (cf.~\cite[$3.2.20,\,3.2.22$]{Fed69}).
	
	We recall the following result from~\cite[Lemma~5.2]{SantilliValentini}.
	
	\begin{lemma}\label{lem approx diff unit normal}
		Let $X \subset \mathbf{R}^{n+1}$ be $\mathcal{H}^n$-measurable and 
		countably $\mathcal{H}^n$-rectifiable of class $C^2$, and let 
		$\nu : X \to \mathbf{S}^n$ be an $(\mathcal{H}^n \restrict X)$-measurable 
		map such that
		\begin{equation}\label{hypotesis}
			\nu(a) \in \textup{Nor}^n(\mathcal{H}^n \restrict X, a) 
			\quad \text{for $\mathcal{H}^n$-a.e.\ $a \in X$.}
		\end{equation}
		Then there exists a countable family of measurable sets 
		$X_i \subseteq X$ such that 
		$\mathcal{H}^n\big(X \setminus \bigcup_{i=1}^\infty X_i\big) = 0$ 
		and $\textup{Lip}(\nu|X_i) < \infty$. Moreover, $\nu$ is 
		$(\mathcal{H}^n \restrict X)$-approximately differentiable at 
		$\mathcal{H}^n$-a.e.\ $a \in X$, and $\textup{ap}\,D\nu(a)$ is a 
		symmetric endomorphism of $\textup{Tan}^n(\mathcal{H}^n \restrict X, a)$.
	\end{lemma}
	
	\subsection{Differential forms and currents} Given a finite-dimensional vector space $ V $ over $\mathbf{R}$, we denote by $ v_1 \wedge \cdots \wedge v_k $ the \textit{simple $ k $-vector} obtained by the \emph{exterior product} of vectors $ v_1, \ldots , v_k $ in $ V $, and by $ \bigwedge_k V $ the vector space generated by all simple $ k $-vectors of $ V $. Every linear map $ L : V \rightarrow V' $ uniquely extends to a linear map 
	\begin{equation}\label{natural map}
		{\textstyle \bigwedge_k} L : {\textstyle \bigwedge_k} V \rightarrow {\textstyle \bigwedge_k} V' 
	\end{equation}
	such that $ [{\textstyle \bigwedge_k} L](v_1 \wedge \cdots \wedge v_k) := L(v_1) \wedge \cdots \wedge L(v_k) $ for every $ v_1, \ldots , v_k \in V $.
	
	The vector space of all \emph{alternating $ k $-linear} maps $ \Phi : V^k \rightarrow \mathbf{R} $ (i.e., so that $ \Phi(v_1, \ldots, v_k) =0 $ whenever $ v_1, \ldots, v_k \in V $ and $ v_i = v_j $ for some $ i \neq j $) is denoted by $\bigwedge^k V $.\,The elements of $\smash{\bigwedge^k V} $ are also called \textit{$k$-covectors} of $V$. There is a canonical isomorphism between $ \smash{\bigwedge^k} V $ and the space of all linear $ \mathbf{R} $-valued maps on $ \bigwedge_k V $. Following\cite[1.4.1]{Fed69}, we write
	$$ \langle \xi, h \rangle := h(\xi) \quad\textrm{for every $ \xi \in {\textstyle\bigwedge_k} V $ and $ h \in {\textstyle\bigwedge^k} V $}. $$
	If $V$ is an inner product space, then $\bigwedge_k V$ and 
	$\smash{\bigwedge^k} V$ inherit natural inner products and the corresponding norms are both denoted by $|\pmb{\cdot}|$  (cf. \cite[1.7.5]{Fed69}).
	
	Let $ U \subseteq \mathbf{R}^d$ be open set and $ k\in\mathbf{N}$. A  \emph{$k$-dimensional differential form} (or simply, a $k$-form) is a smooth map $ \smash{\phi : U \rightarrow \bigwedge^k \mathbf{R}^d} $ (with the convention $\smash{\bigwedge^0 \mathbf{R}^d := \mathbf{R} }$). Following \cite[$4.1.7$]{Fed69}, we denote by $ \mathcal{E}^k(U) $  the space of all smooth $ k $-forms on $ U $, and by denote by $ \mathcal{D}^k(U) $ the space of all smooth $ k $-forms with compact support in $ U $, both equipped with the standard locally convex topologies described in \cite[$4.1.1$]{Fed69}.
	We denote by
	\begin{equation*}
		E_{n+1}:=\pmb{e}_1\wedge\cdots\wedge\pmb{e}_{n+1}\in\textstyle{\bigwedge_{n+1}}\mathbf{R}^{n+1}\quad\!\textup{and}\quad E'_{n+1}:=\pmb{e}'_1\wedge\cdots\wedge\pmb{e}'_{n+1}\in\textstyle{\bigwedge^{n+1}}\mathbf{R}^{n+1}
	\end{equation*}
	the \emph{standard orientation} and \emph{standard volume form} of $\mathbf{R}^{n+1}$, respectively. Given $\phi \in \mathcal{E}^k(U)$, 
	we denote by $d\phi$ its \emph{exterior derivative} 
	(cf.~\cite[$4.1.6$]{Fed69}). If $f : U \rightarrow \mathbf{R}^m$ is 
	smooth and $\psi$ is a $k$-form on an open set $V \subseteq \mathbf{R}^m$ 
	with $f(U) \subseteq V$, the \emph{pull back} $f^{\#}\psi$ is the $k$-form 
	on $U$ defined by
	$$ \langle v_1 \wedge \cdots \wedge v_k, f^{\#} \psi(x) \rangle := \langle \big[\textstyle\bigwedge_{k}\Der f(x)\big](v_1 \wedge \cdots \wedge v_k), \psi\big(f(x)\big) \rangle$$
	for every $ x \in U $ and $ v_1, \ldots , v_k \in \mathbf{R}^d $. We refer to \cite[4.1.6]{Fed69} for the basic properties of the operator $ f^\# $.
	
	Functions from a subset of $ U $ to $ \bigwedge_k\mathbf{R}^m $ are called \emph{$ k $-vectorfields}.
	
	Let $k\geq0$ be an integer. A \emph{$ k $-current} is a continuous, real-valued linear operator on $ \mathcal{D}^k(U) $. The space of all $ k $-currents on $ U $ is denoted by $ \mathcal{D}_k(U) $.
	
	A sequence $ \{T_\ell\}_{\ell\in\mathbf{N}} \subset \mathcal{D}_k(U)$ \emph{weakly$^\ast$ converges} to $ T\in \mathcal{D}_k(U) $ if
	$$ T_\ell(\phi) \to T(\phi) \quad \textup{for every $ \phi \in \mathcal{D}^k(U) \, .$} $$
	The \emph{support} of $T \in \mathcal{D}_k(U)$ is
	$$
	\spt(T) := U \setminus \bigcup\big\{ V \subseteq U : 
	\text{$V$ open and $T(\phi) = 0$ whenever 
		$\phi \in \mathcal{D}^k(U)$ with $\spt(\phi) \subset V$} \big\}\,,
	$$
	and its \emph{boundary} is the $(k-1)$-current 
	$\partial T \in \mathcal{D}_{k-1}(U)$ defined by
	$$
	\partial T(\phi) := T(d\phi) \quad \text{for every } \phi \in \mathcal{D}^k(U)\,,
	$$
	note that $\spt(\partial T) \subseteq \spt(T)$.

	If $T \in \mathcal{D}_k(U)$ has compact support in $U$, then $T$ extends 
	uniquely to a continuous linear functional on $\mathcal{E}^k(U)$. Indeed, for any choice of 
	$\gamma \in \mathcal{D}^0(U)$ with 
	$\spt(T) \subseteq \mathrm{interior}\,[\gamma^{-1}(\{1\})]$, 
	the extension is
	$$
	T(\varphi) := T(\gamma\,\varphi) \quad \text{for every } 
	\varphi \in \mathcal{E}^k(U)\,.
	$$
	Given $\psi \in \mathcal{E}^h(U)$ and $T \in \mathcal{D}_k(U)$ with 
	$h \leq k$, the \emph{restriction of $T$ to $\psi$} is
	$$
	(T \restrict \psi)(\phi) := T(\psi \wedge \phi) 
	\quad \text{for every } \phi \in \mathcal{D}^{k-h}(U)\,.
	$$
	Given $T \in \mathcal{D}_k(U)$ and a smooth map $f : U \rightarrow V$ 
	onto an open set $V \subseteq \mathbf{R}^m$ such that $f|\spt(T)$ is 
	proper, we define the \emph{push-forward} $f_{\#}T \in \mathcal{D}_k(V)$ by
	\begin{equation}\label{eq push forward current}
		f_\# T(\phi) := T\big(\gamma\, f^\#\phi\big) 
		\quad \text{for every } \phi \in \mathcal{D}^k(V)\,,
	\end{equation}
	for any choice of $\gamma \in \mathcal{D}^0(U)$ such that 
	$f^{-1}(\spt(\phi)) \cap \spt(T) \subseteq 
	\mathrm{interior}\,[\gamma^{-1}(\{1\})]$. 
	Note that $\spt(f_{\#}T) \subseteq f(\spt(T))$. 
	When $\spt(T)$ is compact in $U$,
	$$
	f_\# T(\phi) = T(f^\# \phi) \quad \text{for every } \phi \in \mathcal{E}^k(V).
	$$ 

	We say that $ T \in \mathcal{D}_k(U) $ is an \emph{integer multiplicity locally rectifiable $ k $-current} if
	$$ T(\phi) = \int_{M}\langle \vv{\eta}(x), \phi(x) \rangle \, d\mathcal{H}^k(x) \quad \textrm{for every $ \phi \in \mathcal{D}^k(U) $} $$
	where $ M \subseteq U $ is $ \mathcal{H}^k $-measurable and countably $ \mathcal{H}^k $-rectifiable, and $ \vv{\eta} $ is an $ (\mathcal{H}^k \restrict M) $-measurable $ k $-vectorfield satisfying:\emph{\begin{enumerate}[label=(\roman*)]
			\item  \textup{$ \int_{K \cap M} | \vv{\eta}|\, d\mathcal{H}^k < \infty $ for every compact subset $ K $ of $ U $;}
			\item  \textup{$ \vv{\eta}(x) $ is simple and $ | \vv{\eta}(x) | $ is a positive integer for $ \mathcal{H}^k $-a.e.\ $ x \in  M ;$}
			\item  \textup{$ \textup{Tan}^k(\mathcal{H}^k \restrict M,x) $ is associated with $ \vv{\eta}(x) $ for $ \mathcal{H}^k $-a.e.\ $ x \in M $.}\end{enumerate}}
	The set  $ M $ is called \emph{carrier} of $T$ and is $ \mathcal{H}^k $-almost uniquely determined by $ T $. Its \emph{canonical representative} is given by (cf. \cite[7]{FedererNote})
	$$W_T:=\{x\in U:0<\Theta^k(\|T\|,x)<\infty\}\,,$$
	where $\|T\|=|\vv{\eta}|\,\mathcal{H}^k\restrict M$. Defining the \emph{multiplicity} and \emph{orientation} of $T$ by $$i_T:=|\vv{\eta}| \qquad \textup{and} \qquad \vv{\eta}_T:=\frac{\vv{\eta}}{|\vv{\eta}|}\,,$$ 
	we obtain the representation
	$$ T(\phi) = \int_{W_T}\langle \vv{\eta}_T(x), \phi(x) \rangle\,i_T(x) \, d\mathcal{H}^k(x) \quad \textrm{for every $ \phi \in \mathcal{D}^k(U) $} $$
	and it is customary to denote this compactly by $T=i_T\,(\mathcal{H}^k\restrict W_T)\wedge\vv{\eta}_T$. We denote the class of integer multiplicity locally rectifiable $ k $-current in $U$ as $\mathcal{R}_k(U)$. Notice that
	$$\spt(T)=\overline{W}_T\quad\textup{and}\quad\mathcal{H}^k\big(\spt(T)\setminus W_T\big)=0\,.$$
	
	Now suppose $T\in\mathcal{R}_k(U)$, such that $T=i_T\,(\mathcal{H}^k\restrict W_T)\wedge\vv{\eta}_T$, and let $ f : U \rightarrow V $ be a Lipschitz map onto an open subset $ V $ of $ \mathbf{R}^m $ such that $f|\spt(T)$ is proper. The push-forward $f_\# T \in \mathcal{D}_k(V)$ is defined by
	\begin{equation}\label{pushfor}
		\!f_\# T(\varphi):=\int_{W_T}\!\langle\big[\textstyle\bigwedge_k\textup{ap}\,Df(x)\big]\big(\vv{\eta}_T(x)\big),\varphi\big(f(x)\big)\rangle \, i_T(x)\,d\mathcal{H}^k(x)\quad\textup{for every $\varphi\in\mathcal{D}^k(V)\,.$}
	\end{equation}
	In particular, the area formula implies that $f_\# T \in \mathcal{R}_k(V)$ and is carried by $f(W_T)$ (cf.~\cite[4.1.30]{Fed69} or 
	\cite[Remark~27.2\,(3)]{SimonBook}).
	
	\subsection{Legendrian currents}\label{LegCyc} Here we introduce the notion of Legendrian cycles and we collect some fundamental facts.
	
	\begin{definition}[Contact $1$-form] We say that $\alpha \in \mathcal{E}^1(\mathbf{R}^{n+1} \times \mathbf{R}^{n+1})$ is the \emph{contact $1$-form} of $\mathbf{R}^{n+1}$ if it acts as follows
		$$ \langle (y,v),\alpha(x,u) \rangle := y \bullet u \quad \textup{for every} \ y,v,x,u\in\mathbf{R}^{n+1} .$$
		Equivalently, $\alpha(x,u) := \sum_{i=1}^{n+1} u\bullet \pmb{e}_i \, \pi_0^{\#}\pmb{e}'_i$ for all $x,u\in\mathbf{R}^{n+1}$.
	\end{definition}
	
	\begin{definition}[Legendrian cycle]\label{def Legendrian cycle}
		\emph{Let $ W \subseteq \mathbf{R}^{n+1} $ be an open set. We say that an integer-multiplicity locally rectifiable $ n $-current $  T $ of $ W \times \mathbf{R}^{n+1} $ is a Legendrian cycle of $W$ if the following conditions hold:
			\begin{enumerate}[label=(\roman*)]
				\item $ \spt(T) \subseteq W \times \mathbf{S}^n $;
				\item $ \partial T =0$\,;
				\item $ T \restrict \alpha =0 $, where  $ \alpha  \in \mathcal{E}^1(\mathbf{R}^{n+1} \times \mathbf{R}^{n+1}) $ is the contact $ 1 $-form of $ \mathbf{R}^{n+1} $.
		\end{enumerate}}
	\end{definition}
	
	\begin{lemma}\label{lem patching of legendrian cycles}
		Suppose $W_1, \ldots , W_m \subseteq \mathbf{R}^{n+1} $ are bounded open sets and $ T \in \mathcal{D}_{n}(\mathbf{R}^{n+1} \times \mathbf{R}^{n+1}) $ such that $ T\restrict (W_i \times \mathbf{R}^{n+1}) $ is a Legendrian cycle of $ W_i $ for every  $ i\in\{1, \ldots , m\} $ and $ \spt(T) $ is a compact subset of $ \bigcup_{i=1}^m W_i \times \mathbf{S}^n $. Then $ T   $ is a Legendrian cycle of $ \mathbf{R}^{n+1} $.
	\end{lemma}
	\begin{proof}
		For every $ i\in\{1, \ldots , m\} $ choose an open set $ V_i $  with compact closure in $ W_i $ and a function $ f_i \in C^\infty_c(\mathbf{R}^{n+1}) $ such that $ \spt(f_i) $ is a compact subset of $ W_i $\,, $ \sum_{i=1}^m f_i(x) = 1 $ for every $ x \in \bigcup_{i=1}^m V_i $ and $ \spt(T) \subseteq \bigcup_{i=1}^m V_i \times \mathbf{S}^n $. Then $ T = T \restrict\sum_{i=1}^m(f_i\circ\pi_0) $, $ \sum_{i=1}^m d(f_i\circ\pi_0) =0 $ on $ \bigcup_{i=1}^m V_i\times\mathbf{S}^n$, 
		$$ (T\restrict\alpha)(\phi) = \sum_{i=1}^m (T\restrict \alpha)\big((f_i\circ\pi_0)\,\phi\big)=0 $$
		and 
		$$ \partial T(\phi) = \sum_{i=1}^m \partial T\big((f_i\circ\pi_0)\,\phi\big) - T\Big(\big({\textstyle\sum_{i=1}^m}d(f_i\circ\pi_0)\big)\wedge \phi \Big) =0 \ ,$$
		for every $ \phi \in \mathcal{D}^{n-1}(\mathbf{R}^{n+1} \times \mathbf{R}^{n+1}) $.
	\end{proof}
	
	\begin{definition}[Lipschitz-Killing forms]\label{Lipschitz Killing}
		\emph{For $ k \in \{0, \ldots , n\} $	the $ k $-th Lipschitz-Killing differential form $ \varphi_k \in \mathcal{E}^n(\mathbf{R}^{n+1} \times \mathbf{R}^{n+1})$  is defined by
			{\allowdisplaybreaks\begin{align*}
					\langle \xi_1 \wedge \cdots \wedge \xi_n, \varphi_k(x,u) \rangle:= \sum_{\sigma \in \Sigma_{n,k}} \langle \pi_{\sigma(1)}(\xi_1)\wedge \cdots \wedge \pi_{\sigma(n)}(\xi_n) \wedge u, \pmb{e}'_1 \wedge \cdots \wedge \pmb{e}'_{n+1}\rangle \, ,
			\end{align*} }
			for every $ \xi_1, \ldots , \xi_n \in \mathbf{R}^{n+1}\times \mathbf{R}^{n+1}$, where 
			$$ \Sigma_{n,k} := \bigg\{ \sigma: \{1, \ldots , n\} \rightarrow \{0,1\}: \sum_{i=1}^n \sigma(i) = n-k \bigg\} \, . $$}
	\end{definition}
	Given a $ C^2 $-diffeomorphism $ F :  U \rightarrow  V $ between open subsets of $ \mathbf{R}^{n+1} $ we define the $ C^1 $-diffeomorphism $ \Psi_F  : \mathbf{R}^{n+1} \times \mathbf{S}^n \rightarrow \mathbf{R}^{n+1} \times \mathbf{S}^n $ by
	\begin{equation}\label{PsiF}
		\Psi_F(x,u) := \bigg(F(x), \frac{(D F(x)^{-1})^\ast(u)}{| (D F(x)^{-1})^\ast(u) |}\bigg) \quad\textrm{for every $(x,u) \in \mathbf{R}^{n+1} \times \mathbf{S}^n $} \, .
	\end{equation} 
	We recall the following result by Joseph Fu \cite[Lemma 3.1]{Fu98} on the exterior derivative of the Lipschitz-Killing forms (for a proof, see also \cite[Lemma 2.4]{SantilliValentini}).
	
	\begin{lemma}\label{Lemma Fu}
		Let $ T $ be a compactly supported Legendrian cycle of $ \mathbf{R}^{n+1} $, $ \smash{\{F_t\}_{t \in (-\epsilon, \epsilon)} }$  is a  local variation\footnote{A \emph{local variation} $\smash{\{F_t\}_{t \in (-\epsilon, \epsilon)}} $ of $ \mathbf{R}^{n+1} $ is a family of smooth maps over $ \mathbf{R}^{n+1} $, varying smoothly in the time variable $t $, such that $ F_0 $ is the identity of $ \mathbf{R}^{n+1} $, $ F_t : \mathbf{R}^{n+1} \rightarrow \mathbf{R}^{n+1} $ is a diffeomorphism for every $ t \in (-\epsilon, \epsilon) $ and the set $ \{x \in \mathbf{R}^{n+1}: F_t(x) \neq x\} $ has compact closure in $\mathbf{R}^{n+1} $ for every $t\in(-\epsilon,\epsilon)$. The \emph{initial velocity vector field} $ V $ of $\{F_t\}_{t \in (-\epsilon, \epsilon)} $ is the smooth vector field $ V $ of $ \mathbf{R}^{n+1} $ defined as 
			$$ V(x) := \lim_{t \to 0} \frac{F_t(x) - x}{t} \quad\textrm{for every $ x \in \mathbf{R}^{n+1} $.} $$} of $ \mathbf{R}^{n+1} $ with initial velocity vector field $ V $ and
		$$\theta_V : (x,u)\in\mathbf{R}^{n+1} \times \mathbf{R}^{n+1} \mapsto V(x) \bullet u\in\mathbf{R} \, .$$
		Then
		$$ \frac{d}{dt}\big((\Psi_{F_t})_{\#}T \big)(\varphi_{i}) \Big|_{t=0} = (n+1-i)\,T(\theta_V \, \varphi_{k-1}) \quad\textrm{for every $ k \in\{ 1, \ldots , n\} $} $$
		and 
		$$  \frac{d}{dt}\big((\Psi_{F_t})_{\#}T \big)(\varphi_{0})  \Big|_{t=0}  =0 \, . $$
	\end{lemma}
	
	The following result describes the approximate tangent space of the carrier of a Legendrian cycle.
	
	\begin{theorem}[\protect{cf.\ \cite[Theorem 9.2]{RatajZaehlebook}}]\label{LegCycCarrier}
		Let $T$ be a Legendrian cycle of $\mathbf{R}^{n+1}$, with carrier $W_T$. For $\mathcal{H}^n$-a.e. $(x,u)\in W_T$ there exist numbers $$-\infty<\kappa_1(x,u)\leq\ldots\leq \kappa_n(x,u)\leq+\infty$$ and vectors $\tau_1(x,u),\ldots,\tau_n(x,u)$ such that $ \{\tau_1(x,u), \ldots , \tau_n(x,u), u\} $ is a positively oriented orthonormal basis of $\mathbf{R}^{n+1} $ $($i.e. $\smash{\langle\textstyle{\bigwedge_{i=1}^n}\tau_{i}(x,u)\wedge u,\pmb{e}'_1\wedge\cdots\wedge\pmb{e}'_{n+1}\rangle=1)}$ and the vectors 
		$$\xi_i(x,u):=\begin{cases}
			\big(\tau_i(x,u),\kappa_i(x,u)\tau_i(x,u)\big) \quad \textup{\emph{if $\kappa_i(x,u)<+\infty$}}\\
			\big(0,\tau_i(x,u)\big) \, \, \, \, \, \quad\quad\quad\quad\quad\quad\textup{\emph{if $\kappa_i(x,u)=+\infty$}}
		\end{cases} \quad i\in\{1,\ldots,n\}$$
		form an orthogonal basis of $ \textup{Tan}^n(\mathcal{H}^n \restrict Q, (x,u)) $ for every $ \mathcal{H}^n $-measurable set $ Q \subseteq W_T $ with $ \mathcal{H}^n(Q) < \infty $ and for $ \mathcal{H}^n $-a.e.\ $(x,u) \in Q $.  
		
		The maps $ \kappa_{1}, \ldots , \kappa_n $ are $ \smash{(\mathcal{H}^n \restrict W_T)} $-almost uniquely determined, as well as the substaces of $\mathbf{R}^{n+1}$ spanned by vectors $\tau_j(x,u)$ belonging to a fixed value among the $\kappa_i(x,u)$ $(i\in\{1,\ldots,n\})$ is uniquely determinated.
	\end{theorem}
	
	\begin{definition} \emph{Given $T$, $\{\tau_1,\ldots,\tau_n\}$ and $\{\xi_1,\ldots,\xi_n\}$ as in Theorem $\ref{LegCycCarrier}$, we define
			\begin{equation*}
				\vv{\xi}_T(x,u):=\frac{\xi_1(x,u)\wedge\cdots\wedge\xi_n(x,u)}{|\xi_1(x,u)\wedge\cdots\wedge\xi_n(x,u)|}\in\bigwedge\nolimits_n(\mathbf{R}^{n+1}\times\mathbf{R}^{n+1}) \, ,
			\end{equation*}
			\begin{equation*}
				\zeta_T(x,u):=\frac{1}{|\xi_1(x,u)\wedge\cdots\wedge\xi_n(x,u)|}\in(0,+\infty)\,,
			\end{equation*}
			for $\mathcal{H}^n$-a.e. $(x,u)\in W_T$.}\end{definition}
	
	\begin{remark}
		\textup{This definition is well-posed, in the sense that $\vv{\xi}_T$ and $\zeta_T$ do not depend, up to a set of $\mathcal{H}^n$-measure zero, on the choice of $\{\tau_1,\ldots,\tau_n\}$ (cf. \cite[Remark 3.1]{Santilli24}).}
	\end{remark}
	
	\begin{definition}[Principal curvatures of Legendrian cycles] \emph{Let $T$ be a Legendrian cycle of $\mathbf{R}^{n+1}$, with carrier $W_T$. We define the principal curvatures of $T$ as $K_{T,i}:=\kappa_i$\,, where $\kappa_{1}\leq\ldots\leq\kappa_{n}$ are the functions defined $\mathcal{H}^n$-a.e. on $W_T$ given by Theorem $\ref{LegCycCarrier}$.}
	\end{definition}

	\begin{definition}[$k$-th mean curvature of Legendrian cycles] \emph{Let $T$ be a Legendrian cycle of $\mathbf{R}^{n+1}$, with carrier $\Sigma_T$. We define the sets
			{\allowdisplaybreaks\begin{align*}
					&W^{(i)}_T:=\big\{(x,u)\in W_T:K_{T,i}(x,u)<+\infty\,,\,K_{T,i+1}(x,u)=+\infty\} \quad \textup{\emph{for}} \ i\in\{1,...,n-1\} \, ,&\\
					&W^{(0)}_T:=\big\{(x,u)\in W_T:K_{T,1}(x,u)=+\infty\} \, ,&\\
					&W^{(n)}_T:=\big\{(x,u)\in\Sigma_T:K_{T,n}(x,u)<+\infty\}\,.
			\end{align*}}
			Then we define the $k$-th mean curvature function of $T$ as
			{\allowdisplaybreaks\begin{align*}
					&H_{T,k}:=\pmb{1}_{W^{(n-k)}_T}+\sum_{i=1}^k \, \sum_{\lambda\in\bigwedge(n-k+i,i)}K_{T,\lambda(1)}\ldots K_{T,\lambda(i)} \, \pmb{1}_{W^{(n-k+i)}_T} \quad \textup{\emph{for}} \ k\in\{1,\dots,n\} \, ,&\\
					&H_{T,0}:=\pmb{1}_{W^{(n)}_T}\,,
			\end{align*}}
			where $\bigwedge(m,h)$, for $h\leq m$, is the set of all increasing functions from $\{1,\dots,h\}$ to $\{1,\dots,m\}$.}
	\end{definition}
	
	The definition of the mean curvature functions is motivated by the following result.
	
	\begin{lemma}[\protect{cf.\ \cite[Lemma 3.2]{Santilli24}}]\label{FuFirst9}
		If $T=i_T(\mathcal{H}^n\restrict W_T)\wedge\vv{\eta}_T$ is a Legendrian cycle of $\mathbf{R}^{n+1}$ and $k\in\{0,\ldots,n\}$, then
		\begin{equation*}
			(T\restrict\varphi_k)(\phi)=\int_{W_T}\phi(x,u)\,i_T(x,u)\,\zeta_T(x,u)\,H_{T,n-k}(x,u) \, d\mathcal{H}^n(x,u)
		\end{equation*}
		for every $\phi\in\mathcal{D}^0(\mathbf{R}^{n+1}\times\mathbf{R}^{n+1})$.
	\end{lemma}
	
	\subsection{Normal bundle} 
	For an arbitrary subset $C \subseteq \mathbf{R}^{n+1} $ we define the \emph{distance function} $\pmb{\delta}_C$ from $C$ as
	\begin{equation*}
		\smash{\pmb{\delta}_C(x):=\inf\big\{|x-a|:a\in C\big\}\in[0,+\infty) \quad\text{for every }x\in\mathbf{R}^{n+1}\, .}
	\end{equation*}
	\begin{definition}[\protect{Proximal unit normal bundle; cf. \cite[p.\,212]{RockafellarWets}}] \emph{Given $C\subseteq\mathbf{R}^{n+1}$, we define the proximal unit normal bundle of $C$ as the set
			\begin{equation*}
				\textup{nor}(C):=\big\{(x,\nu)\in \overline{C}\times\mathbf{S}^n:\pmb{\delta}_C(x+s\nu)=s \ \textup{\emph{for some}} \ s>0\big\} \, .
			\end{equation*}
			Moreover, we define
			{\allowdisplaybreaks\begin{align*}
					\qquad\qquad  &\textup{nor}(C,x):=\big\{\nu\in\mathbf{S}^n:(x,\nu)\in\textup{nor}(C)\big\}\quad \textup{\emph{for}} \ x\in \overline{C} \,, &\\
					&\textup{nor}(C)\restrict E:=\big\{(x,\nu)\in\textup{nor}(C):x\in E\big\} \quad\textup{\emph{for}} \ E\subseteq\mathbf{R}^{n+1} \, .
		\end{align*}}}
	\end{definition}
	\begin{remark}
		Notice that $ \textup{nor}(C) = \textup{nor}(\overline{C}) $. We recall that $ \textup{nor}(C) $ is a Borel set and it is always countably $ \mathcal{H}^n $-rectifiable (cf.\,\cite[Remark 4.3]{SantilliAnnali}\footnote{The proximal unit normal bundle of a closed set $ C $ in \cite{SantilliAnnali} is denoted with $N(C)$.}). Moreover, we say that $\textup{nor}(C)$ satisfies the \emph{Lusin $(N)$-property} on $E\subset\mathbf{R}^{n+1}$ if $\mathcal{H}^n\big(\textup{nor}(C)\restrict E)=0$, whenever $\mathcal{H}^n(E)=0$.
	\end{remark}
	
	It is proved in \cite{SantilliAnnali} that there exists a subset $\widetilde{\textup{nor}}(C)$ of $\textup{nor}(C)$ (cf.\,\cite[Definition 4.4]{SantilliAnnali}), where $\smash{\mathcal{H}^n\big(\textup{nor}(C)\setminus\widetilde{\textup{nor}}(C)\big)=0}$ (cf.\,\cite[Remark 4.5]{SantilliAnnali}), such that for every $(x,u)\in\widetilde{\textup{nor}}(C)$ there exists a linear subspace $\smash{T_C(x,u)}$ of $u^\perp$ and a symmetric bilinear form $$Q_C(x,u):T_C(x,u)\times T_C(x,u)\to\mathbf{R}\,,$$
	whose eigenvalues can be used to provide an explicit representation of the approximate tangent space of $\textup{nor}(C)$ at $\mathcal{H}^n$-a.e. points. To this end, for every $(x,u)\in\widetilde{\textup{nor}}(C)$, we define
	$$-\infty<\kappa_{1}(x,u)\leq\ldots\leq\kappa_{n}(x,u)\leq+\infty$$
	in the following way:
	{\allowdisplaybreaks\begin{align*}
			\quad&\textup{\emph{\ \ \ \ \ \ \ \ \ \ \ \ \ \ \ \ \ \ \   $\kappa_{1}(x,u),\ldots,\kappa_{m}(x,u)$ are the eigenvalues of $Q_C(x,u)$\,,}}&\\
			&\textup{\ \ \ \ \ \ \ \ \ \ \emph{where $m=\dim T_C(x,u)$, and $\kappa_{i}(x,u)=+\infty$ for any $i\in\{m+1,\ldots,n\}$\,.}}
	\end{align*}}
	The following lemma is a simple extension of well known results for sets of positive reach (cf.\ \cite[Proposition 4.23 and Lemma 4.24]{RatajZaehlebook}). 
	\begin{lemma}\label{lem: Santilli20}
		Suppose $ C \subseteq \mathbf{R}^{n+1} $. Then, for $ \mathcal{H}^n $-a.e.\!\,\,\,$(x,u)\in\textup{nor}(C)$, there exist the vectors $ \tau_1(x,u), \ldots , \tau_n(x,u)\in\mathbf{R}^{n+1}$ such that: $ \{\tau_1(x,u), \ldots , \tau_n(x,u), u\} $ is a positively oriented orthonormal basis of $\mathbf{R}^{n+1} $ $($i.e. $\smash{\langle\textstyle{\bigwedge_{i=1}^n}\tau_{i}(x,u)\wedge u,\pmb{e}'_1\wedge\cdots\wedge \pmb{e}'_{n+1}\rangle=1)}$ and the vectors 
		$$ \xi_i(x,u) := \bigg(\frac{1}{\sqrt{1 + \kappa_i(x,u)^2}} \,\tau_i(x,u), \frac{\kappa_i(x,u)}{\sqrt{1 + \kappa_i(x,u)^2}}\,\tau_i(x,u) \bigg)\quad\text{for }i \in\{ 1, \ldots, n \}$$
		form an orthonormal basis of $ \textup{Tan}^n(\mathcal{H}^n \restrict Q, (x,u)) $ for every $ \mathcal{H}^n $-measurable set $ Q \subseteq\textup{nor}(C) $ with $ \mathcal{H}^n(Q) < \infty $ and for $ \mathcal{H}^n $-a.e.\ $(x,u) \in Q $ $($we set $ \frac{1}{\infty} = 0 $ and $ \frac{\infty}{\infty} = 1 )$.  
		
		The maps $ \kappa_{1}, \ldots , \kappa_n $ can be chosen to be $\smash{\big(\mathcal{H}^n \restrict \textup{nor}(C)\big)} $-measurable and they are $ \smash{\big(\mathcal{H}^n \restrict \textup{nor}(C)\big)} $-almost uniquely determined.
	\end{lemma}
	
	\begin{proof}
		The existence part of the statement and the measurability property are discussed in \cite[Section 4]{SantilliAnnali} and \cite[Section 3]{HugSantilli} (see in particular \cite[Remark 3.7]{HugSantilli}). While uniqueness can be proved as in \cite[Lemma 4.24]{RatajZaehlebook}. 
	\end{proof}
	
	\begin{definition}\label{curv}
		\emph{Suppose $ C \subseteq \mathbf{R}^{n+1} $, we denote by $ \kappa_{C,1}, \ldots , \kappa_{C,n} $ the $ \smash{\big(\mathcal{H}^n \restrict \textup{nor}(C)\big)} $-measurable maps given by Lemma $\ref{lem: Santilli20}$.}
	\end{definition}
	
	\begin{definition}[Reach function of $C$] \emph{Suppose $ C \subseteq \mathbf{R}^{n+1} $, we define the reach function of $C$, at $(x,u)\in\textup{nor}(C)$, as
			\begin{equation*}
				\pmb{r}_C(x,u):=\sup\big\{r>0:\pmb{\delta}_C(x+ru)=r\big\}\,.
		\end{equation*}}
	\end{definition}
	
	\begin{remark}\textup{Let $C\subseteq\mathbf{R}^{n+1}$. Given  $(x,u)\in\widetilde{\textup{nor}}(C)$ and $r>0$ such that $\pmb{\delta}_C(x+ru)=r$, we have that (cf. \cite[Lemma 4.8]{SantilliAnnali})
			$$Q_C(x,u)(\tau,\tau)\geq-\frac{|\tau|^2}{r}\quad \textup{whenever $\tau\in T_C(x,u)$}\,.$$
			Therefore, for every $(x,u)\in\widetilde{\textup{nor}}(C)$, we deduce that
			$$-\frac{1}{\pmb{r}_C(x,u)}\leq\kappa_{C,i}(x,u)\leq+\infty \quad\textup{for every $i\in\{1,\ldots,n\}$}\,.$$}
	\end{remark}
	
	\begin{definition}\label{defn-vect}\emph{For $\mathcal{H}^n$-a.e. $(x,u)\in\textup{nor}(C)$ we define
			\begin{equation*}
				\vv{\xi}_C(x,u):=\frac{\xi_{C,1}(x,u)\wedge\cdots\wedge\xi_{C,n}(x,u)}{|\xi_{C,1}(x,u)\wedge\cdots\wedge\xi_{C,n}(x,u)|}\in\bigwedge\nolimits_n(\mathbf{R}^{n+1}\times\mathbf{R}^{n+1}) \, ,
			\end{equation*}
			\begin{equation*}
				\zeta_C(x,u):=|\xi_{C,1}(x,u)\wedge\cdots\wedge\xi_{C,n}(x,u)|^{-1}\in(0,+\infty)\,,
			\end{equation*}
			where, for any $i\in\{1,\ldots,n\}$ and with the notations of Lemma $\ref{lem: Santilli20}$ and Definition $\ref{curv}$, we set
			\begin{equation*}\label{Ufset43}
				\xi_{C,i}(x,u):=\begin{cases}
					\big(\tau_{i}(x,u),\kappa_{C,i}(x,u)\,\tau_{i}(x,u)\big) \quad\textup{\emph{if}} \ \kappa_{C,i}(x,u)<+\infty\\
					\big(0,\tau_{i}(x,u)\big) \quad\quad\quad\quad\quad\quad\quad\ \ \textup{\emph{if}} \ \kappa_{C,i}(x,u)=+\infty
				\end{cases}
			\end{equation*}
			for $\mathcal{H}^n$-a.e. $(x,u)\in\textup{nor}(C)$.}
	\end{definition}
	
	\begin{remark}
		\textup{The definitions of $\vv{\xi}_C$ and $\zeta_C$ do not depend on the 
			choice of $\{\tau_1, \ldots, \tau_n\}$, if the latter is consistent with 
			Lemma~\ref{lem: Santilli20}, in the sense that for any two choices 
			$\{\tau_1, \ldots, \tau_n\}$ and $\{\tau_1', \ldots, \tau_n'\}$ satisfying 
			Lemma~\ref{lem: Santilli20} at $(x,u) \in \textup{nor}(C)$, we have
			\begin{equation}\label{contradiction}
				\bigwedge_{i=1}^n\xi_{C,i}(x,u)=\bigwedge_{i=1}^n\xi_{C,i}'(x,u)\,,
			\end{equation} 
			where the vectors $\{\xi_{C,1}(x,u),\ldots,\xi_{C,n}(x,u)\}$ and $\{\xi_{C,1}'(x,u),\ldots,\xi_{C,n}'(x,u)\}$, given by
			$$\xi_{C,i}:=
			\begin{cases}
				\big(\tau_{i},\kappa_{C,i}\,\tau_{i}\big) \quad  \textup{if} \ \kappa_{C,i}<+\infty\\
				(0,\tau_{i}) \quad\quad\, \, \, \, \, \, \ \textup{if} \ \kappa_{C,i}=+\infty
			\end{cases}\text{and} \quad \xi'_{C,i}:=
			\begin{cases}
				\big(\tau'_{i},\kappa_{C,i}\,\tau'_{i}\big) \quad\textup{if} \ \kappa_{C,i}<+\infty\\
				(0,\tau'_{i}) \quad\quad\, \, \, \, \, \, \, \, \textup{if} \ \kappa_{C,i}=+\infty
			\end{cases}\!\!\!\!\!\!,$$
			form two orthogonal bases of $\textup{Tan}^n(\mathcal{H}^n \restrict Q, (x,u)). $}
	\end{remark}

	\begin{lemma}\label{lem legendrian property of the normal bundle}
		Let $C \subseteq \mathbf{R}^{n+1}$ and let $\alpha$ be the contact $1$-form 
		of $\mathbf{R}^{n+1}$ $($cf.\ Definition~\ref{def Legendrian cycle}$)$. 
		If $Q \subseteq \textup{nor}(C)$ is $\mathcal{H}^n$-measurable with 
		$\mathcal{H}^n(Q) < \infty$, then
		$$
		\langle \xi, \alpha(x,u) \rangle = 0
		$$
		for every $\xi \in \textup{Tan}^n\big(\mathcal{H}^n \restrict Q, (x,u)\big)$ 
		and $\mathcal{H}^n$-a.e.\ $(x,u) \in Q$.
	\end{lemma}
	
	\begin{proof}
		This follows from the definition of $\alpha$ and Lemma~\ref{lem: Santilli20}.
	\end{proof}
	\begin{definition}\label{Ufset13}
		\emph{Given $C\subseteq\mathbf{R}^{n+1}$, we define the following subsets of $C$
			{\allowdisplaybreaks\begin{align*}
					&\ \ \ \ \ \ \ \ \ \ \ \quad\quad N_1(C):=\big\{x\in C:\mathcal{H}^0\big(\textup{nor}(C,x)\big)=1\big\} \, , &\\
					&\ \ \ \ \ \ \ \ \ \ \ \quad\quad N_2(C):=\big\{x\in C:\mathcal{H}^0\big(\textup{nor}(C,x)\big)=2\big\} \, ,&\\
					&\ \ \ \ \ \ \ \ \ \ \ \quad\quad N_\infty(C):=\big\{x\in C:\mathcal{H}^0\big(\textup{nor}(C,x)\big)=\infty\big\} \, .
		\end{align*}}}
	\end{definition}
	\begin{definition}[Distance cone of $C$]\emph{Given a closed set $C \subseteq \mathbf{R}^{n+1}$ and $x \in \pi_0\big(\textup{nor}(C)\big)$, 
			we define the distance cone of $C$ at $x$ as
			$$
			\textup{Dis}(C,x) := \big\{v \in \mathbf{R}^{n+1} : \pmb{\delta}_C(x+v) = |v|\big\}.
			$$
			We denote by $\textup{aff\,Dis}(C,x)$ the affine hull of $\textup{Dis}(C,x)$.}
	\end{definition}
	\begin{remark}
		\textup{We notice that $\textup{Dis}(C,x)$ is a closed convex subset of $\textup{Nor}(C,x)$, for any choice of $x\in\pi_0\big(\textup{nor}(C)\big)$ (cf. \cite[Theorem 4.8 (2)]{Fed59}). Furthermore, the following property holds
			\begin{equation}\label{propDis}
				v\in\textup{Dis}(C,x) \ \ \ \Rightarrow \ \ \  tv\in\textup{Dis}(C,x) \,, \ \textup{for every} \ t\in[0,1]
			\end{equation}
			indeed $B_{|tv|}(x+tv)\cap C\subseteq B_{|v|}(x+v)\cap C=\emptyset$ for any $v\in\textup{Dis}(C,x)$ and $t\in[0,1]$. Moreover
			\begin{equation}\label{norDis}
				\textup{nor}(C,x)=\Big\{\medmath{\frac{u}{|u|}}:u\in\textup{Dis}(C,x)\setminus\{0\}\Big\}\,.
		\end{equation}}
	\end{remark}

	\begin{lemma}\label{lemmaDis}
		Let $C\subseteq\mathbf{R}^n$ be a closed set. If there exists $m\in\{0,\ldots,n-1\}$ such that
		\begin{equation}\label{affhull1}
			\dim\big(\textup{aff}\,\textup{Dis}(C,x)\big)=n-m+1 \quad \textup{for} \ x\in\pi_0\big(\textup{nor}(C)\big)\,,
		\end{equation}
		then \begin{equation*}
			\mathcal{H}^{n-m}\big(\textup{nor}(C,x)\big)\in(0,+\infty)\,.
		\end{equation*}
		Similarly, if
		\begin{equation}\label{affhull2}
			\dim\big(\textup{aff}\,\textup{Dis}(C,x)\big)=1 \quad\textup{for} \ x\in\pi_0\big(\textup{nor}(C)\big)\,,
		\end{equation}
		then \begin{equation*}
			\mathcal{H}^{0}\big(\textup{nor}(C,x)\big)\in\{1,2\}\,.
		\end{equation*}
	\end{lemma}
	\begin{proof} Let $x\in\pi_0\big(\textup{nor}(C)\big)$ such that (\ref{affhull1}) holds, first we notice that $\mathcal{H}^{n-m}\big(\textup{nor}(C,x)\big)<\infty$. In fact
		$$\textup{nor}(C,x)\subseteq \textup{aff}\,\textup{Dis}(C,x)\cap\mathbf{S}^n\cong\mathbf{S}^{n-m}\,.$$

		Now we prove that $\mathcal{H}^{n-m}\big(\textup{nor}(C,x)\big)>0$, let $\epsilon\in(0,1)$ and we introduce the Lipschitz mappings $\rho_\epsilon$ defined as
		$$\rho_\epsilon:v\in\textup{Dis}_\epsilon(C,x)\mapsto\frac{v}{|v|}\in\textup{nor}(C,x)\,,\quad\textup{Dis}_\epsilon(C,x):=\textup{Dis}(C,x)\cap\big(B_{\frac{1}{\epsilon}}(0)\setminus B_\epsilon(0)\big)$$
		where $\textup{Dis}_\epsilon(C,x)$ is $\mathcal{H}^{n-m+1}$-rectifiable (cf. (\ref{affhull1})) and with strictly positive $\mathcal{H}^{n-m+1}$-measure. The last property follows from the fact that $\textup{Dis}(C,x)$ is convex and $\smash{\textup{aff}\,\textup{Dis}(C,x)}$ has dimension $n-m+1$, hence $\textup{Dis}(C,x)$ contains an $(n-m+1)$-dimensional simplex (cf. \cite[Problems 3.5.4 (14)]{Leonard}). We notice that
		$$\mathcal{H}^1\big(\rho_\epsilon^{-1}(u)\big)>0 \quad\textup{for every} \ u\in\textup{Im}(\rho_\epsilon)$$ indeed for each $u\in\textup{Im}(\rho_\epsilon)$ there exists $s\in(\epsilon,\epsilon^{-1})$ such that $su\in\textup{Dis}(C,x)$ and applying (\ref{propDis}) we deduce that $\{tu:t\in[\epsilon,s]\}\subseteq\rho_\epsilon^{-1}(u)$, furthermore by \cite[Theorem 3.2.31]{Fed69} we infer that $\textup{Im}(\rho_\epsilon)$ is also $\mathcal{H}^{n-m}$-rectifiable. Overall, by coarea formula for rectifiable sets (cf. \cite[Theorem 3.2.22]{Fed69}), we obtain
		\begin{equation}\label{coarea}
			\int_{\textup{Dis}_\epsilon(C,x)}J^{\textup{Dis}(C,x)}_{n-m}\rho_\epsilon\,d\mathcal{H}^{n-m+1}=\int_{\textup{Im}(\rho_\epsilon)}\mathcal{H}^1\big(\rho_\epsilon^{-1}(u)\big)\,d\mathcal{H}^{n-m}(u)\,.
		\end{equation}We remember that our goal is to prove that $\mathcal{H}^{n-m}\big(\textup{nor}(C,x)\big)>0$ and to this purpose, since $\smash{\textup{Im}(\rho_\epsilon)\subseteq\textup{nor}(C,x)}$, we use (\ref{coarea}) to show that $\smash{\mathcal{H}^{n-m}\big(\textup{Im}(\rho_\epsilon)\big)>0}$. By virtue of the previous discussion, we have only to prove that
		$$J^{\textup{Dis}(C,x)}_{n-m}\rho_\epsilon(z)>0 \quad\textup{for} \ \mathcal{H}^{n-m+1}\textup{-a.e.} \ z\in\textup{Dis}_\epsilon(C,x)\,.$$
		By the locality property of the approximate tangent spaces (cf. \cite[$3.2.16$]{Fed69}) and (\ref{affhull1}), applying \cite[Lemma $3.2.17$]{Fed69} we deduce that
		\begin{align*}
			\textup{Tan}^{n-m+1}\big(\mathcal{H}^{n-m+1}\restrict \textup{Dis}_\epsilon(C,x), z\big)&=\textup{Tan}^{n-m+1}\big(\mathcal{H}^{n-m+1}\restrict\textup{aff}\,\textup{Dis}(C,x) , z\big)&\\
			&=\textup{aff}\,\textup{Dis}(C,x) \quad \textup{for} \ \mathcal{H}^{n-m+1}\textup{-a.e.} \ z\in\textup{Dis}_\epsilon(C,x)
		\end{align*}
		therefore
		\begin{align*}
			&J^{\textup{Dis}(C,x)}_{n-m}\rho_\epsilon(z) =\sup \big\{ \big|[{\textstyle \bigwedge_{n-m} D \rho_\epsilon(z)}](\xi) \big| :\xi \in {\textstyle \bigwedge_{n-m} \textup{aff}\,\textup{Dis}(C,x)}, \, | \xi | = 1   \big\}
		\end{align*}
		for $\mathcal{H}^{n-m+1}\textup{-a.e.} \ z\in\textup{Dis}_\epsilon(C,x)$. For every fixed $z\in\textup{Dis}_\epsilon(C,x)$, since $\ker D \rho_\epsilon(z)=\textup{span}\{z\}$, we deduce from (\ref{affhull1}) that there exist $v_1,\ldots,v_{n-m}\in\textup{aff\,Dis}(C,x)$ linearly independent and such that $v_i\perp z$ for every $i\in\{1,\ldots,n-m\}$, hence
		$$J^{\textup{Dis}(C,x)}_{n-m}\rho_\epsilon(z)\geq \bigg|\big[{\textstyle \bigwedge_{n-m} D \rho_\epsilon(z)}\big]\bigg(\frac{v_1\wedge\cdots\wedge v_{n-m}}{|v_1\wedge\cdots\wedge v_{n-m}|}\bigg)\bigg|>0$$for $\mathcal{H}^{n-m+1}\textup{-a.e.} \ z\in\textup{Dis}_\epsilon(C,x)$. Overall we conclude that $\mathcal{H}^{n-m}\big(\textup{nor}(C,x)\big)\in(0,+\infty)$.
		
		Let now $\smash{x\in\pi_0\big(\textup{nor}(C)\big)}$ such that (\ref{affhull2}) holds, we prove that $\smash{\mathcal{H}^{0}\big(\textup{nor}(C,x)\big)\in\{1,2\}}$. Again by \cite[Problems 3.5.4 (14)]{Leonard}, we infer that $\textup{Dis}(C,x)$ contains a segment containing the origin. Namely, there exists $u\in\mathbf{S}^n$ such that either $\smash{\{tu:t\in[0,s]\}\subseteq\textup{Dis}(C,x)}$ for some $s>0$, or
		$$\big\{tu:t\in[-s_1,s_2]\big\}\subseteq\textup{Dis}(C,x) \quad\textup{for some} \ s_1,s_2>0\,.$$In the first case we infer that $\mathcal{H}^{0}\big(\textup{nor}(C,x)\big)=1$, while in the second $\mathcal{H}^{0}\big(\textup{nor}(C,x)\big)=2$.
	\end{proof}
	
	\begin{definition}[$m$-th stratum of $C$]\label{stratum}\emph{Suppose $C\subseteq\mathbf{R}^{n+1}$. For every $m\in\{0,\ldots,n\}$ we define the $m$-th stratum of $C$ by
			$$ C^{(m)}:=\big\{x\in C:0<\mathcal{H}^{n-m}\big(\textup{nor}(C,x)\big)<\infty\big\}\,.$$}
	\end{definition}
	
	\begin{remark}\label{propertyStrat} \textup{We notice that the family $\smash{\{C^{(0)},\ldots,C^{(n)}\}}$ is disjoint. Moreover, it is proved (cf. \cite[Remark 5.2]{SantilliAnnali} and \cite[Theorem 4.12]{MenSan}) that every $C^{(m)}$ is a Borel set, also countably $m$-rectifiable and countably $\mathcal{H}^n$-rectifiable of class $2$.}
	\end{remark}
	\begin{corollary}\label{corollaryStrat}
		Let $C\subseteq\mathbf{R}^n$ be a closed set, then:
		\begin{enumerate}[label=(\roman*)]
			\item $\pi_0\big(\textup{nor}(C)\big)=\bigcup_{m =0}^{n}C^{(m)}$;
			\item $N_1(C)\cup N_2(C)= C^{(n)}$;
			\item $N_\infty(C)=\bigcup_{m =0}^{n-1}C^{(m)}.$
		\end{enumerate}
	\end{corollary}
	\begin{proof} To prove \emph{(i)}, since every $C^{(m)}\subseteq\pi_0\big(\textup{nor}(C)\big)$, we show that $\pi_0\big(\textup{nor}(C)\big)\subseteq\bigcup_{m =0}^{n}C^{(m)}$. Let $\smash{x\in\pi_0\big(\textup{nor}(C)\big)}$, namely there exist $\nu\in\textup{nor}(C,x)$ and $s>0$ so that $\pmb{\delta}_C(x+t\nu)=t$ for any $t\in[0,s]$, therefore $\{t\nu:t\in[0,s]\}\subseteq\textup{Dis}(C,x)$. Hence $\smash{\dim\big(\textup{aff}\,\textup{Dis}(C,x)\big)\in\{1,\ldots,n+1\}}$, then by Lemma \ref{lemmaDis} the desired result follows.
		
		To prove \emph{(ii)}, we notice that $\smash{N_1(C)\cup N_2(C)\subseteq C^{(n)}}$. To show that $\smash{C^{(n)}\subseteq N_1(C)\cup N_2(C)}$ assume that there exists $\smash{x\in C^{(n)}\setminus\big(N_1(C)\cup N_2(C)\big)}$, namely $\mathcal{H}^0\big(\textup{nor}(C,x)\big)\in\mathbf{N}\setminus\{0,1,2\}$. This means that $\textup{nor}(C,x)$ contains at least two linearly independent unit vectors $u_1,u_2\in\mathbf{S}^n$, then from (\ref{propDis}) we infer that $\{tu_i:t\in[0,s_i]\big\}\subseteq \textup{Dis}(C,x)$ for any $i\in\{1,2\}$ and for some $s_1,s_2\in\mathbf{R}^+$. We deduce that $\smash{\dim\big(\textup{aff}\,\textup{Dis}(C,x)\big)\in\{2,\ldots,n+1\}}$ and by Lemma \ref{lemmaDis} we infer $\smash{\mathcal{H}^0\big(\textup{nor}(C,x)\big)=\infty}$, which is a contradiction. 
		
		To prove \emph{(iii)}, we notice that $\bigcup_{m =0}^{n-1}C^{(m)}\subseteq N_\infty(C)$. Now assume that $x\in N_\infty(C)$, again $\textup{nor}(C,x)$ contains at least two linearly independent unit vectors, so $\smash{\dim\big(\textup{aff}\,\textup{Dis}(C,x)\big)}$ is at least $2$ and applying Lemma \ref{lemmaDis} we infer that $\smash{x\in\bigcup_{m =0}^{n-1}C^{(m)}}$.
	\end{proof}
	
	\begin{remark}
		\textup{As a consequence of Lemma~\ref{lemmaDis} and 
			Corollary~\ref{corollaryStrat},
			\begin{equation*}
				\textup{\emph{if $x\in N_2(C)$, then there exists $u\in\mathbf{S}^n$ such that $\textup{nor}(C,x)=\{u,-u\}$.}}
			\end{equation*}
			Again, if there exist $u_1,u_2\in\mathbf{S}^n$ linearly independent such that $\textup{nor}(C,x)=\{u_1,u_2\}$, then from (\ref{propDis}) and by Lemma \ref{lemmaDis} we infer that $\smash{x\in C^{(m)}}$ for some $m\in\{0,\ldots,n-1\}$, which is a contradiction since $\smash{x\in C^{(n)}}$ (cf. Corollary \ref {corollaryStrat} \emph{(ii)}).}
	\end{remark}
	
	\begin{lemma}\label{LemmaSelection} Assume that $C\subseteq\mathbf{R}^{n+1}$ is closed, then the following statements hold.
		\begin{enumerate}[label=(\roman*)]
			\item $N_1(C),\,N_2(C)$ and $N_\infty(C)$ are Borel subset of $C$\,.
			\item Given $F$ a $C^2$-diffeomorphism of $\mathbf{R}^{n+1}$, then
			\begin{equation*}
				F\big(N_i(C)\big)=N_i\big(F(C)\big)\quad \textup{\emph{for}} \ i=1,2,\infty\,.
			\end{equation*}
			\item The multivalued function $$\textup{nor}(C,\pmb{\cdot})|N_2(C):N_2(C)\to\mathcal{P}(\mathbf{S}^n)$$ admits an $\mathcal{H}^n$-measurable selection $\nu_C:N_2(C)\to\mathbf{S}^n$ $($namely $\nu_C(x)\in\textup{nor}(C,x)$ for every $x\in N_2(C))$. Moreover,
			\begin{equation}\label{propertyselected}
				\nu_C(p)\in\textup{Nor}^n(\mathcal{H}^n\restrict C,p) \quad\textup{\emph{for every} $p\in N_2(C)$}\,.
			\end{equation}
		\end{enumerate}
	\end{lemma}
	\begin{proof} To prove \emph{(i)}, we notice immediately that $N_\infty(C),\,N_1(C)\cup N_2(C)$ and $\pi_0\big(\textup{nor}(C)\big)$ are Borel subsets of $C$ (cf.\,Remark \ref{propertyStrat} and Corollary \ref{corollaryStrat}). Then, if we consider a countable dense subset $\mathscr{D}$ of $\mathbf{S}^n$, our intention is to obtain the desired result by proving that
		\begin{equation}\label{representation}
			X\setminus N_\infty(C)=N_2(C)\,,
		\end{equation}where $X$ is the Borel subset of $C$ defined as
		$$X:=\bigcup_{h=1}^\infty\,\bigcap_{k=1}^\infty\,\bigcup_{j=k}^\infty\,\bigcup_{(s,u,t,v)\in\mathcal{F}_{h,j}}\!\!C_{k,(s,u,t,v)}\,,$$
		\begin{align*}
			&C_{k,(s,u,t,v)}:=\bigg\{x\in C:\bigg|\frac{\pmb\delta_C(x+su)}{s}-1\bigg|\leq k^{-1}\,,\,&\\
			&\quad\quad\quad\quad\quad\quad\quad\quad\quad\quad\quad\quad\quad\quad\bigg|\frac{\pmb\delta_C(x-tv)}{t}
			-1\bigg|\leq k^{-1}\bigg\} \quad\textup{for $k\in\mathbf{N}^+$, $(s,u,t,v)\in\mathcal{F}_{h,j}$}
		\end{align*}
		and
		{\allowdisplaybreaks\begin{align*}
				\quad&\mathcal{F}_{h,j}:=\big\{(s,u,t,v)\in\mathbf{Q}^+\times\mathscr{D}\times\mathbf{Q}^+\times\mathscr{D}:h^{-1}\leq s,t\leq h\,,&\\
				&\quad\quad\quad\quad\quad\quad\quad\quad\quad\quad\quad\quad\quad\quad\quad\quad\quad\quad\quad\quad|(s,u)-(t,v)|\leq j^{-1}\big\}\quad \textup{for} \  h,j\in\mathbf{N}^+\,.
		\end{align*}}
		
		To prove \eqref{representation}, let $x \in N_2(C)$, that is, assume that 
		there exists $u \in \mathbf{S}^n$ so that $\pmb{\delta}_C(x \pm su) = s$ 
		for some $s > 0$. Since $\pmb{\delta}_C$ is $1$-Lipschitz and 
		$\mathbf{Q}^+ \times \mathscr{D}$ is dense in $\mathbf{R}^+ \times \mathbf{S}^n$, 
		we can find sequences $\{s_k\}_{k \in \mathbf{N}} \subset \mathbf{Q}^+$ and 
		$\{u_k\}_{k \in \mathbf{N}} \subset \mathscr{D}$ such that
		$$
		s_k \xrightarrow[k\to\infty]{} s\,, \quad 
		u_k \xrightarrow[k\to\infty]{} u \quad \textup{and} \quad 
		\frac{\pmb{\delta}_C(x \pm s_k u_k)}{s_k} \xrightarrow[k\to\infty]{} 1\,,
		$$
		from which we readily infer that $x \in X$. Since $x \notin N_\infty(C)$, 
		we deduce that $N_2(C) \subseteq X \setminus N_\infty(C)$.

		Now, let $x\in X\setminus N_\infty(C)$. Namely there exists $h\in\mathbf{N}^+$ so that, for every $k\in\mathbf{N}^+$, there exist an integer $m_k\geq k$ and $(s_k,u_k,t_k,v_k)\in\mathbf{Q}^+\times\mathscr{D}\times\mathbf{Q}^+\times\mathscr{D}$ such that
		$$\bigg|\frac{\pmb\delta_C(x+s_ku_k)}{s_k}-1\bigg|\leq k^{-1}\,,\quad\bigg|\frac{\pmb\delta_C(x-t_kv_k)}{t_k}
		-1\bigg|\leq k^{-1}$$
		$$h^{-1}\leq s_k,t_k\leq h\,,\quad|(s_k,u_k)-(t_k,v_k)|\leq\frac{1}{m_k}\,.$$
		Up to a subsequence we can assume that 
		$$u_k\xrightarrow[k\to\infty]{}u\in\mathbf{S}^n\,,\quad v_k\xrightarrow[k\to\infty]{}v\in\mathbf{S}^n\,,$$
		$$s_k\xrightarrow[k\to\infty]{}s\in\mathbf{R}^+\,,\quad t_k\xrightarrow[k\to\infty]{}t\in\mathbf{R}^+\,,$$
		then we infer that $\pmb\delta_C(x\pm su)=s$. Since $x\notin N_1(C)\cup N_\infty(C)$ and
		\begin{equation}\label{LemmaSelection1}
			\pi_0\big(\textup{nor}(C)\big)=N_1(C)\cup N_2(C)\cup N_\infty(C)\,,
		\end{equation}
		we deduce that $x\in N_2(C)$, namely $X\setminus N_\infty(C)\subseteq N_2(C)$. From (\ref{representation}), since $N_\infty(C)$ and $X$ are Borel subset of $C$, we conclude that $N_2(C)$ is a Borel subset of $C$. Then, from (\ref{LemmaSelection1}), also $N_1(C)$ is a Borel subset of $C$.
		
		To prove \emph{(ii)}, we notice that $\smash{F(C^{(m)})=F(C)^{(m)}}$ for any $\smash{m\in\{0,\ldots,n\}}$ (cf. \cite[Lemma 2.1]{SantilliRectVarCurv}). Hence, from Corollary \ref{corollaryStrat}, we infer that
		$$F\big(N_\infty(C)\big)=N_\infty\big(F(C)\big)\quad \textup{and} \quad F\big(N_1(C)\big)\cup F\big(N_2(C)\big)=N_1\big(F(C)\big)\cup N_2\big(F(C)\big)$$
		and if we prove that
		\begin{equation}\label{LemmaSelection2}
			\smash{F\big(N_1(C)\big)\cap N_2\big(F(C)\big)=\emptyset\,,}
		\end{equation}
		\begin{equation}\label{LemmaSelection3}
			\smash{F\big(N_2(C)\big)\cap N_1\big(F(C)\big)=\emptyset\,,}
		\end{equation}we conclude also that $\smash{F\big(N_1(C)\big)=N_1\big(F(C)\big)}$ and $\smash{F\big(N_2(C)\big)=N_2\big(F(C)\big)}$. To prove (\ref{LemmaSelection2}) assume that there exist $x\in N_1(C)$, $u\in\mathbf{S}^n$ and $s>0$ such that $B_s(F(x)\pm su)\cap F(C)=\emptyset$, namely $\Omega_{\pm}\cap C=\emptyset$ where $\smash{\Omega_{\pm}:=F^{-1}\big[B_s(F(x)\pm su)\big]}$ are two disjoint $C^2$-regular domains with $x\in\partial\Omega_+\cap\partial\Omega_-$ and $\textup{Tan}(\partial\Omega_+,x)=DF^{-1}(x)(u^\perp)=\textup{Tan}(\partial\Omega_-,x)\in\pmb{\textup{G}}(n+1,n)$ by \cite[$3.1.21$]{Fed69}. Since $\Omega_\pm$ are $C^2$-regular domains, so they satisfy the interior sphere condition\footnote{\textbf{Interior sphere condition.} Given $\Omega\subseteq\mathbf{R}^n$ an open set, we say that $\Omega$ satisfies the \emph{interior sphere condition} if for every $x\in\partial\Omega$ there exist $\nu\in\mathbf{S}^{n-1}$ and $r>0$ such that $B_r(x+r\nu)\subseteq\Omega$, notice that $x\in\partial B_r(x+r\nu)$. The interior sphere condition always holds if $\Omega$ is a $C^2$-regular domain (cf. \cite[Remark 4.3.8]{Qing}), in this situation we have $\textup{Tan}(\partial\Omega,x)=\textup{Tan}(\partial B_r(x+r\nu),x)=\nu^\perp$ (the proof is the same as that performed for (\ref{tangent})).}, there exists $r>0$ such that
		$$B_r(x+r\nu)\subseteq\Omega_+ \quad\textup{and} \quad B_r(x-r\nu)\subseteq\Omega_-\,,$$ 
		where $\nu\in\mathbf{S}^n$ is chosen so that $\nu^\perp=\textup{Tan}(\partial\Omega_+,x)=\textup{Tan}(\partial\Omega_-,x)$. Hence $B_r(x\pm r\nu)\cap C=\emptyset$, namely a contradiction since $x\in N_1(C)$. Similarly, to prove (\ref{LemmaSelection3}), assume that there exist $x\in C$, $v\in\mathbf{S}^n$ and $s>0$ such that $B_s(x\pm sv)\cap C=\emptyset$ and $\smash{F(x)\in N_1\big(F(C)\big)}$. Namely, $\Omega_{\pm}\cap F(C)=\emptyset$ where $\smash{\Omega_{\pm}:=F\big[B_s(x\pm sv)\big]}$ are two disjoint $C^2$-regular domains so that $F(x)\in\partial\Omega_+\cap\partial\Omega_-$ and $\textup{Tan}\big(\partial\Omega_+,F(x)\big)=\textup{Tan}\big(\partial\Omega_-,F(x)\big)\in\pmb{\textup{G}}(n+1,n)$, again by the $C^2$-regularity of $\Omega_\pm$ we obtain a contradiction since $\smash{F(x)\in N_1\big(F(C)\big)}$.
		
		To prove \emph{(iii)}, we intend to apply \cite[Theorem \textsl{III}$.6$]{CastaingValadierBook}. That is, if $$N_2(C;U):=\big\{x\in N_2(C):\textup{nor}(C,x)\cap U\neq\emptyset\big\}$$
		is an $\mathcal{H}^n$-measurable subset of $C$ for any open set $U$ in $\mathbf{S}^n$, then the multivalued function $\smash{\textup{nor}(C,\pmb{\cdot}):N_2(C)\to\mathcal{P}(\mathbf{S}^n)}$ admits an $\mathcal{H}^n$-measurable selection $\nu_C:N_2(C)\to\mathbf{S}^n$ (in order to apply \cite[Theorem \textsl{III}$.6$]{CastaingValadierBook}, we notice that for any $\smash{x\in N_2(C)}$ the set $\textup{nor}(C,x)$ is compact and hence complete). Let $U$ be an open set in $\mathbf{S}^n$, we notice that 
		\begin{equation}\label{N^C_2(U)}
			N_2(C;U)=\pi_0\big(\textup{nor}(C)\cap(N_2(C)\times U)\big)
		\end{equation}
		where $\textup{nor}(C)\cap(N_2(C)\times U)$ is a countably $\mathcal{H}^n$-rectifiable Borel subset of $\mathbf{R}^{n+1}$, hence with $\sigma$-finite $\mathcal{H}^n$-measure (cf. \cite[Lemma 15.5 (1)]{MattilaBook}). Namely there exists a family $\{X_i\}_{i\in\mathbf{N}}$ of $\mathcal{H}^n$-measurable sets, with finite $\mathcal{H}^n$-measure, such that $\textup{nor}(C)\cap(N_2(C)\times U)=\bigcup_{i =1}^\infty X_i$. Since $\mathcal{H}^n$ is a Borel-regular outer measure, for every $i\in\mathbf{N}$ (cf. \cite[Proposition 6.2 \emph{(ii)}]{GiaquintaModica_book})
		$$\mathcal{H}^n(X_i)=\sup\big\{\mathcal{H}^n(F):F\subseteq X_i\,,\, F \ \textup{closed}\big\}\,,$$
		so there exists a sequence $\smash{\{F_{i,j}\}_{j\in\mathbf{N}}}$ of closed subsets of $X_i$ such that $\mathcal{H}^n(X_i\setminus\bigcup_{j =1}^\infty F_{i,j})=0$. Since any closed set in the Euclidean space is countable union of compact sets, there exists a sequence $\smash{\{K_{i,j}\}_{j\in\mathbf{N}}}$ of compact subsets of $X_i$ such that $\smash{\mathcal{H}^n(X_i\setminus\bigcup_{j =1}^\infty K_{i,j})=0}$, overall
		\begin{equation*}
			\mathcal{H}^n\Big(\big[\textup{nor}(C)\cap\big(N_2(C)\times U\big)\big]\setminus\bigcup_{i,j =1}^\infty K_{i,j}\Big)=0\,.
		\end{equation*}
		We infer that $N_2(C;U)$ is an $\mathcal{H}^n$-measurable subset of $C$, indeed from (\ref{N^C_2(U)}) we obtain
		$$N_2(C;U)=\bigcup_{i,j=1}^\infty\pi_0(K_{i,j})\cup\pi_0\Big(\big[\textup{nor}(C)\cap\big(N_2(C)\times U\big)\big]\setminus\bigcup_{i,j =1}^\infty K_{i,j}\Big)\,,$$
		where any $\smash{\pi_0(K_{i,j})}$ is compact and $\smash{\pi_0\big(\big[\textup{nor}(C)\cap\big(N_2(C)\times U\big)\big]\setminus\bigcup_{i,j=1}^\infty K_{i,j}\big)}$ is $\mathcal{H}^n$-negligible.
		
		The proof of (\ref{propertyselected}) follows immediately as a consequence of $\textup{nor}(C,x)\subseteq\textup{Nor}(C,x)$ (cf. \cite[Theorem 4.8 (2)]{Fed59}) and $\textup{Tan}^n(\mathcal{H}^n\restrict C,x)\subseteq\textup{Tan}(C,x)$ (cf.\,\cite[$1.1.4$]{RatajZaehlebook}), for every $x\in C$. The proof is complete.
	\end{proof}
	
	\subsection{Fine properties of $W^{2,n} $-graphs}
	We recall some definitions and results given in \cite{SantilliValentini}.
	
	\begin{definition} \textit{Given an open set $U \subseteq \mathbf{R}^n$ and $f \in C^0(U)$, we define 
			$\mathcal{S}^\ast(f)$ $(\mathcal{S}_\ast(f)$, resp.$)$ as the set of $x \in U$ for which there exists a 
			polynomial $P$, of degree at most $2$, with $P(x) = f(x)$ and
			$$
			\limsup_{y \to x} \frac{f(y) - P(y)}{|y-x|^2} < \infty\quad  \bigg( \liminf_{y \to x} \frac{f(y) - P(y)}{|y-x|^2} > - \infty, \ \textrm{resp.}\bigg)\,.
			$$
			Finally, $\mathcal{S}(f)$ is the set of $x \in U$ for which there exists 
			such a $P$ with
			$$
			\lim_{y \to x} \frac{f(y) - P(y)}{|y-x|^2} = 0\,.
			$$}
	\end{definition}

	The next result establishes that every $\smash{W^{2,n}}$-function is pointwise 
	twice differentiable $\mathcal{L}^n$ almost everywhere. It follows from a more general 
	theorem of Calder\'on and Zygmund asserting that functions in $W^{2,p}(U)$ 
	with $p > \frac{n}{2}$ admit a second-order Taylor expansion 
	$\mathcal{L}^n$-a.e.\ on $U$ (cf.~\cite{CalderonZygmund} or 
	\cite[Proposition~2.2]{Caffarelli96}). It can also be obtained by adapting 
	the method of~\cite{Trudinger89} (cf.~\cite[Theorem~2.1.1]{Val25}).
	
	\begin{theorem}\label{W2n twice diff}
		Let $U \subseteq \mathbf{R}^n$ be open and $f \in C^0(U) \cap W^{2,n}_{\textup{loc}}(U)$. 
		Then $\mathcal{L}^n\big(U \setminus \mathcal{S}(f)\big) = 0$.
	\end{theorem}
	We now recall a result by Ulrich Menne (cf.\,\cite[Appendix B]{menne2024sharplowerboundmean}), for a proof see also \cite[Lemma 2.1.3]{Val25}.
	
	\begin{lemma}[Menne]\label{W2n basic estimate}\!Given $ a \in \mathbf{R}^n $\!, $ r > 0 $,\,$ f \in C^0\big(B_r(a)\big) \cap W^{2,n}\big(B_r(a)\big) $\,and $ g \in C^2\big(B_r(a)\big) $ such that 
		$$ g(a) = f(a) \quad \textup{\emph{and}} \quad f(x) \geq g(x)\textrm{ for every $ x \in B_r(a) $\,,} $$
		then $ f $ is pointwise differentiable at $ a $ with $ D f(a) = D g(a) $.
	\end{lemma}
	
	\begin{definition}\label{DefPhi1}
		\emph{Let $\mathbf{S}^{n}_+:=\{z\in\mathbf{S}^n:z\bullet\pmb{e}_{n+1}>0\}$ denote the 
			open upper hemisphere, and let $\psi : \mathbf{R}^n \rightarrow \mathbf{S}^{n}_+$ be the bilipschitz diffeomorphism defined by
			\[
			\psi(y) := \frac{(-y, 1)}{\sqrt{1 + |y|^2}}\,.
			\] 
			Given an open set $U \subseteq \mathbf{R}^n$ and $f \in C^0(U) \cap W^{2,n}_\mathrm{loc}(U)$, we define
			\[
			\Phi^\pm_f(x) := \pigl(x, f(x), \pm\psi\big(\nabla f(x)\big)\pigl) \quad \text{for every } x \in \textup{$\textsf{Diff}$}(f)\,.
			\]}
	\end{definition}
	
	The following result is proved using a Rad\'o-Reichelderfer type argument (cf. \cite[Lemma 3.4]{SantilliValentini} or \cite[Lemma 2.1.4]{Val25}).
	
	\begin{lemma}\label{W2n functions Lusin} Given $ U\subseteq\mathbf{R}^n$ an open set and $ f \in C^0(U) \cap W^{2,n}_{\textup{loc}}(U)$, we have $ \mathcal{H}^n\big(\Phi^\pm_f(Z)\big) =0 $ for every $ Z \subset \mathcal{S}^\ast(f)\cup\mathcal{S}_\ast(f) $ such that $ \mathcal{L}^n(Z) =0 $.
	\end{lemma}
	
	We say that a mapping $g : U \to \mathbf{R}^m$ with $m \geq n$, 
	defined on an open set $U \subseteq \mathbf{R}^n$, satisfies 
	the \emph{Lusin $(N)$-condition} if
	$$
	\mathcal{H}^n\big(g(E)\big) = 0 \quad \text{whenever } 
	E \subset U \text{ satisfies } \mathcal{L}^n(E) = 0\,.
	$$
	From Lemma~\ref{W2n functions Lusin} it easily follows that if 
	$f \in C^0(U) \cap W^{2,n}_{\textup{loc}}(U)$, then both $\nabla f$ 
	and $\overline{\nabla f}$ satisfy the Lusin $(N)$-condition on $\mathcal{S}(f)$. 
	However, the Lusin $(N)$-property does not generally hold for the gradient of $W^{2,n}$-functions. In fact, T.~Roskovec 
	(cf.\ \cite{Roskovec}), using a Cesari-type construction, provides an example 
	of a function $f \in C^1([-1,1]^n)$ such that
	$${\nabla}f\in W^{1,n}\big((-1,1)^n;\mathbf{R}^n\big)\quad\textup{and} \quad [-1,1]^n\subseteq{\nabla}f([-1,1]\times\{0\}^{n-1})\,.$$
	In other words, $\nabla f$ is a $(C^0\cap W^{1,n})$-regular vector field, but it does not satisfy the Lusin $(N)$-property since it maps a segment into an $n$-cube. We also recall that T.~Toro \cite{Toro} constructs a $W^{2,n}$-function with a dense set of singular points, and J.~Fu \cite[p.\,2260]{FuAlexandrov} observes that the gradient of a $W^{2,n}$-function may have a graph dense in $\mathbf{R}^n \times \mathbf{R}^n$.
	
	\begin{theorem}[\protect{cf.\,\cite[Theorem 7.6]{SantilliValentini}}]\label{Nabelpunktsatz} Let $U \subseteq \mathbf{R}^n$ be open and 
		$f \in W^{2,1}_{\textup{loc}}(U)$. If \,$\overline{f}$ satisfies the 
		Lusin $(N)$-condition, then $\Gamma := \overline{f}(U)$ is 
		$\mathcal{H}^n$-rectifiable of class $2$.
	\end{theorem}
	
	\begin{remark}\label{rmk Lusin condition (N) for f}
		\textup{Let $f \in W^{2,p}_{\textup{loc}}(U)$ with 
			$\frac{n}{2} < p < n$, the Sobolev embedding theorem 
			(cf.~\cite[Theorem~7.26]{GilTru_Book}) ensures that 
			$\smash{f \in W^{1,p^\ast}_{\textup{loc}}(U)}$ with $p^\ast > n$.\,Therefore, $\overline{f}$ satisfies the Lusin $(N)$-condition 
			by~\cite[Theorem~1]{MarcusMizel}.}
	\end{remark}

	The next lemma provides an important property of the normal bundle of the graph of a $W^{2,n}$-function (cf. \cite[Lemma 3.5]{SantilliValentini} or \cite[Lemma 2.1.5]{Val25}).
	
	\begin{lemma}\label{lem no vertical touching balls}
		Let $ U \subseteq \mathbf{R}^n $ be open, $ \gamma > \frac{1}{2} $, $ f \in C^{0, \gamma}(U) $, $ x \in U $, and $ \nu \in \mathbf{S}^n \subseteq \mathbf{R}^n \times \mathbf{R} $ be such that $ B^{n+1}_s\big(\overline{f}(x)+s\nu\big) \cap \overline{f}(U) = \varnothing $ for some $ s > 0 $, where $\overline{f}$ is the graph map of $f$.
		It follows that $ \nu \notin \mathbf{R}^n \times \{0\} $. In particular, this holds 
		whenever $ f \in C^0(U)\cap W^{2,n}_{\textup{loc}}(U) $.
	\end{lemma}
	
	\begin{definition}\label{catodef}
		\emph{Let $ U \subseteq \mathbf{R}^n $ be an open set and $ f : U \to \mathbf{R} $. We define the cato-graph and the epi-graph of $f$ as
			$$ C_f := \{(x,u) \in U \times \mathbf{R} :  u \leq f(x)\} \quad\text{and}\quad E_f := \{(x,u) \in U \times \mathbf{R} :  u \geq f(x)\}\,.$$ 
			Moreover, denoting by $\Gamma_f$ the graph of $f$, we define
			$$ N_f := \textup{nor}(C_f) \cap (U \times \mathbf{R} \times \mathbf{R}^{n+1})\,, $$
			$$ M_f := \textup{nor}(E_f) \cap (U \times \mathbf{R} \times \mathbf{R}^{n+1})\,, $$
			$$N(\Gamma_f):=\textup{nor}(\Gamma_f)\cap (U \times \mathbf{R} \times \mathbf{R}^{n+1})\,.$$}
	\end{definition}
	
	The statements of the following lemma concerning $N_f$ are proved 
	in~\cite[Lemma~3.6 and Theorem~3.7]{SantilliValentini} 
	or~\cite[Lemma~2.1.6 and Theorem~2.1.2]{Val25}. The analogous statements for $M_f$ follow by the same argument.
	
	\begin{lemma}\label{W2n functions normal bundle}
		Let $ U\subset\mathbf{R}^n $ be a bounded open set and let $ f \in C^0(U)\cap W^{2,n}(U)$. Then, for every $ A \subseteq U $,
		\begin{align}\label{W2n functions normal bundle par}
			\quad\quad\quad\quad\quad \ \ \ &N_f \cap (A \times \mathbf{R} \times \mathbf{R}^{n+1})  =  \Phi^+_f\big[ A \cap \mathcal{S}^\ast(f) \big]\,,&\\
			\label{W2n functions normal bundle par2}
			&M_f \cap (A \times \mathbf{R} \times \mathbf{R}^{n+1})  =  \Phi^-_f\big[ A \cap \mathcal{S}_\ast(f) \big]\,.
		\end{align}
		Furthermore, there exist two Borel $n$-vectorfields, $\smash{\vv{\eta}_{N_f}}$ defined on $N_f$ and $\smash{\vv{\eta}_{M_f}}$ defined on $M_f$, such that
		$$(\mathcal{H}^n \restrict N_f) \wedge \vv{\eta}_{N_f} \!\quad\textup{\emph{and}} \quad (\mathcal{H}^n \restrict M_f) \wedge \vv{\eta}_{M_f}$$
		are Legendrian cycles of $ U \times \mathbf{R} $, where $\mathcal{H}^n(N_f)<\infty$ and $\mathcal{H}^n(M_f)<\infty$. In particular, for $ \mathcal{H}^n $-a.e.\ $(z,\nu) \in N_f $ we have
		$$\begin{cases}
			| \vv{\eta}_{N_f}(z, \nu)| =1 \ \text{and} \ \textrm{$ \vv{\eta}_{N_f}(z, \nu) $ is simple}\\
			\textrm{$\textup{Tan}^n\big(\mathcal{H}^n \restrict N_f, (z, \nu)\big) $ is associated with $ \vv{\eta}_{N_f}(z, \nu) $}\\
			\langle \big[ {\textstyle\bigwedge_n} \pi_0\big]\big(\vv{\eta}_{N_f}(z,\nu)\big) \wedge \nu , \pmb{e}'_1 \wedge \cdots \wedge \pmb{e}'_{n+1} \rangle > 0\,,
		\end{cases}$$
		
		the analogous properties hold for $\mathcal{H}^n$-a.e.\ 
		$(z,\nu) \in M_f$, with $\vv{\eta}_{M_f}$ in place of $\vv{\eta}_{N_f}$.
		
		\end{lemma}
		
		\begin{remark}\label{EpiAndGraph} We readily infer that
			\begin{equation}\label{norEpiCato}
				N(\Gamma_f)=N_f\cup M_f \quad\textup{and} \quad N_f\cap M_f=\emptyset\end{equation}
			and, from (\ref{W2n functions normal bundle par}) and (\ref{W2n functions normal bundle par2}),
			\begin{equation}\label{W2n functions normal bundle par3}
				\textup{nor}(\Gamma_f)\cap (A \times \mathbf{R} \times \mathbf{R}^{n+1})  =  \Phi^+_f\big[ A \cap \mathcal{S}^\ast(f) \big]\cup\Phi^-_f\big[ A \cap \mathcal{S}_\ast(f) \big]
			\end{equation}
			for every $A\subseteq U$.  In particular, for any 
			$Z \subset \Gamma$ with $\mathcal{H}^n(Z) = 0$,
			{\allowdisplaybreaks\begin{align*}
					N(\Gamma)\restrict Z&=\textup{nor}(\Gamma)\restrict\big(\pi(Z)\times\mathbf{R}\big)&\\
					&=\Phi^+_f\big(\pi(Z)\cap\mathcal{S}^\ast(f)\big)\cup\Phi^-_f\big(\pi(Z)\cap\mathcal{S}_\ast(f)\big)
			\end{align*}}where $\pi$ denotes the canonical projection onto the first $n$ components, 
			so that $\mathcal{L}^n\big(\pi(Z)\big) = 0$. Since $\Phi^+_f$ and $\Phi^-_f$ 
			satisfy the Lusin $(N)$-condition on $\mathcal{S}^\ast(f)$ and 
			$\mathcal{S}_\ast(f)$, respectively (cf.~Lemma~\ref{W2n functions Lusin}), 
			we conclude that $\mathcal{H}^n\big(N(\Gamma) \restrict Z\big) = 0$, namely
			$$\textup{$N(\Gamma)$ \emph{satisfies the Lusin $(N)$-property on $\Gamma$.}}$$Furthermore, the Borel $n$-vectorfield
			\begin{equation}\label{etaGraph}
				\vv{\eta}_\Gamma:=\vv{\eta}_{N_f}\pmb{1}_{N_f}+\vv{\eta}_{M_f}\pmb{1}_{M_f}\quad\textup{on} \ N(\Gamma)
			\end{equation}
			satisfies, for $\mathcal{H}^n$-a.e. $(x,\nu)\in N(\Gamma)$,
			$$ | \vv{\eta}_\Gamma(z, \nu)| =1\, , \quad \textrm{$ \vv{\eta}_\Gamma(z, \nu) $ is simple}\, , $$
			$$ \textrm{$\textup{Tan}^n\big(\mathcal{H}^n \restrict N(\Gamma), (z, \nu)\big) $ is associated with $ \vv{\eta}_\Gamma(z, \nu) $} \,,$$
			\begin{equation}\label{projectionproperty}
				\langle \big[ {\textstyle\bigwedge_n} \pi_0\big]\big(\vv{\eta}_\Gamma(z,\nu)\big) \wedge \nu\, , \pmb{e}'_1 \wedge \cdots \wedge \pmb{e}'_{n+1} \rangle > 0\,.
			\end{equation}
			Overall, from Theorem \ref{W2n functions normal bundle} and (\ref{W2n functions normal bundle par3}), we conclude that:
			\begin{align}\nonumber
				&\,\quad\quad\quad\quad \, \mathcal{N}_\Gamma:=\big(\mathcal{H}^n \restrict 
				N(\Gamma)\big) \wedge \vv{\eta}_\Gamma \ \textup{\emph{is a Legendrian cycle of $ U \times \mathbf{R} $, where}}&\\ \label{LusinNorGraph}
				&\ \, \, \quad\quad\textup{\emph{$N(\Gamma)$ has finite $\mathcal{H}^n$-measure and satisfies the Lusin $(N)$-property on $\Gamma$.}}
			\end{align}
			We denote $\mathcal{N}_\Gamma$ as the \emph{Legendrian cycle associated with $\Gamma$}.
		\end{remark}
		
		\section{Legendrian cycles and Reilly-type variational formulae on $F_nW^{2,n}$-sets}
		
		\subsection{Introduction to $F_nW^{2,n}$-sets} We introduce the class of $ F_nW^{2,n} $-functions following  Ambrosio, Gobbino and Pallara (cf. \cite{AmGobPal}), who developed an idea of De Giorgi.
		
		\begin{definition}[$F_nW^{2,n}$-functions]\label{Ufset5} \emph{Given $\Omega\subseteq\mathbf{R}^{n+1}$ an open set and a function $\pmb\bigiotab:\Omega\to\mathbf{N}$, we say that $\pmb\bigiotab\in F_nW^{2,n}(\Omega)$ if for any $z\in\{\pmb\bigiotab>0\}$ there exist an open neighborhood $U\subseteq \Omega$ of $z$, also a positive integer $q(z)$ and a family $\{\Gamma_1,\ldots,\Gamma_{q(z)}\}$ of subset of $\mathbf{R}^{n+1}$, such that
				\begin{equation}\label{Ufset7}
					\pmb\bigiotab(x)=\sum_{i=1}^{q(z)}\pmb{1}_{\Gamma_i}(x) \quad\textup{\emph{for any} }x\in U
				\end{equation}
				where every $\Gamma_i$ satisfies the following property:
				{\allowdisplaybreaks\begin{align}\nonumber
						&\quad\quad\quad\textup{\emph{there exist $p_i\in \Gamma_i\cap U$, $\eta_i\in\mathbf{S}^n$, a set $V_i$ open in $\eta^\perp_i$ so that $0\in V_i$,}}&\\ \nonumber
						&\quad\quad\quad\textup{\emph{$f_i\in C^0(V_i)\cap W^{2,n}(V_i)$ such that $f_i(0)=0$ which also induces the map}}&\\ \nonumber
						&\quad\quad\quad\quad\quad\quad\quad\quad\quad\quad\ \ \overline{f_i}:x\in V_i\mapsto x+f_i(x)\,\eta_i\in\mathbf{R}^{n+1},&\\ \nonumber
						&\quad\quad\quad\textup{\emph{in such a way that}}&\\ \label{Ufset4}
						&\quad\quad\quad\quad\quad\quad\quad\quad\quad\quad\ \ \ \ \quad\quad\Gamma_i\cap U=\overline{f_i}(V_i)+p_i \,.
				\end{align}}We say $f_i$ a graph map of $\Gamma_i\cap U$, on $V_i$.}
		\end{definition}
		
		\begin{remark}\label{Ufset3}
			Given an open set $U \subset \mathbf{R}^{n+1}$ and 
			$f \in C^0(U) \cap W^{2,n}(U)$, consider $\Gamma := \textup{graph}(f)$. 
			Since $\overline{f}$ satisfies the Lusin $(N)$-condition 
			(cf.~Remark~\ref{rmk Lusin condition (N) for f}), 
			Theorem~\ref{Nabelpunktsatz} implies that $\Gamma$ is 
			$\mathcal{H}^n$-rectifiable of class $2$.
		\end{remark}
		
		\begin{definition}[$F_nW^{2,n}$-sets]\label{Ufset5.5}
			\emph{We say that a closed set $\mathcal{S}\subset\mathbf{R}^{n+1}$ is a $F_nW^{2,n}$-set if 
				\begin{equation*}
					\mathcal{S}=\{\pmb\bigiotab>0\}
				\end{equation*}
				for some $\pmb\bigiotab\in F_nW^{2,n}(\mathbf{R}^{n+1})$. We call $\pmb\bigiotab$ a multiplicity function of $\mathcal{S}$.}
		\end{definition}

		\begin{remark}\label{Ufset6}
			Equivalently, a closed set $\mathcal{S} \subset \mathbf{R}^{n+1}$ is a 
			$F_nW^{2,n}$-set if for every $z \in \mathcal{S}$ there exist a positive 
			integer $q(z)$, an open neighborhood $U \subset \mathbf{R}^{n+1}$ of $z$ 
			(which we may assume bounded without loss of generality) and a family 
			$\{\Gamma_i\}_{i=1}^{q(z)}$ of subsets of $\mathbf{R}^{n+1}$ such that
			\begin{equation}\label{Ufset1}
				\mathcal{S} \cap U = \bigcup_{i=1}^{q(z)} (\Gamma_i \cap U)\,,
			\end{equation}
			where every $\Gamma_i \cap U$ coincides with the graph of a 
			$(C^0 \cap W^{2,n})$-function.
		\end{remark}
		
		\begin{definition}[$\mathscr{W}^{2,n}$-sets]\label{def W2n sets} \emph{We say that that a closed set $\mathcal{S}\subset\mathbf{R}^{n+1}$ is a $\mathscr{W}^{2,n}$-set if there exists a pair $(\mathcal{S}',F)$ satisfying the following properties:
				\begin{enumerate}[label=(\roman*)]
					\item $\mathcal{S}'$ is a $F_nW^{2,n}$-set;
					\item $F(\mathcal{S}')=\mathcal{S}$, where $F$ is a $C^2$-diffeomorphism of $\mathbf{R}^{n+1}$.
			\end{enumerate}}
		\end{definition}
		
		\begin{remark}\label{tangentWsets} \textup{Given a $\mathscr{W}^{2,n}$-set $\mathcal{S}$ with associated pair 
				$(\mathcal{S}', F)$, assume that \eqref{Ufset1} holds in an open 
				neighborhood $U \subset \mathbf{R}^{n+1}$ of $z \in \mathcal{S}'$, 
				where every $\Gamma_i \cap U$ coincides with the graph of a 
				$(C^0 \cap W^{2,n})$-function. Then
				\begin{equation}\label{tanUnion}
					\textup{Tan}\big(\mathcal{S}, F(z)\big) = 
					\bigcup_{i=1}^{q(z)} \textup{Tan}\big(F(\Gamma_i), F(z)\big)\,.
				\end{equation}
				To prove \eqref{tanUnion}, first notice that 
				$\bigcup_{i=1}^{q(z)} \textup{Tan}\big(F(\Gamma_i), F(z)\big) 
				\subseteq \textup{Tan}\big(\mathcal{S}, F(z)\big)$, 
				since \eqref{Ufset1} implies $\textup{Tan}\big(F(\Gamma_i), F(z)\big) 
				\subseteq \textup{Tan}\big(\mathcal{S}, F(z)\big)$ for every 
				$i \in \{1, \ldots, q(z)\}$. Now let 
				$u \in \textup{Tan}\big(\mathcal{S}, F(z)\big) \cap \mathbf{S}^n$,
				then there exists $\{z_j\}_{j \in \mathbf{N}} \subset 
				\mathcal{S} \setminus \{F(z)\}$ such that
				$$
				z_j \xrightarrow[j\to\infty]{} F(z) \qquad \text{and} \qquad 
				\frac{z_j - F(z)}{|z_j - F(z)|} \xrightarrow[j\to\infty]{} u\,.
				$$
				Since $\mathcal{S} \cap F(U)$ is a finite union of $F(\Gamma_i) \cap F(U)$, 
				there exist $h \in \{1, \ldots, q(z)\}$ and a subsequence 
				$\{z_{j_k}\}_{k \in \mathbf{N}} \subseteq \{z_j\}_{j \in \mathbf{N}}$ 
				with $\{z_{j_k}\}_{k \in \mathbf{N}} \subset F(\Gamma_h) \setminus \{F(z)\}$ 
				and
				$$
				z_{j_k} \xrightarrow[k\to\infty]{} F(z) \qquad \text{and} \qquad 
				\frac{z_{j_k} - F(z)}{|z_{j_k} - F(z)|} \xrightarrow[k\to\infty]{} u\,,
				$$
				namely $u \in \textup{Tan}\big(F(\Gamma_h), F(z)\big)$. Hence 
				$\textup{Tan}\big(\mathcal{S}, F(z)\big) \subseteq 
				\bigcup_{i=1}^{q(z)} \textup{Tan}\big(F(\Gamma_i), F(z)\big)$.}
		\end{remark}
		
		\begin{lemma}\label{graphProperties}
			Let $U \subseteq \mathbf{R}^n$ be an open set, $f \in C^0(U) \cap W^{2,n}(U)$, 
			$\Gamma := \textup{graph}(f)$, and $x \in U$. The following properties hold.
			\begin{enumerate}[label=(\roman*)]
				\item Let $g \in C^2(U)$ be such that $f(x) = g(x)$ and $\Gamma$ is 
				contained either in the epi-graph or in the cato-graph of $g$. Setting 
				$\Sigma := \textup{graph}(g)$ and $z := \overline{f}(x) = \overline{g}(x)$,
				\begin{equation}\label{tangent*}
					\textup{Tan}(\Gamma, z) = \textup{Tan}(\Sigma, z)\,.
				\end{equation}
				\item Denoting by $E_f$ and $C_f$ the epi-graph and cato-graph of $f$, 
				respectively,
				\begin{equation}\label{W2,nDom0.15}
					\overline{f}\big(\mathcal{S}^\ast(f) \cap \mathcal{S}_\ast(f)\big) 
					= N_2(\Gamma)\,,
				\end{equation}
				\begin{equation*}
					\overline{f}\big(\mathcal{S}^\ast(f)\big) = N_1(C_f) 
					\quad \text{and} \quad 
					\overline{f}\big(\mathcal{S}_\ast(f)\big) = N_1(E_f)\,.
				\end{equation*}
			\end{enumerate}
		\end{lemma}
		\begin{proof}
			To prove \emph{(i)}, since $\Sigma$ is a $C^2$-regular graph and therefore 
			satisfies the two-sides sphere condition\footnote{\textbf{Two-sides sphere 
					condition.} We say that the graph $\Gamma$ of a continuous function $f$ 
				satisfies the \emph{two-sides sphere condition} if for every $x \in \Gamma$ 
				there exist $\nu \in \mathbf{S}^{n-1}$ and $r > 0$ such that 
				$B_r(x+r\nu) \subseteq E_f$ and $B_r(x-r\nu) \subseteq C_f$, where $E_f$ 
				and $C_f$ are the epi-graph and the cato-graph of $f$, respectively. This 
				always holds if $\Gamma$ is a $C^2$-regular graph 
				(cf.~\cite[Remark~4.3.8]{Qing}), in which case 
				$\textup{Tan}(\Gamma,x) = \textup{Tan}(\partial B_r(x \pm r\nu),x) = \nu^\perp$ 
				(the proof is the same as that of \eqref{tangent}).}, 
			there exist $\nu \in \mathbf{S}^n$ and $s > 0$ such that
			$$
			B_s(z+s\nu) \cap \Gamma = \emptyset \quad \text{and} \quad 
			B_s(z+s\nu) \cap \Sigma = \emptyset\,.
			$$
			More precisely, we claim that
			\begin{equation}\label{tangent}
				\textup{Tan}(\Gamma,z) = \textup{Tan}(\Sigma,z) = \nu^\perp.
			\end{equation}
			First, applying Lemma~\ref{lem no vertical touching balls}, we deduce that 
			$\nu \notin \mathbf{R}^n \times \{0\}$. Then, by the implicit function 
			theorem, there exist $\delta > 0$, $h \in C^\infty\big(B_\delta(x)\big)$ 
			and an open neighborhood $V \subset U \times \mathbf{R}$ of $z$ such that 
			$V \cap \partial B_s(z+s\nu) = \textup{graph}(h)$, $\overline{h}(x) = z$, 
			and (up to a sign)
			$$
			\nu = \frac{(-\nabla h(x), 1)}{\sqrt{1 + |\nabla h(x)|^2}}.
			$$
			Therefore, by Lemma~\ref{W2n basic estimate} and~\cite[$3.1.21$]{Fed69}, 
			we conclude that
			$$
			\textup{Tan}(\Gamma,z) = D\overline{f}(x)[\mathbf{R}^n] = 
			D\overline{g}(x)[\mathbf{R}^n] = \textup{Tan}(\Sigma,z)
			$$
			and
			$$
			D\overline{f}(x)[\mathbf{R}^n] = D\overline{h}(x)[\mathbf{R}^n] = 
			\big\{(v, \nabla h(x) \bullet v) : v \in \mathbf{R}^n\big\} = \nu^\perp.
			$$
			
			To prove \emph{(ii)}, we establish only \eqref{W2,nDom0.15}; the remaining 
			identities follow by the same argument. To show 
			$\overline{f}\big(\mathcal{S}^\ast(f) \cap \mathcal{S}_\ast(f)\big) 
			\subseteq N_2(\Gamma)$, let $x \in \mathcal{S}^\ast(f) \cap \mathcal{S}_\ast(f)$ 
			be arbitrary and set $z := \overline{f}(x)$. From the definitions of 
			$\mathcal{S}^\ast(f)$ and $\mathcal{S}_\ast(f)$, we infer that 
			$z \in N_2(\Gamma) \cup N_\infty(\Gamma)$, and by statement \emph{(i)} we 
			deduce that $\textup{Tan}(\Gamma,z) \in \mathbf{G}(n+1,n)$. If 
			$z \in N_\infty(\Gamma)$, then since 
			$\textup{nor}(\Gamma,z) \subseteq \textup{Nor}(\Gamma,z) \cap \mathbf{S}^n$ 
			(cf.~\cite[Theorem~4.8\,(2)]{Fed59}), we would have
			$$
			\infty = \mathcal{H}^0\big(\textup{nor}(\Gamma,z)\big) \leq 
			\mathcal{H}^0\big(\textup{Nor}(\Gamma,z) \cap \mathbf{S}^n\big) = 2\,,
			$$
			a contradiction. Hence $z \in N_2(\Gamma)$, which is the desired result.
			
			Now we show that $N_2(\Gamma) \subseteq \overline{f}\big(\mathcal{S}^\ast(f) 
			\cap \mathcal{S}_\ast(f)\big)$. Let $z \in N_2(\Gamma)$, so that $z = \overline{f}(x)$ for some $x \in U$, 
			and let $\nu \in \mathbf{S}^n$ and $s > 0$ be such that
			\begin{equation}\label{W2,nDom0.2}
				B^{n+1}_s(z \pm s\nu) \cap \Gamma = \emptyset \quad \text{and} \quad 
				\pi'\big(B^{n+1}_s(z \pm s\nu)\big) \subseteq U,
			\end{equation}
			where $\pi'$ denotes the canonical projection 
			onto the first $n$ components (in what follows, $\pi''$ denotes the 
			canonical projection onto the last component). By 
			Lemma~\ref{lem no vertical touching balls}, $\nu \in \mathbf{S}^n \setminus 
			(\mathbf{R}^n \times \{0\})$, i.e.\ $\pi''(\nu) \neq 0$. Our goal is to 
			prove that
			\begin{equation}\label{W2,nDom0.25}
				B^{n+1}_s(z-s\nu) \subset C_f \quad \text{and} \quad 
				B^{n+1}_s(z+s\nu) \subset E_f \quad (\text{or vice versa})\,,
			\end{equation}
			where $C_f$ and $E_f$ are the cato-graph and epi-graph of $f$.\,Indeed,\,\eqref{W2,nDom0.25} together with the implicit function theorem 
			(using $\pi''(\nu) \neq 0$) and the Taylor expansion gives 
			$N_2(\Gamma) \subseteq \overline{f}\big(\mathcal{S}^\ast(f) \cap 
			\mathcal{S}_\ast(f)\big)$.
			
			If \eqref{W2,nDom0.25} does not hold, then $B^{n+1}_s(z-s\nu)$ and 
			$B^{n+1}_s(z+s\nu)$ are both contained in $C_f$ or both in $E_f$, 
			otherwise \eqref{W2,nDom0.2} would be contradicted. It is not 
			restrictive to assume that
			$$
			B^{n+1}_s(z \pm s\nu) \subset E_f \quad \text{and} \quad \pi''(\nu) > 0\,.
			$$
			Focusing on $B^{n+1}_s(z-s\nu)$, we notice that:
			\begin{enumerate}
				\item $x = \pi'(z) \in \pi'\big(B^{n+1}_s(z-s\nu)\big)$\,;
				\item for every $\xi \in \{y \in \partial B^{n+1}_s(z-s\nu) : 
				(y-z+s\nu) \bullet \pmb{e}_{n+1} \leq 0\}$,
				$$
				f(x) = \pi''(z) = \pi''(z-s\nu) + s\,\pi''(\nu) > 
				\pi''(z-s\nu) \geq \pi''(\xi)\,;
				$$
				\item since $B^{n+1}_s(z-s\nu) \subset E_f$,
				$$
				f(y) \leq \min\big\{t \in \mathbf{R} : y + t\,\pmb{e}_{n+1} \in 
				\overline{B}^{n+1}_s(z-s\nu)\big\} \quad 
				\text{for every } y \in \pi'\big(\overline{B}^{n+1}_s(z-s\nu)\big).
				$$
			\end{enumerate}
			Hence, if $\widehat{z} \in \partial B^{n+1}_s(z-s\nu)$ is such that
			$$
			\big\{y \in \partial B^{n+1}_s(z-s\nu) : 
			(y-z+s\nu) \bullet \pmb{e}_{n+1} \leq 0\big\} \cap (\pi')^{-1}(x) 
			= \{\widehat{z}\}\,,
			$$
			we reach a contradiction since $f(x) > \pi''(\widehat{z}) \geq f(x)$, 
			contradicting the fact that $f$ is a graph.
		\end{proof}
		
		\begin{definition}\emph{Let $\mathcal{S}$ be a $F_nW^{2,n}$-set and assume that \eqref{Ufset1} 
				holds in an open neighborhood $U \subseteq \mathbf{R}^{n+1}$ of 
				$z \in \mathcal{S}$. We define the map
				\begin{equation}\label{newmap}
					\pmb{\mathcal{I}}:x\in\mathcal{S}\cap U\mapsto\big\{i\in\{1,\dots,q(z)\}:x\in\Gamma_i\cap U\big\}\in\mathbf{N}\,.\end{equation}}
		\end{definition}
		\begin{remark}
			\textup{We notice that $\mathcal{H}^0\big(\pmb{\mathcal{I}}(x)\big)=\pmb\bigiotab(x)$ for every $x\in\mathcal{S}\cap U$.}
		\end{remark}
		
		Now we collect some fine properties of $\mathscr{W}^{2,n}$-sets.
		
		\begin{lemma}\label{Ufset14} Let $\mathcal{S}$ be a $\mathscr{W}^{2,n}$-set with associated pair 
			$(\mathcal{S}', F)$ and assume that, in an open neighborhood 
			$U \subseteq \mathbf{R}^{n+1}$ of $z \in \mathcal{S}'$, there exists 
			a family $\{\Gamma_1,\ldots,\Gamma_{q(z)}\}$ of subsets of $\mathbf{R}^{n+1}$ 
			such that
			\begin{equation}\label{Ufset1.5}
				\mathcal{S}' \cap U = \bigcup_{i=1}^{q(z)} (\Gamma_i \cap U),
			\end{equation}
			where every $\Gamma_i \cap U$ coincides with the graph of a 
			$(C^0 \cap W^{2,n})$-function. Then$:$
			\begin{enumerate}[label=(\roman*)]
				\item $K\cap\mathcal{S}$ is $\mathcal{H}^n$-rectifiable of class $2$ for every compact set $K\subset\mathbf{R}^{n+1}$\,$;$
				\item $\textup{nor}(\mathcal{S})$ satisfies the Lusin $(N)$-property, namely
				$$\smash{Z\subset\mathcal{S} \ \ \textup{\emph{s.t.}} \ \  \mathcal{H}^n(Z)=0 \ \ \Rightarrow \ \ \mathcal{H}^n\big(\textup{nor}(\mathcal{S})\restrict Z\big)=0 \, ;}$$
				\item If $(w, \nu) \in \textup{nor}(\mathcal{S})$, then $ \textup{Tan}(\mathcal{S}, w) =  \nu^\perp\,;$
				\item $N_\infty(\mathcal{S})=\emptyset \, ;$
				\item $\mathcal{H}^n\big(N_1(\mathcal{S})\big)=0 \, $;
				\item $\mathcal{H}^n\Big(\big(\mathcal{S}'\cap U\big)\setminus\bigcup_{i=1}^{q(z)}\big(N_2(\Gamma_i)\cap U\big)\Big)=0 \,;$
				\item for each $i \in \{1, \ldots, q(z)\}$, let 
				$\nu_i : N_2(\Gamma_i) \cap U \to \mathbf{S}^n$ be such that
				$$\textup{nor}(\Gamma_i\,,x)=\big\{\nu_i(x),-\nu_i(x)\big\} \quad\textup{\emph{for $x\in N_2(\Gamma_i)\cap U$}}.$$
				Then, for $\mathcal{H}^n$-a.e. $x\in\mathcal{S'}\cap U$, we have
				\begin{equation*}\label{prop1}
					\nu_i(x)=\pm\nu_j(x)\quad \textup{\emph{and}} \quad B^{n+1}_s\big(x\pm\nu_i(x)\big)\cap\mathcal{S'}=\emptyset\textup{\emph{ for some $s>0$}}
				\end{equation*}
				whenever $i,j\in\pmb{\mathcal{I}}(x)$. Moreover, for $\mathcal{H}^n$-a.e. $(x,u)\in\textup{nor}(\mathcal{S}')\restrict U$, we have
				\begin{equation}\label{Ufset32}
					\sum_{i=1}^{q(z)}\pmb{1}_{\textup{nor}(\Gamma_i)}(x,u)=\pmb\bigiotab(x) \ \pmb{1}_{\textup{nor}(\mathcal{S}')}(x,u)
				\end{equation}
				therefore
				\begin{equation}\label{Ufset34}
					\mathcal{H}^n\Big(\Big[\bigcup_{i=1}^{q(z)}\textup{nor}(\Gamma_i)\restrict U\Big]\setminus\big[\textup{nor}(\mathcal{S}')\restrict U\big]\Big)=0\,;
				\end{equation}
				\item  $\mathcal{H}^n\big(\mathcal{S} \setminus N_2(\mathcal{S})\big) =0\,.$
			\end{enumerate}
		\end{lemma}
		
		\begin{proof}
			Since $F(\mathcal{S}')=\mathcal{S}$, the assertion \emph{(i)} follows from (\ref{Ufset1.5}) and Remark \ref{Ufset3}.
			
			To prove \emph{(ii)}, we consider the $\smash{C^1}$-diffeomorphism 
			\begin{equation*}
				\Psi_F:(x,y)\in\mathbf{R}^{n+1}\times\mathbf{S}^{n}\mapsto \bigg(F(x), \frac{(D F(x)^{-1})^\ast(y)}{| (D F(x)^{-1})^\ast(y) |}\bigg)\in\mathbf{R}^{n+1}\times\mathbf{S}^n
			\end{equation*}
			for which we have (cf. \cite[Lemma 2.1]{SantilliRectVarCurv})
			\begin{equation*}
				\Psi_F\big(\textup{nor}(\mathcal{S}')\big) = \textup{nor}(\mathcal{S}) \, .
			\end{equation*}
			Hence, for every set $Z\subseteq\mathbf{R}^{n+1}$ we deduce that
			\begin{align*}\mathcal{H}^n\big(\textup{nor}(\mathcal{S})\restrict Z\big)=\mathcal{H}^n\Big(\Psi_F\big(\textup{nor}(\mathcal{S}')\restrict F^{-1}(Z)\big)\Big)
			\end{align*}then, to prove the Lusin $(N)$-property on $\textup{nor}(\mathcal{S})$, it is enough to prove the same property on $\textup{nor}(\mathcal{S}')$. To this aim, let us consider an arbitrary $z\in\mathcal{S'}$ and we choose an open neighborhood $ U \subseteq \mathbf{R}^{n+1} $ of $z$ such that (\ref{Ufset1.5}) holds, then we readily infer that
			\begin{equation}\label{Ufset33}
				\textup{nor}(\mathcal{S}')\restrict U\subseteq\bigcup_{i=1}^{q(z)}\textup{nor}(\Gamma_i)\restrict U\,.
			\end{equation}Hence, from (\ref{Ufset33}) and (\ref{LusinNorGraph}), we deduce that
			\begin{equation*}\label{preLusin}
				\mathcal{H}^n\big(\textup{nor}(\mathcal{S}')\restrict (U\cap Z)\big)\leq\sum_{i=1}^{q(z)}\mathcal{H}^n\big(\textup{nor}(\Gamma_i)\restrict (U\cap Z)\big)=0
			\end{equation*}
			for every $\mathcal{H}^n$-negligible set $Z\subset\mathcal{S}$. Since $ U $ is arbitrarily chosen, the Lusin $(N)$-property on $\textup{nor}(\mathcal{S'})$ easily follows.
			
			Now we prove \emph{(iii)}, hence we choose $(w, \nu) \in \textup{nor}(\mathcal{S}) $, where $ w = F(z) $. There exist an open neighborhood $U$ of $z$ and  $\smash{\{\Gamma_i\}_{i=1}^{q(z)}}$ such that (\ref{Ufset1.5}) holds and
			$$	\mathcal{S}\cap F(U)=\bigcup_{i=1}^{q(z)}\big(F(\Gamma_i)\cap F(U)\big)\,,$$
			where every $\Gamma_i\cap U$ is the graph of a $(C^0\cap W^{2,n})$-function. We claim that, if $i\in\{1,\ldots,q(z)\}$ such that $ z \in \Gamma_i\cap U$, then $ \textup{Tan}\big(F(\Gamma_i), w\big) = \nu^\perp $. This clearly proves \emph{(iii)}, since (cf. (\ref{tanUnion}))
			$$ \textup{Tan}(\mathcal{S}, w) = \bigcup_{i=1}^{q(z)} \textup{Tan}\big(F(\Gamma_i), w\big) $$
			and $ \textup{Tan}\big(F(\Gamma_i), w\big) = \emptyset $ if $ w \notin F(\Gamma_i)$ (indeed $ w \in F(U)$ and $F(\Gamma_i)$ is a closed set in $F(U)$). If $ z \in \Gamma_i \cap U $, then there exists $ r > 0 $ such that
			\begin{equation*}
				B^{n+1}_{r}\big(w+ r\nu\big)\cap\big(F(\Gamma_i)\cap F(U)\big)=\emptyset\,.
			\end{equation*}	
			The domain $ \Omega : = F^{-1}\big[ B^{n+1}_{r}\big(w+ r\nu\big)  \big] $ is $ C^2 $-regular, $ z \in \partial \Omega$ and $ \textup{Tan}(\partial \Omega, z) = DF^{-1}(z)(\nu^\perp)$ by \cite[$3.1.21$]{Fed69}. Since $ \Omega \cap \Gamma_i \cap U = \emptyset$,
			it follows from Lemma \ref{graphProperties} \emph{(i)} that 
			$$\textup{Tan}(\Gamma_i, z) = \textup{Tan}(\partial \Omega, z )\,,$$
			namely $ \textup{Tan}(\Gamma_i, z) = DF^{-1}(z)(\nu^\perp) $. Again by \cite[$3.1.21$]{Fed69}, we infer that
			$$ \textup{Tan}\big(F(\Gamma_i), w\big) = D F(z)\big[ \textup{Tan}(\Gamma_i, z)\big] = \nu^\perp\,. $$
			
			Clearly \emph{(iv)} follows from \emph{(iii)}.
			
			Now we prove $\emph{(v)}$ and $\emph{(vi)}$. To show that $ \mathcal{H}^n(N_1(\mathcal{S})) =0 $, since $ N_1(\mathcal{S}) = F(N_1(\mathcal{S}')) $ (cf. Lemma \ref{LemmaSelection} \emph{(ii)}), it is enough to prove that $ \mathcal{H}^n(N_1(\mathcal{S}')) =0 $. Let us consider an arbitrary $z\in\mathcal{S'}$ and we choose a bounded open neighborhood $ U \subset \mathbf{R}^{n+1} $ of $z$ such that 
			\begin{equation*}\label{Ufset2}
				\mathcal{S}'\cap U=\bigcup_{i=1}^{q(z)}\big(\Gamma_i\cap U\big)\,,
			\end{equation*} 
			where every $$\Gamma_i\cap U =\overline{f_i}(V_i)+p_i$$
			for some $V_i$ open in $\eta_i^\perp$ with $ 0 \in V_i $, $p_i\in\Gamma_i\cap U$, $f_i\in C^0(V_i)\cap W^{2,n}(V_i)$ with $ f_i(0) =0 $ and
			$$\overline{f_i}:x\in V_i\mapsto x+f_i(x)\eta_i\in\mathbf{R}^{n+1}\,.$$
			Therefore, for every $i\in \{1,\ldots,q(z)\}$, by Lemma \ref{graphProperties} \emph{(ii)} we have that
			\begin{equation*}\label{Ufset19} N_1(\Gamma_i ) \cap U \subseteq\overline{f_i}\Big(\big(\mathcal{S}^\ast(f_i)\cup\mathcal{S}_\ast(f_i)\big)\setminus\big(\mathcal{S}^\ast(f_i)\cap\mathcal{S}_\ast(f_i)\big)\Big)+p_i
			\end{equation*}
			and by Theorem \ref{W2n twice diff} and Remark \ref{rmk Lusin condition (N) for f} we infer
			$$\mathcal{H}^n\big(N_1(\Gamma_i ) \cap U \big)\leq\mathcal{H}^n\bigg(\overline{f_i}\Big(\big(\mathcal{S}^\ast(f_i)\cup\mathcal{S}_\ast(f_i)\big)\setminus\big(\mathcal{S}^\ast(f_i)\cap\mathcal{S}_\ast(f_i)\big)\Big)\bigg)=0\,.$$
			Noting that $\smash{N_1(\mathcal{S}')\cap U\subseteq \textstyle\bigcup_{i=1}^{q(z)}N_1(\Gamma_i) \cap U }$, we obtain that $\smash{\mathcal{H}^n\big(N_1(\mathcal{S}')\cap U\big)=0}$. Since $ U $ is arbitrarily chosen we conclude that $ \mathcal{H}^n(N_1(\mathcal{S}')) =0 $. To prove \emph{(vi)}, applying Lemma \ref{graphProperties} \emph{(ii)} we deduce that
			{\allowdisplaybreaks\begin{align*}
					(\mathcal{S}'\cap U)\setminus\bigcup_{i=1}^{q(z)}\big(N_2(\Gamma_i)\cap U\big)&\subseteq\bigcup_{i=1}^{q(z)}\Big[\big(\Gamma_i\setminus N_2(\Gamma_i)\big)\cap U\Big]&\\
					&=\bigcup_{i=1}^{q(z)}\Big[\,\overline{f}_i\Big(V_i\setminus\big(\mathcal{S}^\ast(f_i)\cap\mathcal{S}_\ast(f_i)\big)\Big)+p_i\Big]
			\end{align*}}then, by Theorem \ref{W2n twice diff} and  Remark \ref{rmk Lusin condition (N) for f}, we obtain the desired result.
			
			To prove \emph{(vii)}, first we recall that the set $\smash{\textstyle{\bigcup_{i=1}^{q(z)}}(N_2(\Gamma_i)\cap U)}$ has full $\mathcal{H}^n$-measure in $\mathcal{ S}'\cap U$ (cf. statement \emph{(vi)}). Then, for $\mathcal{H}^n$-a.e. $x\in N_2(\Gamma_i)\cap N_2(\Gamma_j)\cap U$ where $i,j\in\{1,\ldots,q(z)\}$, by Lemma \ref{graphProperties} \emph{(i)} (cf. (\ref{tangent})) and by the locality property of the approximate tangent spaces (cf. \cite[3.2.16]{Fed69}) we have that
			\begin{align*}
				\quad\quad\quad\quad\quad\quad\nu_i(x)^\perp&=\textup{Tan}(\Gamma_i\,,x)=\textup{Tan}^n(\mathcal{H}^n\restrict\Gamma_i\,,x)&\\
				&=\textup{Tan}^n(\mathcal{H}^n\restrict\Gamma_j\,,x)=\textup{Tan}(\Gamma_j\,,x)=\nu_j(x)^\perp\,,
			\end{align*}
			hence
			\begin{equation*}
				\nu_i(x)=\pm\nu_j(x) \quad\textup{for $\mathcal{H}^n$-a.e. $x\in N_2(\Gamma_i)\cap N_2(\Gamma_j)\cap U$}\,.
			\end{equation*}
			So there exists a map $\nu$, with values in $\mathbf{S}^n$ and defined $\mathcal{H}^n$-a.e. on $\mathcal{S'}\cap U$ in such a way that
			$$\nu(x)\in\{\nu_i(x),-\nu_i(x)\} \quad\textup{if $x\in N_2(\Gamma_i)\cap U$}$$
			for $i\in\{1,\ldots,q(z)\}$, such that for $\mathcal{H}^n$-a.e. $x\in\mathcal{S'}\cap U$ we have
			$$\smash{B^{n+1}_s\big(x\pm s\nu(x)\big)\cap\mathcal{ S}'=\emptyset} \quad\textup{for some $s>0$}$$namely (cf. statement \emph{(iv)})
			\begin{equation}\label{N2property}
				\smash{\mathcal{H}^n\big((\mathcal{S}'\cap U)\setminus N_2(\mathcal{S}')\big)=0\,.}
			\end{equation}
			Henceforth, for every $i\in\{1,\ldots,q(z)\}$, we have that
			$$\textup{nor}(\mathcal{S'},x)=\textup{nor}(\Gamma_i\,,x)=\{\nu(x),-\nu(x)\} \quad\textup{for $\mathcal{H}^n$-a.e. $x\in\Gamma_i\cap U$}$$
			thus, by the Lusin $(N)$-property on $\textup{nor}(\mathcal{S'})$
			$$\pmb{1}_{\textup{nor}(\Gamma_i)}(x,u)=\pmb{1}_{\textup{nor}(\mathcal{S}')}(x,u)\,\pmb{1}_{\Gamma_i}(x) \quad\textup{for $\mathcal{H}^n$-a.e. $(x,u)\in\textup{nor}(\mathcal{S}')\restrict U$}\,.$$
			Therefore, from the definition of $\pmb\bigiotab$ (cf.\,(\ref{Ufset7})), we conclude that 
			$$\sum_{i=1}^{q(z)}\pmb{1}_{\textup{nor}(\Gamma_i)
			}(x,u)=\pmb\bigiotab(x) \ \pmb{1}_{\textup{nor}(\mathcal{S}')
			}(x,u) \quad\textup{for $\mathcal{H}^n$-a.e. $(x,u)\in\textup{nor}(\mathcal{S}')\restrict U$}$$
			hence
			$$\mathcal{H}^n\Big(\Big[\bigcup_{i=1}^{q(z)}\textup{nor}(\Gamma_i)\restrict U\Big]\setminus\big[\textup{nor}(\mathcal{S}')\restrict U\big]\Big)=0\,.$$    
			
			To prove \emph{(viii)}, as usually it is sufficient to show the assertion on $\mathcal{S'}$ (cf.\,Lemma \ref{LemmaSelection}\,\emph{(ii)}) and so we conclude from (\ref{N2property}).The proof is complete.\end{proof}
		
		\subsection{Curvatures on $\mathscr{W}^{2,n}$-sets}\label{curvaturenotions}
		Let $\mathcal{S}$ be a $\mathscr{W}^{2,n}$-set. Since $N_2(\mathcal{S})$ 
		has full $\mathcal{H}^n$-measure in $\mathcal{S}$, 
		Lemma~\ref{LemmaSelection}\,\emph{(iii)} implies that the multivalued 
		function
		$$
		\textup{nor}(\mathcal{S}, \pmb{\cdot}) : \mathcal{S} \to \mathcal{P}(\mathbf{S}^n)
		$$
		admits an $\mathcal{H}^n$-measurable selection 
		$\nu_\mathcal{S} : \mathcal{S} \to \mathbf{S}^n$, moreover 
		(cf.\,\eqref{propertyselected})
		$$
		\nu_\mathcal{S}(p) \in \textup{Nor}^n(\mathcal{H}^n \restrict \mathcal{S}, p) 
		\quad \text{for $\mathcal{H}^n$-a.e.\ $p \in \mathcal{S}$}.
		$$
		Applying Lemma~\ref{lem approx diff unit normal}, we infer that 
		$\nu_\mathcal{S}$ is $(\mathcal{H}^n \restrict \mathcal{S})$-approximately 
		differentiable at $\mathcal{H}^n$-a.e.\ $p \in \mathcal{S}$ and 
		$\textup{ap}\,D\nu_\mathcal{S}(p)$ is a symmetric endomorphism of 
		$\textup{Tan}^n(\mathcal{H}^n \restrict \mathcal{S}, p)$. We refer to 
		$\nu_\mathcal{S}$ as a \emph{selected unit-normal vector field} on 
		$\mathcal{S}$.
		\begin{definition}[Approximate principal curvatures]
			\emph{Given a $\mathscr{W}^{2,n}$-set $\mathcal{S}$ and a selected unit-normal 
				vector field $\nu_\mathcal{S}$ on $\mathcal{S}$, the approximate 
				principal curvatures of $\mathcal{S}$ with respect to $\nu_\mathcal{S}$ 
				are the $(\mathcal{H}^n \restrict \mathcal{S})$-measurable maps
				$$
				\rchi_{\mathcal{S},1}, \ldots, \rchi_{\mathcal{S},n}
				$$
				defined so that $\rchi_{\mathcal{S},1}(p) \leq \cdots \leq 
				\rchi_{\mathcal{S},n}(p)$ are the eigenvalues of 
				$\textup{ap}\,D\nu_\mathcal{S}(p)$ for $\mathcal{H}^n$-a.e.\ $p \in \mathcal{S}$.}
		\end{definition}
		\begin{definition}[$ r $-th elementary symmetric function]\label{symmetric function}
			\emph{Let $ r \in \{1, \ldots , n\} $. The $r$-th elementary symmetric function $ \sigma_r : \mathbf{R}^n \rightarrow \mathbf{R} $ is  defined by 
				$$ \sigma_r(t_1, \ldots , t_n) :=  \frac{1}{{n \choose r}}\sum_{\lambda \in \Lambda_{n, r}} t_{\lambda(1)}\cdots \,t_{\lambda(r)}\quad\text{for every }(t_1, \ldots , t_n) \in \mathbf{R}^n\,, $$
				where $ \Lambda_{n,r} $ is the set of all increasing functions from $ \{1, \ldots , r\} $ to $ \{1, \ldots , n\} $. We also set 
				$$ \sigma_{0}(t_1, \ldots , t_n) := 1 \quad \text{for every $ (t_1, \ldots , t_n) \in \mathbf{R}^n $\,.} $$}
		\end{definition}
		\begin{definition}[$k$-th mean curvature function]\label{rth mean curvature function3}
			\emph{Given a $\mathscr{W}^{2,n}$-set $\mathcal{S}$, a selected unit-normal 
				vector field $\nu_\mathcal{S}$ on $\mathcal{S}$, and $k \in \{0, \ldots, n\}$, 
				the $k$-th mean curvature function of $\mathcal{S}$ with respect 
				to $\nu_\mathcal{S}$ is
				\begin{align*}
					H_{\mathcal{S},k}(p) &:= \sigma_k\big(\rchi_{\mathcal{S},1}(p), 
					\ldots, \rchi_{\mathcal{S},n}(p)\big) \\
					&\,= \frac{1}{\binom{n}{k}} \sum_{\lambda \in \Lambda(n,k)} 
					\rchi_{\mathcal{S},\lambda(1)}(p) \cdots \rchi_{\mathcal{S},\lambda(k)}(p)
					\quad \text{for $\mathcal{H}^n$-a.e.\ $p \in \mathcal{S}$}\,.
			\end{align*}}
		\end{definition}
		
		\begin{definition}\label{A_k}
			\emph{Given a $\mathscr{W}^{2,n}$-set $\mathcal{S}$ with associated pair 
				$(\mathcal{S}', F)$, where
				$$
				\mathcal{S}' = \{\pmb{\iota} > 0\} \quad \text{for some } 
				\pmb{\iota} \in F_nW^{2,n}(\mathbf{R}^{n+1})\,,
				$$
				and a selected unit-normal vector field $\nu_\mathcal{S}$ on $\mathcal{S}$, 
				we define
				$$
				\mathcal{A}_k(\mathcal{S}) := \int_{\mathcal{S}} H_{\mathcal{S},k}(x)\, 
				\pmb{\iota}\big(F^{-1}(x)\big)\, d\mathcal{H}^n(x) 
				\quad \text{for } k \in \{0, \ldots, n\}\,.
				$$}
		\end{definition}
		
		\begin{lemma}\label{Ufset15}
			Let $\mathcal{S}$ be a compact $\mathscr{W}^{2,n}$-set with associated 
			pair $(\mathcal{S}', F)$, and let $\nu_\mathcal{S}$ be a selected 
			unit-normal vector field on $\mathcal{S}$. Then the following statements hold:
			\begin{enumerate}[label=(\roman*)]
				\item $\mathcal{H}^n\big(\overline{\nu}_\mathcal{S}(Z)\big) = 0$ 
				for any $Z \subseteq N_2(\mathcal{S})$ with $\mathcal{H}^n(Z) = 0\,;$
				\item $\mathcal{H}^n\big(\textup{nor}(\mathcal{S}) \setminus 
				\big[\overline{\nu}_\mathcal{S}\big(N_2(\mathcal{S})\big) \cup 
				\overline{-\nu}_\mathcal{S}\big(N_2(\mathcal{S})\big)\big]\big) = 0$ 
				and $\mathcal{H}^n\big(\textup{nor}(\mathcal{S})\big) < \infty\,;$
				\item for any $i \in \{1, \ldots, n\}$ (cf.\,Definition~\ref{curv}),
				\begin{equation}\label{pmEigen}
					\rchi_{\mathcal{S},i}(x) = \kappa_{\mathcal{S},i}
					\big(x, \nu_\mathcal{S}(x)\big) = 
					-\kappa_{\mathcal{S},i}\big(x, -\nu_\mathcal{S}(x)\big) 
					\quad \text{for $\mathcal{H}^n$-a.e.\ $x \in \mathcal{S}$}\,.
				\end{equation}
				In particular, $\kappa_{\mathcal{S},i}(x,u) < +\infty$ for 
				$\mathcal{H}^n$-a.e.\ $(x,u) \in \textup{nor}(\mathcal{S})$.
			\end{enumerate}
		\end{lemma}
		
		\begin{proof} Assertion \emph{(i)} follows readily from the Lusin $(N)$-condition 
			on $\textup{nor}(\mathcal{S})$ (cf.~Lemma~\ref{Ufset14}\,\emph{(ii)}).
			
			To prove \emph{(ii)}, we recall the definition of $\Psi_F$ from 
			\eqref{PsiF} and notice that (cf.\,\cite[Lemma~2.1]{SantilliRectVarCurv})
			$$
			\Psi_F\big(\textup{nor}(\mathcal{S}')\big) = \textup{nor}(\mathcal{S}),
			$$
			hence we deduce that $\mathcal{H}^n\big(\textup{nor}(\mathcal{S})\big) < \infty$ 
			from \eqref{Ufset1} and \eqref{LusinNorGraph}. Moreover, since 
			$N_2(\mathcal{S})$ has full $\mathcal{H}^n$-measure in $\mathcal{S}$ 
			(cf.\,Lemma~\ref{Ufset14}\,\emph{(vii)}), by the Lusin $(N)$-property 
			on $\textup{nor}(\mathcal{S})$ we obtain
			$$
			\mathcal{H}^n\big(\textup{nor}(\mathcal{S}) \setminus 
			\big[\overline{\nu}_\mathcal{S}\big(N_2(\mathcal{S})\big) \cup 
			\overline{-\nu}_\mathcal{S}\big(N_2(\mathcal{S})\big)\big]\big) = 
			\mathcal{H}^n\big(\textup{nor}(\mathcal{S}) \setminus 
			\textup{nor}(\mathcal{S}) \restrict N_2(\mathcal{S})\big) = 0\,.
			$$
			
			To prove \emph{(iii)}, we first apply Lemma~\ref{lem approx diff unit normal} 
			to find a countable family of $\mathcal{H}^n$-measurable sets 
			$X_i \subseteq N_2(\mathcal{S})$ with 
			$\mathcal{H}^n\big(\mathcal{S} \setminus \bigcup_{i=1}^\infty X_i\big) = 0$ 
			and $\textsf{Lip}(\nu_\mathcal{S}|X_i) < \infty$ for every $i \in \mathbf{N}$. 
			We then define $Y_i$ as the set of points $x \in X_i$ satisfying all of 
			the following:
			\begin{itemize}
				\item[$\medmath\bullet$] $\nu_\mathcal{S}$ is $(\mathcal{H}^n \restrict \mathcal{S})$-approximately 
				differentiable at $x$\,;
				\item[$\medmath\bullet$] $\Theta^n\big(\mathcal{H}^n \restrict (\mathcal{S} \setminus X_i), x\big) = 0$ 
				(this holds $\mathcal{H}^n$-a.e.\ on $X_i$, cf.~\cite[2.10.19\,(4)]{Fed69})\,;
				\item[$\medmath\bullet$] $\textup{Tan}^n(\mathcal{H}^n \restrict \mathcal{S}, x) \in \mathbf{G}(n+1,n)$\,;
				\item[$\medmath\bullet$] $\textup{Tan}^n\big(\mathcal{H}^n \restrict \textup{nor}(\mathcal{S}), 
				\overline{\nu}_\mathcal{S}(x)\big) = \textup{Tan}^n\big(\mathcal{H}^n \restrict 
				\overline{\nu}_\mathcal{S}\big(N_2(\mathcal{S})\big), 
				\overline{\nu}_\mathcal{S}(x)\big) \in \mathbf{G}(2n+2,n)$\,;
				\item[$\medmath\bullet$] $\textup{Tan}^n\big(\mathcal{H}^n \restrict \textup{nor}(\mathcal{S}), 
				\overline{-\nu}_\mathcal{S}(x)\big) = \textup{Tan}^n\big(\mathcal{H}^n \restrict 
				\overline{-\nu}_\mathcal{S}\big(N_2(\mathcal{S})\big), 
				\overline{-\nu}_\mathcal{S}(x)\big) \in \mathbf{G}(2n+2,n)$\,.
			\end{itemize}
			Since $ \overline{\nu}_\mathcal{S}| X_i $ and $ \overline{\medmath{-}\nu}_\mathcal{S}| X_i $ are bi-lipschitz and $\smash{\textup{Tan}^n\big(\mathcal{H}^n\restrict \overline{\medmath{\pm}\nu}_\mathcal{S}\big(N_2(\mathcal{S})\big), (x,u)\big)}\in\mathbf{G}(2n+2,n)$
			for $ \mathcal{H}^n $-a.e.$\, (x,u) \in \overline{\medmath{\pm}\nu}_\mathcal{S}\big(N_2(\mathcal{S})\big)$ (cf.\,\cite[Theorem 3.2.19]{Fed69}), by the locality property of approximate tangent spaces (cf.\,\cite[$3.2.16$ p.252]{Fed69}) we infer
			\begin{equation*}
				\mathcal{H}^n(X_i \setminus Y_i) =0 \quad \textrm{for every $ i\in\mathbf{N}\,. $}
			\end{equation*} 
			It follows from \emph{(i)}  and \emph{(ii)} that 
			{\allowdisplaybreaks\begin{align}\label{W2n domains eq2.5}\nonumber
					&\mathcal{H}^n\Big(\textup{nor}(\mathcal{S}) \setminus \bigcup_{i=1}^\infty \big[\overline{\nu}_\mathcal{S}(Y_i)\cup\overline{\medmath{-}\nu}_\mathcal{S}(Y_i)\big]\Big)&\\ \nonumber
					&\quad=\mathcal{H}^n\Big(\big[\overline{\nu}_\mathcal{S}\big(N_2(\mathcal{S})\big)\cup\overline{\medmath{-}\nu}_\mathcal{S}\big(N_2(\mathcal{S})\big)\big] \setminus \bigcup_{j=1}^\infty \big[\overline{\nu}_\mathcal{S}(Y_j)\cup\overline{\medmath{-}\nu}_\mathcal{S}(Y_j)\big]\Big)&\\ \nonumber
					&\quad=\mathcal{H}^n\Big(\bigcup_{i=1}^\infty\big[\overline{\nu}_\mathcal{S}(X_i)\cup\overline{\medmath{-}\nu}_\mathcal{S}(X_i)\big] \setminus \bigcup_{j=1}^\infty \big[\overline{\nu}_\mathcal{S}(Y_j)\cup\overline{\medmath{-}\nu}_\mathcal{S}(Y_j)\big]\Big)&\\ \nonumber
					&\quad\leq\sum_{i=1}^\infty\mathcal{H}^n\Big(\big[\overline{\nu}_\mathcal{S}(X_i)\cup\overline{\medmath{-}\nu}_\mathcal{S}(X_i)\big] \setminus  \big[\overline{\nu}_\mathcal{S}(Y_i)\cup\overline{\medmath{-}\nu}_\mathcal{S}(Y_i)\big]\Big)&\\
					&\quad\leq\sum_{i=1}^\infty\mathcal{H}^n\big(\overline{\nu}_\mathcal{S}(X_i) \setminus  \overline{\nu}_\mathcal{S}(Y_i)\Big)+\sum_{i=1}^\infty\mathcal{H}^n\big(\overline{\medmath{-}\nu}_\mathcal{S}(X_i) \setminus  \overline{\medmath{-}\nu}_\mathcal{S}(Y_i)\big)=0\,.
			\end{align}}
			In addition
			{\allowdisplaybreaks\begin{align}\label{W2n domains eq3}\mathcal{H}^n\Big(\mathcal{S}\setminus\bigcup_{i=1}^\infty Y_i\Big)\leq\sum_{j=1}^\infty\mathcal{H}^n(X_j\setminus Y_j)=0\,.
			\end{align}}
			Fix $x \in Y_i$. Then there exists a map $ g : \mathbf{R}^{n+1} \rightarrow \mathbf{R}^{n+1}$ pointwise differentiable at $ x $ such that $g(x)=\nu_\mathcal{ S}(x)$, $\Theta^n(\mathcal{H}^n \restrict \mathcal{S} \setminus \{g = \nu_\mathcal{S}\}, x)=0$, and
			$$\textup{ap}\,D \overline{\nu}_\mathcal{S}(x) := D \overline{g}(x) |  \textup{Tan}^n(\mathcal{H}^n\restrict \mathcal{S}, x)\,,\quad\textup{ap}\,D \overline{\medmath{-}\nu}_\mathcal{S}(x) := D \overline{\medmath{-}g}(x) |  \textup{Tan}^n(\mathcal{H}^n\restrict \mathcal{S}, x)\,.$$
			Since $ \overline{\medmath{}g} | X_i \cap \{\overline{\medmath{}g} = \overline{\medmath{}\nu}_\mathcal{S}\} $ and $ \overline{\medmath{-}g} | X_i \cap \{\overline{\medmath{-}g} = \overline{\medmath{-}\nu}_\mathcal{S}\} $ are bi-lipschitz, $ \textup{ap}\,D \overline{\medmath{}\nu}_\mathcal{S}(x) $ and $ \textup{ap}\,D \overline{\medmath{-}\nu}_\mathcal{S}(x) $ are injective, and  
			$$ \textup{Tan}^n(\mathcal{H}^n \restrict \mathcal{S}, x) = \textup{Tan}^n\big(\mathcal{H}^n \restrict X_i \cap \{\overline{\medmath{}g} = \overline{\medmath{}\nu}_\mathcal{S}\}, x\big)\,, $$
			we infer from \cite[Lemma B.2]{SantilliAnnali} that
			{\allowdisplaybreaks\begin{align*}
					\textup{ap}\,D \overline{\medmath{}\nu}_\mathcal{S}(x)\big[\textup{Tan}^n(\mathcal{H}^n\restrict \mathcal{S}, x)\big]&=D \overline{\medmath{}g}(x)\big[\textup{Tan}^n(\mathcal{H}^n\restrict \mathcal{S}, x)\big]&\\
					&=D \overline{\medmath{}g}(x)\big[\textup{Tan}^n\big(\mathcal{H}^n \restrict X_i \cap \{\overline{\medmath{}g} = \overline{\medmath{}\nu}_\mathcal{S}\}, x\big)\big]&\\
					&\subseteq\textup{Tan}^n\big(\mathcal{H}^n\restrict\overline{\medmath{}g}(X_i \cap \{\overline{\medmath{}g} = \overline{\medmath{}\nu}_\mathcal{S}\}),\overline{\medmath{}g}(x)\big)&\\
					&\subseteq \textup{Tan}^n\big(\mathcal{H}^n\restrict \overline{\medmath{}\nu}_\mathcal{S}\big(N_2(\mathcal{S})\big), \overline{\medmath{}\nu}_\mathcal{S}(x)\big)&\\
					&=\textup{Tan}^n\big(\mathcal{H}^n\restrict \textup{nor}(\mathcal{S}), \overline{\medmath{}\nu}_\mathcal{S}(x)\big)\in\mathbf{G}(2n+2,n)\,,
			\end{align*}}
			and similarly
			$$\textup{ap}\,D \overline{\medmath{-}\nu}_\mathcal{S}(x)\big[\textup{Tan}^n(\mathcal{H}^n\restrict \mathcal{S}, x)\big]\subseteq\textup{Tan}^n\big(\mathcal{H}^n\restrict \textup{nor}(\mathcal{S}), \overline{\medmath{-}\nu}_\mathcal{S}(x)\big)\in\mathbf{G}(2n+2,n)\,.$$
			Since $\textup{Tan}^n(\mathcal{H}^n\restrict \mathcal{S}, x)\in\mathbf{G}(n+1,n)$ and $\textup{ap}\,D \overline{\medmath{\pm}\nu}_\mathcal{S}(x)$ are injective, we infer
			$$\textup{ap}\,D \overline{\medmath{}\nu}_\mathcal{S}(x)\big[\textup{Tan}^n(\mathcal{H}^n\restrict \mathcal{S}, x)\big]=\textup{Tan}^n\big(\mathcal{H}^n\restrict \textup{nor}(\mathcal{S}), \overline{\medmath{}\nu}_\mathcal{S}(x)\big)\,,$$
			$$\textup{ap}\,D \overline{\medmath{-}\nu}_\mathcal{S}(x)\big[\textup{Tan}^n(\mathcal{H}^n\restrict \mathcal{S}, x)\big]=\textup{Tan}^n\big(\mathcal{H}^n\restrict \textup{nor}(\mathcal{S}), \overline{\medmath{-}\nu}_\mathcal{S}(x)\big)\,.$$
			Hence, if $ \{\tau_1, \ldots , \tau_n\} $ is an orthonormal basis of $ \textup{Tan}^n(\mathcal{H}^n \restrict \mathcal{S}, x) $ with $$ \textup{ap}\,D \nu_\mathcal{S}(x)(\tau_i) = \rchi_{\mathcal{S}, i}(x) \tau_i \quad \textup{for $ i \in\{ 1, \ldots , n \}$}\,,$$we conclude that 
			$$ \Bigg\{\bigg(\frac{1}{\sqrt{1 + \rchi_{\mathcal{S}, i}(x)^2}}\tau_i\,, \frac{\rchi_{\mathcal{S}, i}(x)}{\sqrt{1 + \rchi_{\mathcal{S}, i}(x)^2}}\tau_i\bigg): i \in\{ 1, \ldots , n \}\Bigg\}\, , $$
			$$ \Bigg\{\bigg(\frac{1}{\sqrt{1 + \rchi_{\mathcal{S}, i}(x)^2}}\tau_i\,, \frac{-\rchi_{\mathcal{S}, i}(x)}{\sqrt{1 + \rchi_{\mathcal{S}, i}(x)^2}}\tau_i\bigg): i \in\{ 1, \ldots , n \}\Bigg\} $$
			are orthonormal basis of $ \textup{Tan}^n\big(\mathcal{H}^n \restrict \textup{nor}(\mathcal{S}), \overline{\nu}_\mathcal{S}(x)\big) $ and $ \textup{Tan}^n\big(\mathcal{H}^n \restrict \textup{nor}(\mathcal{S}), \overline{\medmath{-}\nu}_\mathcal{S}(x)\big) $, respectively. Since $x \in Y_i$ 
			was arbitrary, thanks to~\eqref{W2n domains eq3}, we deduce from the 
			uniqueness stated in Lemma~\ref{lem: Santilli20} 
			(cf.\,Definition~\ref{curv}) that
			$$
			\rchi_{\mathcal{S},i}(x) = \kappa_{\mathcal{S},i}\big(x,\nu_\mathcal{S}(x)\big) 
			= -\kappa_{\mathcal{S},i}\big(x,-\nu_\mathcal{S}(x)\big) 
			\quad \text{for $\mathcal{H}^n$-a.e.\ $x \in \mathcal{S}$}\,.
			$$
			From~\eqref{W2n domains eq2.5} it also follows that
			$$\kappa_{\mathcal{S},i}\big(x,u)\in\big\{\rchi_{\mathcal{S},i}(x),-\rchi_{\mathcal{S},i}(x)\big\} \quad \textrm{for $ \mathcal{H}^n $-a.e.\ $ (x,u) \in \textup{nor}(\mathcal{S})$}\,.$$
			The proof is complete.
		\end{proof}
		
		\begin{definition}
			\emph{Let $\mathcal{S}$ be a compact $\mathscr{W}^{2,n}$-set. We define
				$$\textup{nor}(\mathcal{S})^{(n)}:=\big\{(x,u)\in\textup{nor}(\mathcal{S}):\kappa_{\mathcal{S},n}(x,u)<+\infty\big\}\,.$$}
		\end{definition}
		
		\begin{remark} \textup{From Lemma \ref{Ufset15}\,\emph{(iii)}, we infer that
				\begin{equation}\label{n-curv2}
					\mathcal{H}^n\big(\textup{nor}(\mathcal{S})\setminus\textup{nor}(\mathcal{S})^{(n)}\big)=0\,.
				\end{equation}
				Furthermore, if $\pi_0:\mathcal{ S}\times\mathbf{S}^n\to\mathcal{S}$ is the canonical projection onto the first factor, then
				{\allowdisplaybreaks\begin{align}\label{Ufset25}\nonumber J^{\textup{nor}(\mathcal{ S})}_n\pi_0(x,u)&=\bigg|\Big[\bigwedge\nolimits_nD^{\textup{nor}(\mathcal{ S})}\pi_0(x,u)\Big]\bigg(\frac{\xi_{\mathcal{ S},1}(x,u)\wedge\cdots\wedge\xi_{\mathcal{ S},n}(x,u)}{\big|\xi_{\mathcal{ S},1}(x,u)\wedge\cdots\wedge\xi_{\mathcal{ S},n}(x,u)\big|}\bigg)\bigg|&\\ \nonumber
						&=\zeta_\mathcal{ S}(x,u) \, \big|\pi_0\big(\xi_{\mathcal{ S},1}(y,u)\big)\wedge\cdots\wedge\pi_0\big(\xi_{\mathcal{ S},n}(x,u)\big)\big|&\\ \nonumber
						&=\zeta_\mathcal{ S}(x,u)&\\
						&=\smash{\textstyle{\prod_{i=1}^{n}}\big(1 + \kappa_{\mathcal{S}, i}(x,u)^2\big)^{-\frac{1}{2}}>0 }
				\end{align}}
				for $\mathcal{H}^n$-a.e. $(x,u)\in\textup{nor}(\mathcal{ S})$, where $\zeta_\mathcal{ S}$ and $\{\xi_{\mathcal{ S},1},\ldots,\xi_{\mathcal{ S},n}\}$ are given in Definition \ref{defn-vect}.}
		\end{remark}
		The following definition is well posed.
		\begin{definition}[$k$-th mean curvature function of $\textup{nor}(\mathcal{S})$]\label{rth mean curvature function4} \emph{Given $\mathcal{S}$ a compact $\mathscr{W}^{2,n}$-set and $ k \in \{0, \ldots, n\} $, the $ k $-th mean curvature function of $\textup{nor}(\mathcal{S})$ is
				\begin{align*} H_{\textup{nor}(\mathcal{S}),k}(x,u)&:= \sigma_k\big(\kappa_{\mathcal{S}, 1}(x,u), \ldots , \kappa_{\mathcal{S}, n}(x,u)\big)&\\
					&\,\,=\frac{1}{{n \choose k}}\sum_{\lambda \in \Lambda(n,k)} \kappa_{\mathcal{S}, 1}(x,u)\cdots \kappa_{\mathcal{S},n}(x,u) \quad \textup{\emph{	for $ \mathcal{H}^n $-a.e.\ $ (x,u)\in\textup{nor}(\mathcal{S}) $\,,}}\end{align*}
				where $ \kappa_{\mathcal{S},1}, \ldots , \kappa_{\mathcal{S},n} $ are given in Definition $\ref{curv}$.}
		\end{definition} 
		
		\subsection{Reilly-type variational formulae}
		In this section let $\mathcal{S}$ be a compact $\mathscr{W}^{2,n}$-set with associated pair $(\mathcal{S}',F)$, that is $\mathcal{S}=F(\mathcal{S}')$ where $\mathcal{S}'=\{\pmb\bigiotab>0\}$ for some $\pmb\bigiotab\in F_nW^{2,n}(\mathbf{R}^{n+1})$ and $F$ is a $C^2$-diffeomorphism of $\mathbf{R}^{n+1}$. There exist a family $\smash{\{U_i\}_{i=1}^N}$ of bounded open neighborhoods of points $\smash{\{z_i\}_{i=1}^N\subset\mathcal{S}}$ such that
		\begin{equation}\label{Ufset39}
			\mathcal{S}'=\bigcup^N_{i=1}(\mathcal{S}'\cap U_i)=\bigcup^N_{i=1}\bigcup_{j=1}^{q(z_i)}(\Gamma^{(i)}_j\cap U_i) \ ,
		\end{equation}
		where every $\Gamma^{(i)}_j\cap U_i$ coincides with the graph of a $(C^0\cap W^{2,n})$-function (cf.\,Definition \ref{Ufset5}). Namely, for each $\smash{i\in\{1,\ldots,N\}}$ and for each $\smash{j\in\{1,\dots,q(z_i)\}}$ there exist $\smash{p^{(i)}_j\in\Gamma^{(i)}_j\cap U_i}$ and $\eta^{(i)}_j\in\mathbf{S}^n$ such that
		\begin{equation*}
			\Gamma^{(i)}_j\cap U_i=\overline{f}^{(i)}_j(V^{(i)}_j)+p^{(i)}_j , \end{equation*}for some sets $V^{(i)}_j$ open in $\smash{(\eta^{(i)}_j)^\perp}$ with $0\in V^{(i)}_j$, and some graph functions $\smash{f^{(i)}_j}$ on $V^{(i)}_j$.
		
		We associate with $\textup{nor}(\mathcal{S'})$ and $\textup{nor}(\mathcal{S})$ the $n$-vectorfields $\vv{\xi}_{\mathcal{S}'}$ and $\vv{\xi}_{\mathcal{S}}$, respectively, defined in according to Definition $\ref{defn-vect}$ and Lemma $\ref{lem: Santilli20}$. In particular, for $\mathcal{H}^n$-a.e. $(x,u)\in\textup{nor}(\mathcal{S'})$,
		\begin{equation*}
			\vv{\xi}_{\mathcal{S}'}(x,u):=\frac{\xi_{{\mathcal{S}'},1}(x,u)\wedge\cdots\wedge\xi_{{\mathcal{S}'},n}(x,u)}{|\xi_{{\mathcal{S}'},1}(x,u)\wedge\cdots\wedge\xi_{{\mathcal{S}'},n}(x,u)|}\in\bigwedge\nolimits_n(\mathbf{R}^{n+1}\times\mathbf{R}^{n+1})\,,
		\end{equation*}
		\begin{equation*}
			\zeta_\mathcal{S'}(x,u):=\frac{1}{|\xi_{\mathcal{S'},1}(x,u)\wedge\cdots\wedge\xi_{\mathcal{S'},n}(x,u)|}\in(0,+\infty)\,,
		\end{equation*}
		where each $\xi_{{\mathcal{S}'},i}$ is given by (notice that $\textup{nor}(\mathcal{S'})^{(n)}$ has full $\mathcal{H}^n$-measure in $\textup{nor}(\mathcal{S'})$; cf.\,(\ref{n-curv2}))
		\begin{equation*}\label{Ufset43}
			\xi_{{\mathcal{S}'},i}(x,u)=
			\big(\tau_{i}(x,u),\kappa_{{\mathcal{S}'},i}(x,u)\tau_{i}(x,u)\big)
		\end{equation*}
		for $\mathcal{H}^n$-a.e. $(x,u)\in\textup{nor}({\mathcal{S}'})$. The maps $ \{\tau_1, \ldots , \tau_n \}$ are defined $ \mathcal{H}^n $-a.e.\ on $ \textup{nor}(\mathcal{S}') $ in such a way that $ \{\tau_{1}(x,u), \ldots , \tau_n(x,u), u\} $ form a positively oriented orthonormal basis of $ \mathbf{R}^{n+1} $ for $\mathcal{H}^n$-a.e. $(x,u)\in\textup{nor}(\mathcal{S}')$ $($cf.\,Lemma $\ref{lem: Santilli20})$, i.e. $\smash{\langle\textstyle{\bigwedge_{i=1}^n}\tau_{i}(x,u)\wedge u,\pmb{e}'_1\wedge\cdots\wedge \pmb{e}'_{n+1}\rangle=1}$.
		
		\begin{lemma}\label{measurabilityXi}
			$\vv{\xi}_{\mathcal{S}'}$ is an $\big(\mathcal{H}^n\restrict\textup{nor}(\mathcal{S'})\big)$-measurable $n$-vectorfield.
		\end{lemma}
		
		\begin{proof} For every fixed $i\in\{1,\ldots,N\}$ and for $j\in\{1,\dots,q(z_i)\}$, we prove that 
			\begin{equation}\label{norBasis4}
				\smash{\vv{\xi}_{\mathcal{S}'}(y,u)=\vv{\eta}_{\Gamma^{(i)}_j\cap U_i}(y,u) \quad\textup{for $\mathcal{H}^n$-a.e. $(y,u)\in\textup{nor}(\Gamma^{(i)}_j)\restrict U_i$}}
			\end{equation}
			where the $n$-vectorfield on the right hand-side of~\eqref{norBasis4} is the Borel $n$-vectorfield given in (\ref{etaGraph}) (cf.\,Remark \ref{EpiAndGraph}).
			By the $\mathcal{H}^n$-rectifiability of the unit normal bundle, the locality property of the approximate tangent spaces, and (\ref{Ufset34}), we infer
			{\allowdisplaybreaks\begin{align*}
					&\smash{\quad\quad\textup{Tan}^n\big(\mathcal{H}^n\restrict\textup{nor}(\mathcal{S}'),(y,u)\big)}&\\
					&\quad\quad\quad\quad\quad=\textup{Tan}^n\big(\mathcal{H}^n\restrict\textup{nor}(\Gamma^{(i)}_j),(y,u)\big)\in\mathbf{G}(2n+2,n) \, ,
			\end{align*}}
			for $\mathcal{H}^n$-a.e. $\smash{(y,u)\in\textup{nor}(\Gamma^{(i)}_j)\restrict U_i}$. For such points $(y,u)$,
			{\allowdisplaybreaks\begin{align*}
					&\smash{\quad\quad\bigwedge\nolimits_n\big[\textup{Tan}^n\big(\mathcal{H}^n\restrict\textup{nor}(\mathcal{S}'),(y,u)\big)\big]}&\\
					&\quad\quad\quad\quad\quad=\bigwedge\nolimits_n\big[\textup{Tan}^n\big(\mathcal{H}^n\restrict\textup{nor}(\Gamma^{(i)}_j),(y,u)\big)\big]
			\end{align*}}
			and
			{\allowdisplaybreaks\begin{align*}
					&\smash{\quad\quad\dim\Big(\bigwedge\nolimits_n\big[\textup{Tan}^n\big(\mathcal{H}^n\restrict\textup{nor}(\mathcal{S}'),(y,u)\big)\big]\Big)}&\\
					&\quad \quad \quad\quad\quad =\dim\Big(\bigwedge\nolimits_n\big[\textup{Tan}^n\big(\mathcal{H}^n\restrict\textup{nor}(\Gamma^{(i)}_j),(y,u)\big)\big]\Big)=1 \, ,
			\end{align*}}
			therefore we readily deduce that
			\begin{equation}\label{norBasis3}
				\smash{\vv{\xi}_{\mathcal{S}'}(y,u)=\pm\vv{\eta}_{\Gamma^{(i)}_j\cap U_i}(y,u) \quad\textup{for $\mathcal{H}^n$-a.e. $(y,u)\in\textup{nor}(\Gamma^{(i)}_j)\restrict U_i$\,.}}
			\end{equation}
			Now we introduce the sets
			{\allowdisplaybreaks\begin{align*}
					&\smash{Z:=\big\{(y,u)\in \textup{nor}(\Gamma^{(i)}_j)\restrict U_i:\vv{\xi}_{\mathcal{S}'}(y,u)=-\vv{\eta}_{\Gamma^{(i)}_j\cap U_i}(y,u)\big\} \, ,}&\\
					&\smash{Z':=Z\cap\textup{nor}(\mathcal{S}')^{(n)} ,}
			\end{align*}}with the aim of proving that $\mathcal{H}^n(Z)=0$. Since $\mathcal{H}^n(Z\setminus Z')=0$ (cf.\,(\ref{n-curv2})\,and\,(\ref{Ufset34})), it suffices to show that $\mathcal{H}^n(Z')=0$. First, note that
			{\allowdisplaybreaks\begin{align*}\nonumber
					\big[\bigwedge\nolimits_n\pi_0\big]\big(\vv{\xi}_{\mathcal{S}'}(y,u)\big)&=\frac{\pi_0\big(\xi_{\mathcal{S}',1}(y,u)\big)\wedge\cdots\wedge\pi_0\big(\xi_{\mathcal{S}',n}(y,u)\big)}{|\xi_{\mathcal{S}',1}(y,u)\wedge\cdots\wedge\xi_{\mathcal{S}',n}(y,u)|}&\\
					&=\zeta_{\mathcal{S}'}(y,u) \pmb\cdot \tau_{1}(y,u)\wedge\cdots\wedge \tau_{n}(y,u) \quad\textup{for $\mathcal{H}^n$-a.e. $(y,u)\in Z'$}
			\end{align*}}and since $\smash{\{\tau_{1}(x,u),\dots,\tau_{n}(x,u),u\}}$ form a positively oriented orthonormal basis of $\mathbf{R}^{n+1}$, we obtain
			$$\smash{\langle\big[\bigwedge\nolimits_n\pi_0\big]\big(\vv{\xi}_{\mathcal{S}'}(y,u)\big)\wedge u, \pmb{e}_1\wedge\cdots\wedge \pmb{e}'_{n+1}\rangle=\zeta_{\mathcal{S}'}(y,u)>0 \quad \textup{for $\mathcal{H}^n$-a.e. $(y,u)\in Z'\,.$}}$$
			This gives a contradiction, from which the desired result follows. Indeed,
			$$\smash{\big[\bigwedge\nolimits_n\pi_0\big]\big(\vv{\xi}_{\mathcal{S}'}(y,u)\big)=-\big[\bigwedge\nolimits_n\pi_0\big]\big(\vv{\eta}_{\Gamma^{(i)}_j\cap U_i}(y,u)\big) \quad\textup{for any} \ (y,u)\in Z'}$$
			where (cf.\,(\ref{projectionproperty}))
			$$\smash{\langle\big[\bigwedge\nolimits_n\pi_0\big]\big(\vv{\eta}_{\Gamma^{(i)}_j\cap U_i}(y,u)\big)\wedge u, \pmb{e}'_1\wedge\cdots\wedge \pmb{e}'_n\rangle>0 \quad \textup{for $\mathcal{H}^n$-a.e. $(y,u)\in Z'$}} \, .$$
			The proof is complete.
		\end{proof}
		
		\begin{definition}\label{LegCycleS'}
			We define $\mathcal{N}_{\mathcal{S}'}\in\mathcal{D}_n(\mathbf{R}^{n+1}\times\mathbf{R}^{n+1})$ by
			\begin{equation}\label{Ufset42}
				\mathcal{N}_{\mathcal{S}'}:=(\pmb\bigiotab\circ\pi_0)\big(\mathcal{H}^n\restrict\textup{nor}(\mathcal{S}')\big)\wedge\vv{\xi}_{\mathcal{S}'} \, ,
			\end{equation}
			where $\pi_0:\mathbf{R}^{n+1}\times\mathbf{R}^{n+1}\to\mathbf{R}^{n+1}$ is the canonical projection on the first factor.\end{definition}
		
		\begin{theorem}\label{LegCycSset} $\mathcal{N}_{\mathcal{S}'}$ is a Legendrian cycle of $\mathbf{R}^{n+1}$, we call it as the Legendrian cycle associated with $\mathcal{S}'$. Moreover, for a selected unit-normal vector field $\nu_{\mathcal{S}'}$ on $\mathcal{S}'$, the following relation hold
			{\allowdisplaybreaks\begin{align}\label{W2n domains0.80}
					&\big(\mathcal{N}_{\mathcal{S}'} \restrict \varphi_{n-k}\big)(\phi)&\\ \nonumber
					&\quad\quad\quad= {n \choose k}\int_{\mathcal{S}'} \big[\phi\big(x, \nu_{\mathcal{S}'}(x)\big) +(-1)^k\,\phi\big(x,-\nu_{\mathcal{S}'}(x)\big)\big]H_{\mathcal{S}', k}(x)\,\pmb\bigiotab(x)\, d\mathcal{H}^n(x)  \, ,
			\end{align}}
			for every $ \phi \in C^\infty(\mathbf{R}^{n+1} \times \mathbf{R}^{n+1}) $ and every $k\in\{0,\ldots,n\}$. 
		\end{theorem}
		\begin{proof} To prove that $\mathcal{N}_{\mathcal{S}'}$ is a Legendrian cycle of $\mathbf{R}^{n+1}$, note that $\mathcal{N}_{\mathcal{S}'}\restrict(U_i\times\mathbf{R}^{n+1})$ is a Legendrian cycle of $U_i$ for every $i\in\{1,\ldots,N\}$: indeed, for every $\smash{\varphi\in\mathcal{D}^n(U_i\times\mathbf{R}^{n+1})}$, from (\ref{Ufset32}) and (\ref{norBasis4}) we obtain
			{\allowdisplaybreaks\begin{align*}
					\big[\mathcal{N}_{\mathcal{S}'}\restrict(U_i\times\mathbf{R}^{n+1})\big](\varphi)&=\int_{\textup{nor}(\mathcal{S}')\restrict U_i}\langle \vv{\xi}_{\mathcal{S}'},\varphi\rangle \, \pmb\bigiotab\circ\pi_0 \, d\mathcal{H}^n&\\
					&=\sum_{j=1}^{q(z_i)}\int_{\textup{nor}(\Gamma^{(i)}_j)\restrict U_i}\langle \vv{\xi}_{\mathcal{S}'},\varphi\rangle \, d\mathcal{H}^n&\\
					&=\sum_{j=1}^{q(z_i)}\int_{\textup{nor}(\Gamma^{(i)}_j)\restrict U_i}\langle \vv{\eta}_{\Gamma^{(i)}_j\cap U_i},\varphi\rangle \, d\mathcal{H}^n&\\
					&=\sum_{j=1}^{q(z_i)}\mathcal{N}_{\Gamma^{(i)}_j\cap U_i}(\varphi) \, ,
			\end{align*}}where every $\mathcal{N}_{\Gamma^{(i)}_j\cap U_i}\in\mathcal{D}_n(U_i\times\mathbf{R}^{n+1})$ is a Legendrian cycle of $U_i$ (cf.\,(\ref{LusinNorGraph}) in Remark \ref{EpiAndGraph}). Applying Lemma \ref{lem patching of legendrian cycles}, we infer that $\mathcal{N}_{\mathcal{S}'}$ is a Legendrian cycle of $\mathbf{R}^{n+1}$.
			
			To prove (\ref{W2n domains0.80}), we recall that (cf.\,(\ref{Ufset25}))
			\begin{equation}\label{Ufset26}\nonumber
				J^{\textup{nor}(\mathcal{S}')}_n\pi_0(y,u)=\zeta_\mathcal{S'}(y,u)>0 \quad\textup{for} \ \mathcal{H}^n\textup{-a.e.} \ (y,u)\in\textup{nor}(\mathcal{S}')^{(n)} \, .
			\end{equation}
			Applying Lemma \ref{FuFirst9}, Lemma \ref{Ufset14} \emph{(ii)} and \emph{(vii)}, Lemma \ref{Ufset15} \emph{(iii)}, (\ref{n-curv2}), and  the area formula for rectifiable sets \cite[Theorem 3.2.22 (3)]{Fed69}, and noting that 
			$\spt(\mathcal{N}_{\mathcal{S}'})$ is compact, we infer 
			{\allowdisplaybreaks\begin{align*}
					&\big(\mathcal{N}_{\mathcal{S}'} \restrict \varphi_{n-k}\big)(\phi)={n \choose k} \int_{\textup{nor}(\mathcal{S}')\restrict N_2(\mathcal{S'})}\phi(y,u)\,\pmb\bigiotab(y)\,H_{\textup{nor}(\mathcal{S}'),k}(y,u)\, J^{\textup{nor}(\mathcal{S}')}_n\pi_0(y,u)\, d\mathcal{H}^n(y,u)&\\
					&\quad\quad={n \choose k} \int_{N_2(\mathcal{S}')}\sum_{(y,u)\in(\pi_0|\textup{nor}(\mathcal{S}')\restrict N_2(\mathcal{S}'))^{-1}(x)}\big[\phi(y,u)\,H_{\textup{nor}(\mathcal{S}'),k}(y,u)\,\pmb\bigiotab(y)\big]\, d\mathcal{H}^n(x)&\\
					&\quad\quad={n \choose k}\int_{\mathcal{S}'} \big[\phi\big(x, \nu_{\mathcal{S}'}(x)\big) +(-1)^k\,\phi\big(x,-\nu_{\mathcal{S}'}(x)\big)\big]H_{\mathcal{S}', k}(x)\,\pmb\bigiotab(x)\, d\mathcal{H}^n(x) \ ,
			\end{align*}}for any $ \phi \in C^\infty(\mathbf{R}^{n+1} \times \mathbf{R}^{n+1}) $ and $ k\in\{0, \ldots, n \}$.
		\end{proof}
		
		Now we consider again the $C^1$-diffeomorphism 
		\begin{equation*}
			\Psi_F:(x,y)\in\mathbf{R}^{n+1}\times\mathbf{S}^{n}\mapsto \bigg(F(x), \frac{(D F(x)^{-1})^\ast(y)}{| (D F(x)^{-1})^\ast(y) |}\bigg)\in\mathbf{R}^{n+1}\times\mathbf{S}^n \, ,
		\end{equation*}
		for which (cf.\,\cite[Lemma 2.1]{SantilliRectVarCurv})
		\begin{equation}\label{PsiF4}
			\Psi_F\big(\textup{nor}(\mathcal{S}')\big) = \textup{nor}(\mathcal{S}) \, .
		\end{equation}
		
		\begin{definition} We define $\mathcal{N}_{\mathcal{S}}\in\mathcal{D}_n(\mathbf{R}^{n+1}\times\mathbf{R}^{n+1})$ by
			$$\mathcal{N}_{\mathcal{S}}:=(\Psi_F)_{\#}(\mathcal{N}_{\mathcal{S}'})\,.$$
		\end{definition}
		
		\begin{theorem}\label{ReillyWset1}
			$\mathcal{N}_{\mathcal{S}}$ is a Legendrian cycle of $\mathbf{R}^{n+1}$, we call it the Legendrian cycle associated with $\mathcal{S}$. In particular
			\begin{equation}\label{Ufset44}
				\mathcal{N}_{\mathcal{S}}=\big(\pmb\bigiotab\circ\pi_0\circ(\Psi_F|\textup{nor}(\mathcal{S}'))^{-1}\big)\,\big(\mathcal{H}^n\restrict\textup{nor}(\mathcal{S})\big)\wedge\vv{\xi}_{\mathcal{S}} \,.
			\end{equation}Moreover, for a selected unit-normal vector field $\nu_{\mathcal{S}}$ on $\mathcal{S}$, the following relations holds
			{\allowdisplaybreaks\begin{align}\label{relations}
					&\big(\mathcal{N}_{\mathcal{S}} \restrict \varphi_{n-k}\big)(\phi)&\\ \nonumber
					&\quad\quad\quad={n \choose k}\int_{\mathcal{S}} \big[\phi\big(x, \nu_{\mathcal{S}}(x)\big) +(-1)^k\,\phi\big(x,-\nu_{\mathcal{S}}(x)\big)\big]H_{\mathcal{S}, k}(x)\,\pmb\bigiotab\big(F^{-1}(x)\big)\, d\mathcal{H}^n(x)
			\end{align}}for every $ \phi \in C^\infty(\mathbf{R}^{n+1} \times \mathbf{R}^{n+1}) $ and $ k\in\{0, \ldots, n \}$.\end{theorem}
		\begin{proof} We introduce $\psi:=\Psi_F|\textup{nor}(\mathcal{S}')$. Recalling (\ref{PsiF4}) and since
			\begin{equation}\label{W2n domains eq9}
				\textup{ap}\,D\psi\big(\psi^{-1}(y,v)\big)=D\Psi_F\big(\psi^{-1}(y,v)\big)\big|\textup{Tan}^n\big(\mathcal{H}^n\restrict\textup{nor}(\mathcal{S}'),\psi^{-1}(y,v)\big)
			\end{equation}
			for $\mathcal{H}^n$-a.e. $(y,v)\in\textup{nor}(\mathcal{S})$, we define the simple $\big(\mathcal{H}^n\restrict\textup{nor}(\mathcal{S})\big)$-measurable $n$-vectorfield
			{\allowdisplaybreaks\begin{align*}
					\vv{\eta}(y,v):&=\frac{\big[\bigwedge\nolimits_n\textup{ap}\,D\psi\big(\psi^{-1}(y,v)\big)\big]\vv{\xi}_{\mathcal{S}'}\big(\psi^{-1}(y,v)\big)}{\Big|\big[\bigwedge\nolimits_n\textup{ap}\,D\psi\big(\psi^{-1}(y,v)\big)\big]\vv{\xi}_{\mathcal{S}'}\big(\psi^{-1}(y,v)\big)\Big|}&\\
					&\ =\frac{\big[\bigwedge\nolimits_n\textup{ap}\,D\psi\big(\psi^{-1}(y,v)\big)\big]\vv{\xi}_{\mathcal{S}'}\big(\psi^{-1}(y,v)\big)}{\emph{J}_n^{\,\textup{nor}(\mathcal{S}')}\psi\big(\psi^{-1}(y,v)\big)} \quad \textup{for} \ \mathcal{H}^n\textup{-a.e.} \ (y,v)\in\textup{nor}(\mathcal{S}) \, .
			\end{align*}}
			By the area formula for rectifiable currents (cf.\,\cite[$4.1.30$]{Fed69}\,or\,\cite[p.\,197]{KrantzParks}) we obtain
			$$(\Psi_F)_{\#}(\mathcal{N}_{\mathcal{S}'})=(\pmb\bigiotab\circ\pi_0\circ\psi^{-1})\big(\mathcal{H}^n\restrict\textup{nor}(\mathcal{S})\big)\wedge\vv{\eta} \, ,$$
			where $\smash{|\vv{\eta}(y,v)|=1}$ and $\smash{\textup{Tan}^n\big(\mathcal{H}^n\restrict\textup{nor}(\mathcal{S}),(y,v)\big)}$ is associated with $\smash{\vv{\eta}(y,v)}$ for $\mathcal{H}^n$-a.e. $(y,v)\in\textup{nor}(\mathcal{S})$. Moreover, $(\Psi_F)_{\#}
			(\mathcal{N}_{\mathcal{S}'})$ is a cycle, and by the shuffle formula 
			(cf.~\cite[$1.4.2$]{Fed69}) and Lemma~\ref{lem legendrian property of 
				the normal bundle} it is also Legendrian.
			
			Since $\tau_{1}(x,u)\wedge\cdots\wedge \tau_{n}(x,u)=(-1)^n\,{\pmb{\ast}u}$\, for $\mathcal{H}^n$-a.e. $(x,u)\in\textup{nor}(\mathcal{S}')$, we infer (cf. (\ref{n-curv2}))
			{\allowdisplaybreaks\begin{align*}\label{W2n domains eq6.6}\nonumber
					&\Big[\bigwedge\nolimits_n \pi_0 \Big]\Big(\vv{\eta}\big(\Psi_F(x,u)\big)\Big)&\\ \nonumber
					&=\frac{\prod_{i=1}^{n}\big(1 + \kappa_{\mathcal{S}', i}(x,u)^2\big)^{-\frac{1}{2}}}{\emph{J}_n^{\,\textup{nor}(\mathcal{S}')}\psi(x,u)} \, \Big[\bigwedge\nolimits_n DF(x) \Big]\big(\tau_{1}(x,u)\wedge\cdots\wedge \tau_{n}(x,u)\big) &\\
					&=(-1)^n \ \frac{\prod_{i=1}^{n}\big(1 + \kappa_{\mathcal{S}', i}(x,u)^2\big)^{-\frac{1}{2}}}{\emph{J}_n^{\,\textup{nor}(\mathcal{S}')}\psi(x,u)} \, \Big[\bigwedge\nolimits_n DF(x) \Big](\pmb{\ast} u) \quad\textup{for} \ \mathcal{H}^n\textup{-a.e.} \ (x,u) \in \textup{nor}(\mathcal{S}')\,.
			\end{align*}}Therefore, from (\ref{PsiF4}) and by \cite[Remark 5.8]{SantilliValentini}, either 
			$$\big\langle \Big[\bigwedge\nolimits_n \pi_0 \Big]\big(\vv{\eta}(y,v)\big)\wedge v,\pmb{e}'_1\wedge\cdots\wedge \pmb{e}'_{n+1}\big\rangle > 0  \quad \textup{for} \ \mathcal{H}^n\textup{-a.e.} \ (y,v)\in\textup{nor}(\mathcal{S})$$
			or 
			$$\big\langle \Big[\bigwedge\nolimits_n \pi_0 \Big]\big(\vv{\eta}(y,v)\big)\wedge v,\pmb{e}'_1\wedge\cdots\wedge \pmb{e}'_{n+1}\big\rangle < 0 \quad\textup{for} \ \mathcal{H}^n\textup{-a.e.} \ (y,v)\in\textup{nor}(\mathcal{S})\,.$$
			In addition, for $\mathcal{H}^n$-a.e. $(y,v)\in\textup{nor}(\mathcal{S})$, $\smash{\textup{Tan}^n\big(\mathcal{H}^n\restrict\textup{nor}(\mathcal{S}),(y,v)\big)}$ is associated with both the $n$-vectorfields $\smash{\vv{\eta}(y,v)}$ and $\vv{\xi}_{\mathcal{S}}(y,v)$, where $($cf.\,Definition $\ref{defn-vect}$ and Lemma $\ref{lem: Santilli20})$
			$$\big\langle \Big[\bigwedge\nolimits_n \pi_0 \Big]\big(\vv{\xi}_{\mathcal{S}}(y,v)\big)\wedge v,\pmb{e}'_1\wedge\cdots\wedge \pmb{e}'_{n+1}\big\rangle >0  \quad\textup{for} \ \mathcal{H}^n\textup{-a.e.} \ (y,v)\in\textup{nor}(\mathcal{S})\,.$$
			Hence, up to a change of sign,
			$$\vv{\eta}(y,v)=\vv{\xi}_{\mathcal{S}}(y,v) \quad \textup{for} \ \mathcal{H}^n\textup{-a.e.} \ (y,v)\in\textup{nor}(\mathcal{S})\,,$$
			hence $(\Psi_F)_{\#}(\mathcal{N}_{\mathcal{S}'}) = \mathcal{N}_{\mathcal{S}}$.
			
			To prove (\ref{relations}), we apply Lemma \ref{FuFirst9}, Lemma \ref{Ufset14} \emph{(ii)} and \emph{(vii)}, Lemma \ref{Ufset15} \emph{(iii)}, (\ref{n-curv2}), and  the area formula for rectifiable sets \cite[Theorem 3.2.22 (3)]{Fed69} to infer 
			{\allowdisplaybreaks\begin{align*}
					\big(\mathcal{N}_{\mathcal{S}} \restrict \varphi_{n-k}\big)(\phi) &={n \choose k} \int_{\textup{nor}(\mathcal{S})\restrict N_2(\mathcal{S})}\phi(y,u)\,\pmb\bigiotab\big(F^{-1}(y)\big)\,H_{\textup{nor}(\mathcal{S}),k}(y,u)\, J^{\textup{nor}(\mathcal{S})}_n\pi_0(y,u)\, d\mathcal{H}^n(y,u)&\\
					&={n \choose k} \int_{N_2(\mathcal{S})}\sum_{(y,u)\in(\pi_0|\textup{nor}(\mathcal{S})\restrict N_2(\mathcal{S}))^{-1}(x)}\!\!\!\big[\phi(y,u)\,H_{\textup{nor}(\mathcal{S}),k}(y,u)\,\pmb\bigiotab\big(F^{-1}(y)\big)\big]\, d\mathcal{H}^n(x)&\\
					&={n \choose k}\int_{\mathcal{S}} \big[\phi\big(x, \nu_{\mathcal{S}}(x)\big) +(-1)^k\,\phi\big(x,-\nu_{\mathcal{S}}(x)\big)\big]H_{\mathcal{S}, k}(x)\,\pmb\bigiotab\big(F^{-1}(x)\big)\, d\mathcal{H}^n(x) \, ,
			\end{align*}}for any $ \phi \in C^\infty(\mathbf{R}^{n+1} \times \mathbf{R}^{n+1}) $ and $ k\in\{0, \ldots, n \}$. The proof is complete.
		\end{proof}

		\begin{remark}\label{FSprop}\textup{Let $G$ be a $C^2$-diffeomorphism of $\mathbf{R}^{n+1}$. Then
				\begin{equation}\label{PushReilly}
					(\Psi_G)_{\#}(\mathcal{N}_{\mathcal{S}}) =(\Psi_{G\circ F})_{\#}(\mathcal{N}_{\mathcal{S}'})= \mathcal{N}_{G(\mathcal{S})}\,.
				\end{equation}
				Moreover, for a selected unit-normal vector field $\nu_{G(\mathcal{S})}$ on $G(\mathcal{S})$, from (\ref{relations}) we obtain
				{\allowdisplaybreaks\begin{align}\label{W2n domains0.801}\nonumber
						\mathcal{N}_{G(\mathcal{S})}(\varphi_{n-k})&= {n \choose k}\int_{G(\mathcal{S})} \big(1+(-1)^k\big)\,H_{G(\mathcal{S}),k}(x)\,\pmb\bigiotab\big((G\circ F)^{-1}(x)\big)\, d\mathcal{H}^n(x)&\\
						&=\begin{cases}
							2\, \displaystyle{n \choose k}\,\mathcal{A}_k(G(\mathcal{S}))& \textup{if $k$ is even}\\0& \textup{if $k$ is odd}
						\end{cases}
				\end{align}}for any $k\in\{0,\ldots,n\}$, where $H_{G(\mathcal{S}),k}$ is the $k$-th mean curvature of $G(\mathcal{S})$ with respect to $\nu_{G(\mathcal{S})}$.
		}\end{remark}

		Now we derive the following extension of Reilly's variational formulae to $\mathscr{W}^{2,n}$-sets.
		
		\begin{theorem}\label{ReillyW-sets}
			Let $ \{F_t\}_{t \in (-\epsilon, \epsilon)} $ be a local variation of $ \mathbf{R}^{n+1} $ with initial velocity vector field $ V $. If $k\in\{1, \ldots , n\}$ is odd, then
			$$\frac{d}{dt}\,\mathcal{A}_{k-1}\big(F_t(\mathcal{S})\big) \Big|_{t =0}=(n-k+1)\int_{\mathcal{S}}  V(x)\bullet\nu_\mathcal{S}(x) \,H_{\mathcal{S}, k}(x)\,\pmb\bigiotab\big(F^{-1}(x)\big)\, d\mathcal{H}^n(x)\,.$$Moreover, if $n$ is even,
			$$\frac{d}{dt}\,\mathcal{A}_n\big(F_t(\mathcal{S})\big) \Big|_{t =0}=0 \, .$$
		\end{theorem}
		
		\begin{proof}
			Combining (\ref{PushReilly}) and (\ref{W2n domains0.801}), we obtain
			\begin{align}\label{formula}\nonumber
				&\big[(\Psi_{F_t})_{\#}(\mathcal{N}_\mathcal{S})\big](\varphi_{n-k+1}) = \mathcal{N}_{F_t(\mathcal{S})}(\varphi_{n-k+1})&\\
				&\quad\quad=\begin{cases}
					2\, \displaystyle{n \choose k-1}\,\mathcal{A}_{k-1}\big(F_t(\mathcal{S})\big)&\textup{if $k\in\{1,\ldots,n+1\}$ is odd}\\
					0&\textup{if $k\in\{1,\ldots,n+1\}$ is even}\,.
				\end{cases}
			\end{align}
			From (\ref{formula}), if $k\in\{1, \ldots , n\}$ is odd, setting $\theta_V(x,y):=V(x)\bullet y$ for $x,y\in\mathbf{R}^{n+1}$ and applying Lemma \ref{Lemma Fu} and (\ref{relations}), we obtain
			{\allowdisplaybreaks\begin{align*}
					&2\,{n \choose k-1}\,\frac{d}{dt}\,\mathcal{A}_{k-1}\big(F_t(\mathcal{S})\big) \Big|_{t =0}&\\
					&\quad\quad=\frac{d}{dt}\big[(\Psi_{F_t})_{\#}(\mathcal{N}_\mathcal{S})\big](\varphi_{n-k+1})\Big|_{t=0}=k \ \mathcal{N}_\mathcal{S} (\theta_V\,\varphi_{n-k})&\\
					&\quad\quad=\big(1+(-1)^{k+1}\big)\,k\,{n \choose k}\int_{\mathcal{S}}  V(x)\bullet\nu_\mathcal{S}(x) \,H_{\mathcal{S}, k}(x)\,\pmb\bigiotab\big(F^{-1}(x)\big)\, d\mathcal{H}^n(x)&\\
					&\quad\quad=2\,(n-k+1){n \choose k-1}\int_{\mathcal{S}} V(x)\bullet\nu_\mathcal{S}(x) \,H_{\mathcal{S}, k}(x)\,\pmb\bigiotab\big(F^{-1}(x)\big)\, d\mathcal{H}^n(x)\,.
			\end{align*}}Moreover, if $n$ is even (i.e. $k=n+1$ odd), from Lemma \ref{Lemma Fu} we conclude
			$$  \frac{d}{dt}\,\mathcal{A}_n\big(F_t(\mathcal{S})\big) \Big|_{t =0}=\frac{d}{dt}\big[(\Psi_{F_t})_{\#}(\mathcal{N}_\mathcal{S})  \big](\varphi_{0}) \Big|_{t=0} = 0 \, .$$
			The proof is complete.
		\end{proof}

		\bibliography{valentini.bib}{}

\begin{thebibliography}{CCKS96}

\bibitem[AFP00]{AFP00}
Luigi Ambrosio, Nicola Fusco, and Diego Pallara.
\newblock {\em Functions of bounded variation and free discontinuity problems}.
\newblock Oxford Mathematical Monographs. The Clarendon Press, Oxford
  University Press, New York, 2000.

\bibitem[AGP98]{AmGobPal}
L.~Ambrosio, M.~Gobbino, and D.~Pallara.
\newblock Approximation problems for curvature varifolds.
\newblock {\em J. Geom. Anal.}, 8(1):1--19, 1998.

\bibitem[CCKS96]{Caffarelli96}
L.~Caffarelli, M.~G. Crandall, M.~Kocan, and A.~Swiech.
\newblock On viscosity solutions of fully nonlinear equations with measurable
  ingredients.
\newblock {\em Comm. Pure Appl. Math.}, 49(4):365--397, 1996.

\bibitem[CV77]{CastaingValadierBook}
C.~Castaing and M.~Valadier.
\newblock {\em Convex analysis and measurable multifunctions}, volume Vol. 580
  of {\em Lecture Notes in Mathematics}.
\newblock Springer-Verlag, Berlin-New York, 1977.

\bibitem[CZ61]{CalderonZygmund}
A.-P. Calder\'on and A.~Zygmund.
\newblock Local properties of solutions of elliptic partial differential
  equations.
\newblock {\em Studia Math.}, 20:171--225, 1961.

\bibitem[Fed59]{Fed59}
Herbert Federer.
\newblock Curvature measures.
\newblock {\em Trans. Amer. Math. Soc.}, 93:418--491, 1959.

\bibitem[Fed69]{Fed69}
Herbert Federer.
\newblock {\em Geometric measure theory}.
\newblock Die Grundlehren der mathematischen Wissenschaften, Band 153.
  Springer-Verlag New York, Inc., New York, 1969.

\bibitem[Fed78]{FedererNote}
Herbert Federer.
\newblock Colloquium lectures on geometric measure theory.
\newblock {\em Bull. Amer. Math. Soc.}, 84(3):291--338, 1978.

\bibitem[Fu98]{Fu98}
Joseph H.~G. Fu.
\newblock Some remarks on {L}egendrian rectifiable currents.
\newblock {\em Manuscripta Math.}, 97(2):175--187, 1998.

\bibitem[Fu11]{FuAlexandrov}
Joseph H.~G. Fu.
\newblock An extension of {A}lexandrov's theorem on second derivatives of
  convex functions.
\newblock {\em Adv. Math.}, 228(4):2258--2267, 2011.

\bibitem[GM04]{GiaquintaModica_book}
Mariano Giaquinta and Giuseppe Modica.
\newblock {\em Mathematical Analysis, Foundations and Advanced Techniques for
  Functions of Several Variables}, volume~25 of {\em Oxford Lecture Series in
  Mathematics and its Applications}.
\newblock Oxford University Press, Oxford, 2004.

\bibitem[GT01]{GilTru_Book}
David Gilbarg and Neil~S. Trudinger.
\newblock {\em Elliptic partial differential equations of second order}.
\newblock Classics in Mathematics. Springer-Verlag, Berlin, 2001.
\newblock Reprint of the 1998 edition.

\bibitem[Han11]{Qing}
Qing Han.
\newblock {\em A basic course in partial differential equations}, volume 120 of
  {\em Graduate Studies in Mathematics}.
\newblock American Mathematical Society, Providence, RI, 2011.

\bibitem[HS22]{HugSantilli}
Daniel Hug and Mario Santilli.
\newblock Curvature measures and soap bubbles beyond convexity.
\newblock {\em Adv. Math.}, 411(part A):Paper No. 108802, 89, 2022.

\bibitem[KP08]{KrantzParks}
Steven~G. Krantz and Harold~R. Parks.
\newblock {\em Geometric integration theory}.
\newblock Cornerstones. Birkh\"auser Boston, Inc., Boston, MA, 2008.

\bibitem[LL16]{Leonard}
I.~E. Leonard and J.~E. Lewis.
\newblock {\em Geometry of convex sets}.
\newblock John Wiley \& Sons, Inc., Hoboken, NJ, 2016.

\bibitem[Mat95]{MattilaBook}
Pertti Mattila.
\newblock {\em Geometry of sets and measures in {E}uclidean spaces}, volume~44
  of {\em Cambridge Studies in Advanced Mathematics}.
\newblock Cambridge University Press, Cambridge, 1995.
\newblock Fractals and rectifiability.

\bibitem[Men24]{menne2024sharplowerboundmean}
Ulrich Menne.
\newblock A sharp lower bound on the mean curvature integral with critical
  power for integral varifolds, 2024.

\bibitem[MM73]{MarcusMizel}
M.~Marcus and V.~J. Mizel.
\newblock Transformations by functions in {S}obolev spaces and lower
  semicontinuity for parametric variational problems.
\newblock {\em Bull. Amer. Math. Soc.}, 79:790--795, 1973.

\bibitem[MS19]{MenSan}
Ulrich Menne and Mario Santilli.
\newblock A geometric second-order-rectifiable stratification for closed
  subsets of {E}uclidean space.
\newblock {\em Ann. Sc. Norm. Super. Pisa Cl. Sci. (5)}, 19(3):1185--1198,
  2019.

\bibitem[Rei72]{Reilly1972}
Robert~C. Reilly.
\newblock Variational properties of mean curvatures.
\newblock In {\em Proceedings of the {T}hirteenth {B}iennial {S}eminar of the
  {C}anadian {M}athematical {C}ongress ({D}alhousie {U}niv., {H}alifax,
  {N}.{S}., 1971), {V}ol. 2}, pages 102--114. Canad. Math. Congr., Montreal,
  QC, 1972.

\bibitem[Ros16]{Roskovec}
Tom\'as Roskovec.
\newblock Construction of {$W^{2,n}(\Omega)$} function with gradient violating
  {L}usin ({N}) condition.
\newblock {\em Math. Nachr.}, 289(8-9):1100--1111, 2016.

\bibitem[RW98]{RockafellarWets}
R.~Tyrrell Rockafellar and Roger J.-B. Wets.
\newblock {\em Variational analysis}, volume 317 of {\em Grundlehren der
  mathematischen Wissenschaften [Fundamental Principles of Mathematical
  Sciences]}.
\newblock Springer-Verlag, Berlin, 1998.

\bibitem[RZ19]{RatajZaehlebook}
Jan Rataj and Martina Z\"{a}hle.
\newblock {\em Curvature measures of singular sets}.
\newblock Springer Monographs in Mathematics. Springer, Cham, 2019.

\bibitem[San20]{SantilliAnnali}
Mario Santilli.
\newblock Fine properties of the curvature of arbitrary closed sets.
\newblock {\em Ann. Mat. Pura Appl. (4)}, 199(4):1431--1456, 2020.

\bibitem[San21]{SantilliRectVarCurv}
Mario Santilli.
\newblock Second order rectifiability of varifolds of bounded mean curvature.
\newblock {\em Calc. Var. Partial Differential Equations}, 60(2):Paper No. 81,
  17, 2021.

\bibitem[San24]{Santilli24}
Mario Santilli.
\newblock {Finite Total Curvature and Soap Bubbles With Almost Constant
  Higher-Order Mean Curvature}.
\newblock {\em International Mathematics Research Notices}, page rnae159, 2024.

\bibitem[Sim83]{SimonBook}
Leon Simon.
\newblock {\em Lectures on geometric measure theory}, volume~3 of {\em
  Proceedings of the Centre for Mathematical Analysis, Australian National
  University}.
\newblock Australian National University, Centre for Mathematical Analysis,
  Canberra, 1983.

\bibitem[SV25]{SantilliValentini}
Mario Santilli and Paolo Valentini.
\newblock Alexandrov sphere theorems for {$W^{2,n}$}-hypersurfaces.
\newblock {\em Proceedings of the Royal Society of Edinburgh: Section A
  Mathematics}, pages 1--49, 2025.

\bibitem[Tor94]{Toro}
Tatiana Toro.
\newblock Surfaces with generalized second fundamental form in {$L^2$} are
  {L}ipschitz manifolds.
\newblock {\em J. Differential Geom.}, 39(1):65--101, 1994.

\bibitem[Tru89]{Trudinger89}
Neil~S. Trudinger.
\newblock On the twice differentiability of viscosity solutions of nonlinear
  elliptic equations.
\newblock {\em Bull. Austral. Math. Soc.}, 39(3):443--447, 1989.

\bibitem[Val25]{Val25}
Paolo Valentini.
\newblock Legendrian cycles and \uppercase{A}lexandrov sphere theorems for
  $\uppercase{W}^{2,n}$-hypersurfaces.
\newblock {\em Ph.D. Thesis in Mathematics, University of L'Aquila}, 2025.
\newblock \url{https://hdl.handle.net/11697/262249}.

\end{thebibliography}
		\bibliographystyle{alpha}

	\end{document}